\documentclass[12pt]{amsart}
\usepackage[margin=2cm]{geometry}
\usepackage{amssymb}

\usepackage{amsmath,amsfonts,amssymb,amsthm,amscd,latexsym,euscript,xypic,  verbatim}
\usepackage{amstext}
\usepackage{graphicx}
\usepackage[active]{srcltx}
\usepackage[all]{xy}
\usepackage{graphics}
\usepackage{color}
\usepackage[utf8]{inputenc}
\usepackage{mathtools}
\usepackage{tikz}
\usepackage{dsfont}
\usepackage{tikz-cd}
\usepackage{float}
\usepackage{multirow}
\usepackage{hyperref}

\tikzset{
  symbol/.style={
    draw=none,
    every to/.append style={
      edge node={node [sloped, allow upside down, auto=false]{$#1$}}}
  }
}

\usetikzlibrary{decorations.pathmorphing}
\swapnumbers
\theoremstyle{plain}
\newtheorem{prop}{Proposition}[section]
\newtheorem{thm}[prop]{Theorem}
\newtheorem{cor}[prop]{Corollary}
\newtheorem{lemma}[prop]{Lemma}

\theoremstyle{definition}
\newtheorem{point}[prop]{}

\newtheorem{Def}[prop]{Definition}

\newtheorem*{Def*}{Definition}

\newtheorem{example}[prop]{Example}
\newtheorem{examples}[prop]{Examples}

\newtheorem*{notation*}{Notation}
\newtheorem*{question*}{Question}

\DeclareMathOperator{\colim}{colim}

\xyoption{2cell}

\begin{document}

\title{Realisations of posets and tameness}

\author{Wojciech Chach\'olski}
\thanks{Mathematics, KTH, S-10044
Stockholm, Sweden\email{wojtek@kth.se} \email{alvinj@kth.se} \email{tombari@kth.se} }
\author{ Alvin Jin }
\author{Francesca Tombari}

\maketitle

\begin{abstract}
We introduce a construction called  realisation which transforms posets into posets.
We show that realisations share several key features with upper semilattices.
For example, we define  local dimensions of points  in a poset and show that  these numbers for realisations behave in a similar way as they do for  upper semilattices. Furthermore, similarly to upper semilattices, realisations
have well behaved discrete approximations which are suitable for capturing  homological properties of functors indexed by them.  These discretisations  are convenient  and  effective for describing tameness of functors. 
Homotopical and homological properties of tame functors,
particularly those indexed by  realisations, are  discussed, 
with emphasis on the use of Koszul complexes to compute Betti diagrams of minimal free resolutions of tame functors indexed by upper semilattices and realisations.
\end{abstract}

\section{Introduction}
Vector space valued functors indexed by posets, also called representations or persistent modules, are ubiquitous in mathematics. They are  commonly used to  
encode how various homological invariants change, for example, in algebraic topology, algebraic geometry, homotopy theory, and more recently in topological data analysis (TDA).
To obtain such a representation one often starts with a functor $\alpha\colon [0,\infty)^r\to 2^X$, where $X$ is a topological space. The functor $\alpha$ parameterises a filtration of the topological space whose homologies give rise to a representation. 
For example, the homologies of the sublevel sets of a function $X\to [0,\infty)^r$ form a representation/persistence module of the poset $[0,\infty)^r$, and so do the homologies of Vietoris-Rips complex (multi) filtrations (see~\cite{MR2476414}). Encoding information in the form of $[0,\infty)^r$ representations is attractive for three reasons:
\begin{itemize}
    \item  metric properties  of $[0,\infty)^r$  can be used to define and study  distances on
functors indexed by $[0,\infty)^r$ (see for example~\cite{algebraic_stability, Lesnick2015}), which are essential for addressing stability of various invariants and can be used for hierarchical stabilisation constructions (see~\cite{ MR4057607, OliverWojtek, MR3735858});
\item  the poset   $[0,\infty)^r$ has well behaved discrete  approximations given by  sublattices of the form ${\mathbb N}^r\hookrightarrow [0,\infty)^r$, which  can be used to provide finite  approximations  of  $[0,\infty)^r$ representaions;
\item  the mentioned discretisations and approximations  have well \allowbreak studied algebraic and homological properties, as the path algebra of the poset ${\mathbb N}^r$ is 
isomorphic to the multigraded polynomial ring in $r$ variables.
\end{itemize}
Generalising these results to functors indexed by other posets 
is a growing research direction  in the applied topology 
community, reflected by an increasing  number of publications on this subject, see for example~\cite{Botnan2020ART,MR4323617,MR4334502, EMiller1, MR3975559}.  However there seems to be lack of explicit examples of posets, not  directly related to $[0,\infty)^r$, for which the mentioned three aspects are tightly intertwined. 
The aim of this article
 is to introduce a family of such posets.
 We define a construction called \textbf{realisation}
(see~\ref{drjfhhj}) which transforms   posets into  posets, and  show that  the realisations of 
finite-type posets (posets where down sets of  elements are finite) have natural discretisations
tightly related to homological algebra properties of functors indexed by them. 
Metric aspects of functors indexed by realisations are going to be the subject of a follow up paper. 

\begin{figure}
    \centering
    \includegraphics[width=10cm]{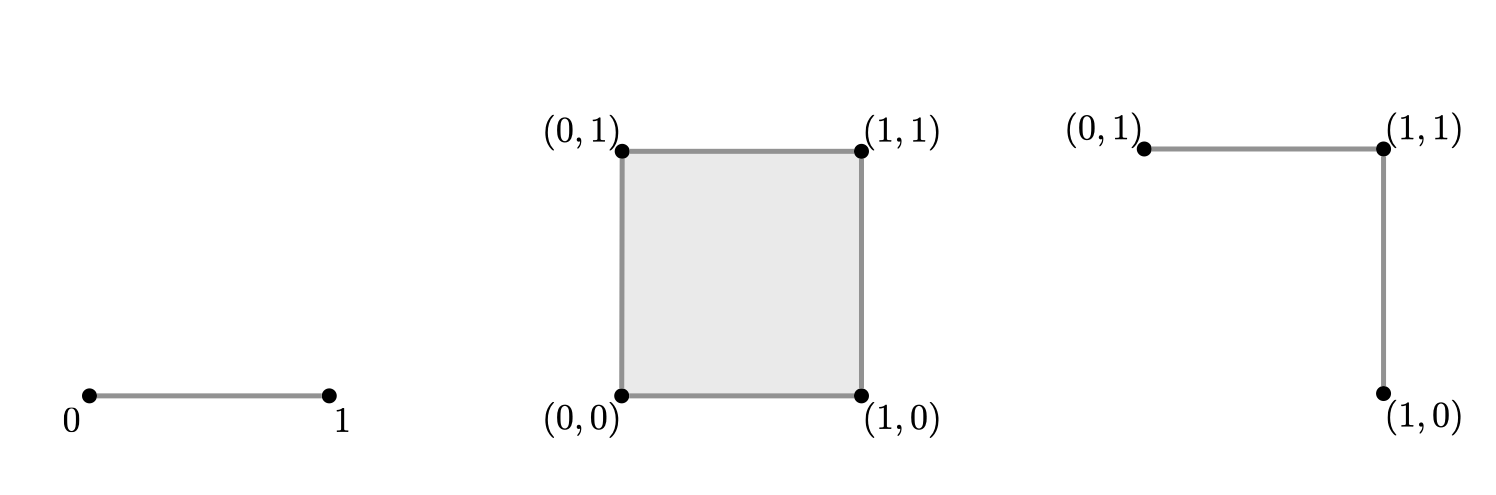}
    \caption{On the left, the realisation of the finite poset $\{0<1\}$, in the center, the one of $\{0<1\}^2$ and, on the right, the one of $\{(1,0), (0,1)<(1,1)\}$. The different tones of gray represent the different dimensions of cube assembled together: $0$-dimensional in black, $1$-dimensional in gray, and $2$-dimensional in lightgray.}
    \label{intro}
\end{figure}

\bigskip
The realisation of a poset $I$ is assembled by posets
of the form $(-1,0)^s$ in the following way. 
For every element $a$ in $I$, and for every finite subset $S$ of parents of $a$ (a parent is an element covered by $a$) that has a common ancestor, consider the poset
$(-1,0)^{|S|}$ with the usual product order. The realisation $\mathcal{R}(I)$ is the disjoint union  of
all these posets  $(-1,0)^{|S|}$ for all $a$ in $I$ and all $S$.
For example, if $S$ is empty, then $(-1,0)^0$ is of size $1$ and we  identify its element with $a$. In this way
we obtain an inclusion $I\subset \mathcal{R}(I)$. 
If $S=\{p\}$ is of size $1$, then we think about the associated subposet  $(-1,0)\subset \mathcal{R}(I) $ as the time to go back from $a$ to its parent $p$.
Figure~\ref{intro} illustrates  the realisations
of the posets $[1]:=\{0<1\}$, $[1]^2$ and the subposest of $[1]$ with elements $\{(1,0),(0,1),(1,1)\}$. 
The greyscale colors indicate some of the  summands $(-1,0)^s$, identifying these posets with, respectively, $[0,1]$, $[0,1]^2$ and the assembly of two copies of $[0,1]$.
We can, hence, see the realisation construction as an assembly of cubes of different dimensions associated with each element of the poset that respects the poset relation. 
An intuitive description of the realisation can be obtained if $I$ is a finite-type distributive (or more generally 
 consistent, see~\ref{adgdgfjdghj})  upper semilattice. 
 For example,  $[0,\infty)^r$ is the realisation of
 $\mathbb{N}^r$, both with the product order. 
For a general poset, however, to define the poset relation on its realisation  is non-trivial. Our strategy is to identify the realisation as a subposet of the Grothendieck construction of a certain  lax functor, which may fail to be a functor in the case $I$ is not a consistent upper semilattice. 
The fact that for an arbitrary  poset this functor is essentially lax is the main  difficulty in expressing the poset relation of the realisation in simple (i.e. functorial) terms.

Our story about  realisations  is divided into two parts.
In Part I, the combinatorial properties of  realisations as posets are discussed.
This part starts with the description of a certain regularity property of realisations of finite-type posets.
This regularity is expressed in terms of dimension and parental dimension, which are two numbers associated with each element of a poset (see Section~\ref{dsfsefseff})
capturing the complexity 
of expressing it as an upper bound of its downset.
Although in general these numbers might differ, they are the same for elements of a distributive
upper semilattice, which we regard as a regular poset.
It turns out that these numbers coincide also for every element in the realisation of a finite-type poset (Corollary~\ref{sdfwefwefwsghj}), that is another reason why we consider the realisation as a regularisation. 
In the first part, we also discuss assumptions under which the realisation of an upper semilattice is an upper semilattice. We show that this is guaranteed under a consistency assumption (see~\ref{adgdgfjdghj}). Key examples of consistent  upper semilattices are distributive upper semilattices. In this context, our main result is Theorem~\ref{dasghfgjk} stating that the realisation of a finite-type  distributive  upper semilattice is a distributive  upper semilattice.

In Part II, we focus on tame functors.
The notion of tameness of $[0,\infty)^r$ representations has been central in both geometry and applied topology, as finiteness properties of such functors guarantee that various invariants can be defined and calculated.
A tame functor is indexed by possibly an infinite poset, but it is only finitely encoded, in the sense that its behaviour is only governed by a finite indexing poset. 
This is why we also refer to them as discretisable functors.
The notion of tameness is, in fact, restrictive and meaningful only when $J$ is an infinite poset, for example if it is a realisation. 
One reason we care about upper semilattices in the first part is because finite-type functions out of them, admit a transfer (see Section~\ref{asfsdfhfhgjmd}). 
This is fundamental to show that the category of tame functors indexed by upper semilattices with values in a category closed under finite limits and colimits is also closed under finite limits and colimits (see Section~\ref{aoceavojcnigt} and, in particular, Proposition~\ref{afsadagdfhbg}).
Another important property of tame functors indexed by a realisation $\mathcal{R}(I)$ is that they can always be discretised by subposet inclusions $\mathcal{R}(I, V)\subset \mathcal{R}(I)$, where $\mathcal{R}(I, V)$ is a finite poset preserving all the properties of $\mathcal{R}(I)$, such as the dimension of the elements and being an upper semilattice.
Moreover, we can show that tame functors  indexed by $\mathcal{R}(I)$, when
$I$ is a finite consistent upper semilattice,
are exactly the functors which are constant on the fibers of the transfer 
$f^!\colon \mathcal{R}(I)\to \mathcal{R}(I,V)_\ast$ (see~\ref{wertgf}).
 This is a standard way of describing tame  functors indexed by $[0,n)^r$ (see~\cite{MR3873177, MR3975559, MR3735858}).

In Section~\ref{afasdfdghfdgjh} we give another motivation for the use of upper semilattices to index tame functors. 
In particular, 
we prove (see Theorem~\ref{afadfhsfgh}) that a model structure on the category $\mathcal{C}$ naturally extends to a model structure on $\mathcal{C}$ valued tame functors indexed by an  upper semilattice.
This enables the use of homotopical and homological algebra tools to study tame functors indexed by  upper semilattices. 
There are situations however when we would like to understand homological properties of tame functors  indexed by realisations which are not upper semilattices. Surprisingly, a lot can be described in such cases as well.
Section~\ref{sdgadfhsgdb} is devoted to present how one might construct minimal resolutions and calculate Betti diagrams using Koszul complexes for  vector space valued tame functors indexed by rather general posets. 
These results become particularly transparent in the case of tame functors indexed by realisations which are described in Theorem~\ref{sdrtyhgf}. 
For example, it is noticeable that the parental dimension of such tame functors can be used to bound their homological dimensions, bringing together the combinatorial study of realisations done in the first part of the article and their use as indexing posets of tame functors. 

\bigskip
\noindent
\textbf{Part I, posets.}

\section{Basic notions}\label{adsgehg}
 In this section, we recall and introduce some basic concepts that might already be familiar. The reader may choose to skip this section and come back to it 
 when a specific reference is made. 
\begin{point}\label{afadfhfgsh}
The standard poset of real  numbers is denoted by $\mathbb{R}$.
Its subposet of  non-negative real numbers is denoted by $[0,\infty)$, of  
natural numbers $\{0,1,\ldots\}$ by $\mathbb{N}$,  the first $n+1$ natural 
numbers by $[n]$, $[a,b]:=\{x\in\mathbb{R} \ |\  a\leq x\leq b\}$, 
and $(a,b]:=\{x\in\mathbb{R} \ |\  a< x\leq b\}$.

The product of posets  $(I,\leq )$ and $(J,\leq)$ is the poset $(I\times J, \leq)$, where   $(x,y)\leq (x_1,y_1)$   if  $x\leq  x_1$ and  $y\leq  y_1$.  
For a set $S$ and a poset $I$, the symbol $I^S$ denotes the poset of all functions $f\colon S\to I$ with  $f\leq g$ if $f(x)\leq g(x)$ for every $x$ in $S$. The poset $I^S$ is isomorphic to the $|S|$-fold product $I^{|S|}$. 

The poset  $[1]^S$ is called   the \textbf{discrete cube} of dimension $|S|$
and the poset $[-1,0]^S$ is called   the \textbf{geometric cube} of dimension $|S|$.

The inclusion poset of all subsets of $S$ is denoted by $2^S$.  The function mapping  $f\colon S\to [1]$  in $[1]^S$ to the subset $f^{-1}(1)\subset S$ is an isomorphism between  the posets $[1]^S$  and  $2^S$. This function is used to identify these two posets, and  $2^S$ is also referred to as the discrete cube of dimension $|S|$.

For a set $S$, consider a disjoint union  $S\coprod\{-\infty,\infty\}$ and the following relation on this enlarged set:  elements in $S$ are incomparable and $-\infty\leq s\leq  \infty$ for all  $s$ in $S$.
This poset is  denoted by $\Sigma S$ and called the \textbf{suspension} of $S$.
This poset is  denoted by $\Sigma S$ and called the \textbf{suspension} of $S$.
\end{point}

\begin{point}\label{dsdgfsfghjdghj}
Let $I$ be a poset and $a$ be its element.
For a subset $S\subset I$, define
$S\leq a:=\{s\in S\  |\  s\leq a\}$, $S< a:=\{s\in S\  |\  s< a\}$, and $a\leq  S\ : =\{s\in S\  |\  a\leq  s\}$.

A poset  $I$ is  called   of \textbf{finite-type} if  $I\leq a$ is finite for every $a$ in $I$.
A finite poset is of finite-type. The  poset  $\mathbb{N}$  is  infinite  and of  finite-type.  The poset  $\mathbb{R}$ is not of finite-type.
An element $x$ in $I$ is called a  \textbf{parent} of $a$ if   $x<a$  and 
  there is no element $y$ in $I$ such that $x<y<a$.
 The term $a$ \textbf{covers} $x$ is also commonly used 
 to say that $x$ is a parent of $a$.
    The symbol ${\mathcal P}_I(a)$ or  ${\mathcal P}(a)$, if $I$ is understood,  denotes the set of all parents of $a$.  
 For example, every element in $\mathbb{R}$  has an empty set of parents.   In the poset $[1]$,
  ${\mathcal P}_{[1]}(0)=\emptyset$  and  ${\mathcal P}_{[1]}(1)=\{0\}$. A poset $I$ is of finite-type if and only if, for every $a$ in $I$,
  the sets  ${\mathcal P}_I(a)$  and $\{n\  |\  \text{there is a sequence } x_0<\cdots < x_n =a\text{ in $I$}\}$ are finite.

Let $I$ and $J$ be posets and $f\colon I\to J$ be a function (not necessarily preserving the poset relations). 
For  an element  $a$ in $J$, define
\[f \leq  a:=\{x\in I\ |\ f(x)\leq a\}\ \ \ \ \ \ \ a\leq f :=\{x\in I\ |\ a\leq f(x)\}.\]
For example, $(\text{id}_I\leq a)\subset I$ coincides with 
 $(I\leq a)\subset I$.

 If $f\leq a$ is finite for every $a$ in $J$,  then  $f\colon I\to J$ is called of \textbf{finite-type}. For example, if $I$ is finite, then every function $f\colon I\to J$ is of finite-type.  
 \end{point}

\begin{point}\label{dsdgfsgjhdgh}
Let $S\subset  I$ be a subset of a poset $I$.
 If  $(a\leq S)=S$ (i.e. $a\leq s$ for all $s$ in $S$), then $a$ is called an \textbf{ancestor}, or lower bound, of $S$.  If  $(a< S)=S$ (i.e. $a< s$ for all $s$ in $S$), then $a$ is called a \textbf{proper ancestor} of $S$.  
Every element is a proper ancestor of the empty subset.

 If $a$ is an ancestor of $S$ for which  there is no other ancestor $b$ of $S$ such that $a<b$, then $a$ is called an \textbf{inf}, or greatest lower bound,  of $S$.
For example, an element $a$ in $I$ is an inf of the empty subset if and only if  
 there is no $b$ in $I$ for which $a<b$. Such an element is also called \textbf{maximal} in $I$. 
 For instance, an element  is an inf of $S$ if and only if it is  maximal in 
 the subposet $\{a\in I\ |\ (a\leq S)=S\}\subset I$.
 An element $a$ in $I$ is an inf of the entire $I$ if and only if  $a\leq x$ for every $x$ in $I$. If  such an element 
 exists, then it is unique and  is  called the \textbf{global minimum} of $I$.
 Two different inf elements of a subset are not comparable.
 
 In general, $S$ can have many inf  elements.  In the case  $S$ has only one inf  element, then this element is called the \textbf{product}, or meet, of $S$ and is denoted either by $\bigwedge_I S$ or   $\bigwedge S$, if $I$ is clear form the context.   Explicitly, the product of $S$ is an element $\nu$ in $I$  such that $(\nu \leq S)=S$ and, for every $a$ in $I$ for which $(a\leq S)=S$,  the relation $a\leq \nu$  holds. The product $\bigwedge_I\{x,y\}$ is also denoted as $x\wedge_I y$ or $x\wedge y$.  For example, if  the product of the empty subset of $I$ exists, then $I$ has a unique maximal element given by  $\bigwedge_I\emptyset $. This element may fail however to be the global maximum of $I$. 
 The element  $\bigwedge_I I $, if it exists, is the global minimum of $I$.
 \end{point}
 
\begin{point}\label{dsdgdfhfgjh}
 Let $S\subset  I$ be a subset of a poset $I$.
 If  $(S\leq a)=S$ (i.e. $s\leq a$ for all $s$ in $S$), then $a$ is called a \textbf{descendent}, or upper bound, of $S$.
Every element is a descendent of the empty subset. 

If $a$ is a descendent of $S$ for which  there is no other descendent $b$ of $S$ such that $b<a$, then $a$ is called a \textbf{sup}, or lowest upper bound, of $S$. 
For example, an element $a$ in $I$ is a sup of the empty subset if and only if there is no $b$ in $I$ for which $b<a$. 
Such an element is also called \textbf{minimal} in $I$. 
For instance, an element  is a sup of $S$ if and only if it is  minimal in the subposet
 $\{a\in I\ |\ (S\leq a)=S\}\subset I$.
 An element $a$ in $I$ is a sup of the entire $I$ if and only if  $x\leq a$ for every $x$ in $I$. If  such an element 
 exists, then it is unique and  is  called the \textbf{global maximum} of $I$.
 Two different sup elements of a subset are not comparable.

 In general, $S$ can have many sup  elements.  In the case  $S$ has only one sup element, then this element is called the \textbf{coproduct}, or join, of $S$ and is denoted either by $\bigvee_I S$ or
 $\bigvee S$, if $I$ is clear form the context.
 Explicitly, the coproduct of $S$ is an element $\nu$ in $I$  such that $(S\leq\nu)=S$, and for every $a$ in $I$ for which $(S\leq a)=S$,  the relation $ \nu\leq a$  holds. The coproduct $\bigvee_I\{x,y\}$ is also denoted as $x\vee_I y$ or $x\vee y$. For example, if 
 the coproduct of the empty subset of $I$ exists, then $I$ has a unique minimal element given by $\bigvee_I\emptyset$. This element may fail however to be the global minimum of $I$. 
  The element $\bigvee_I I$, if it exists, is the global maximum of $I$. 
 
 A product can be expressed as a coproduct: $\bigwedge_I  S = \bigvee_I \left(\bigcap _{x\in S}(I\le x)\right)$,  where the equality should be read as follows: the coproduct on the 
 right   exists if and only if the product on the left  exists,  in which case they are equal.
 Thus, if every  subset of $I$ has the  coproduct, then every subset  has also the product.
  \end{point}
  
 \begin{point}\label{sdfgbsdhgfsh}
  All subsets of $[0,\infty)$ have  products. All subsets of $[-1,0]$ have  products and  coproducts. 
  
  Let $S$ be a set.
  If all subsets of a poset $I$ have  products, then the same is true for all subsets of
  $I^S$, where the product of $T\subset I^S$ is given by the function
  mapping $x$ in $S$ to the product $\bigwedge_{I}\{f(x)\ |\ f\in T\}$ in $I$.
  If all subsets of $I$ have  coproducts, then the same is true for all subsets of
  $I^S$, where the coproduct of $T\subset I^S$ is given by the function
  $x\mapsto \bigvee_{I}\{f(x)\ |\ f\in T\}$.

 In the suspension  $\Sigma S$ (see~\ref{afadfhfgsh}), the element 
 $-\infty$ is its global minimum, and the element $\infty$ is its global maximum.
 The elements $-\infty$ and $\infty$ are respectively the product and coproduct of any subset  $U\subset S$ of size at least $2$. 
 The suspension $\Sigma S$ is of finite-type if and only if $S$ is finite.

 Every subset $Z$ of the discrete cube  $2^S$ (see~\ref{afadfhfgsh})  has  coproduct given by the union$\bigvee Z=\bigcup_{\sigma\in Z}\sigma$, and product 
 given by the intersection $\bigwedge Z=\bigcap_{\sigma\in Z}\sigma$. The subset $\emptyset\subset S$  is the global minimum of $2^S$
 and $S\subset S$ is the global maximum.
Note also that every parent of $U\subset S$ in $2^S$
is of the form $U\setminus \{x\}$ for $x$ in $U$. Thus,
the function $x\mapsto U\setminus \{x\}$ is a bijection between
$U$ and $\mathcal{P}_{2^S}(U)$ (see~\ref{dsdgfsfghjdghj}).
  \end{point}

   \begin{point}\label{dasgfsdfhhfk}
 A poset that has a global minimum is called \textbf{unital}. For example $[-1,0]$ is unital and
 $\mathbb{R}$ is not. If $I$ is unital, then so is $I^S$, for every set $S$, with its  global minimum  given by  the constant function mapping every $s$  in $S$ to the global minimum of $I$.
 A function $f\colon I\to J$ between unital posets is called  \textbf{unital} if it maps the global minimum in $I$ to the global minimum in $J$.
 
 For a poset $I$, the symbol $I_{\ast}$  denotes the  poset formed by adding an additional   element $-\infty$ to  $I$  and setting
$-\infty< x$ for all   $x$ in $I$. The element  $-\infty$  in $I_\ast$ is its  global  minimum.  

A function  $f\colon I\to J$ extends uniquely
to a unital function $f_\ast\colon I_\ast\to J_\ast$ 
for which  the  diagram
\[\begin{tikzcd}[row sep=small]
I\ar{r}{f} \ar[hook]{d} & J\ar[hook]{d}\\
I_\ast\ar{r}{f_\ast} & J_\ast
\end{tikzcd}\]
commutes.
Note that $(f_\ast\leq -\infty)=\{-\infty\}$, and  $(f_\ast\leq a)  = (f \leq a )\cup \{-\infty\}$ if $a$ is in $J$.
Thus, a function $f\colon I\to J$  is of finite-type if and only if $f_\ast\colon I_\ast\to J_\ast$ is of finite-type. Furthermore, $f_\ast\leq a$ is non-empty for every $a$ in $J_\ast$.
Moreover, the function $J^I\to J_{\ast}^{I_{\ast}}$, $f\mapsto f_\ast$,
is   injective and preserves the poset relations: if $f\leq g$, then $f_\ast \leq g_\ast$.
 \end{point}

\begin{point}\label{sdgsdgjdghkfhjk}
 For a poset  $(I,\leq)$, the same symbol $I$ denotes also  the  category  whose set of
objects is  $I$ and  where  
 $\text{mor}_I(x, y)$  is either empty, if $x\not\leq y$, or has cardinality $1$, if $x\leq y$.  
 
 A \textbf{functor} between poset categories  $I$ and  $J$ is  a function $f\colon I\to J$  preserving the poset relations:
 if $x\leq y$ in $I$, then $f(x)\leq f(y)$ in $J$.  
 The category whose objects are posets and morphisms are functors is denoted by  $\text{Posets}$.
 
 A function $f\colon I\to J$ is called a \textbf{subposet inclusion} if, 
 $x\leq y$ happens in $I$  if and only if $f(x)\leq f(y)$ in $J$. 
A subposet inclusion is always an injective functor. Not all injective functors however are subposet inclusions. 
A subposet inclusion 
 $f\colon I\to J$ induces an isomorphism between $I$ and the  subposet $f(I)\subset J$.
  
 The symbol $\text{Fun}(I,J)$  denotes the subposet of $J^{I}$ (see~\ref{afadfhfgsh})
 whose elements are functors. Since in this poset $f\leq g$ if $f(x)\leq g(x)$ for all $x$ in $I$, the category  $\text{Fun}(I,J)$ coincides  with the category whose morphisms are natural transformations. 
 More generally, for a poset $I$ and a category $\mathcal{C}$, the symbol $\text{Fun}(I,\mathcal{C})$ denotes the category of functors $I\to\mathcal{C}$ with natural transformations as morphisms. 
If $f\colon I\to J$ is a functor between posets, then the precomposition with $f$ functor is denoted by
$(-)^f\colon \text{Fun}(J,\mathcal{C})\to \text{Fun}(I,\mathcal{C})$.

The global minimum (see~\ref{dsdgfsgjhdgh}) of a poset coincides with the initial object in the associated  category.
  If $I$ is  a unital  poset (it has the global minimum) and  $\mathcal{C}$ is a category with an initial object, then a functor  $f\colon I\to \mathcal{C}$ is called  \textbf{unital} if it maps the global minimum  in $I$ to an initial object in $ \mathcal{C}$.
  The symbol $\text{Fun}_\ast(I, \mathcal{C})$  denotes the category whose objects  are  unital functors from   $I$ to  $\mathcal{C}$ and morphisms are natural transformations.

Let  $I$ be a poset and  $\mathcal{C}$ be a category with a unique initial object.  Any functor $F\colon I\to \mathcal{C}$ can be extended uniquely to a unital functor 
$F_\ast\colon I_\ast\to \mathcal{C}$  whose composition with $I\subset I_\ast$ is $F$. 
The association $F\mapsto F_\ast$ is an isomorphism of categories between $\text{Fun}(I, \mathcal{C})$ and $\text{Fun}_{\ast}(I_\ast, \mathcal{C})$.
 We use this isomorphism  to identify these categories.
  \end{point}

\begin{point}\label{asdfsgadfsfhjejk}
Let $I$ be a poset. A \textbf{lax functor} $T\colon I\rightsquigarrow \text{Posets}$
(\cite[Definition 1]{MR347936})
assigns to every element $a$ in $I$ a poset $T_a$ and to every relation $a\leq b$ in $I$
a functor $T_{a\leq b}\colon T_a\to T_b $. These functors are required to satisfy two   conditions. First, $T_{a\leq a}$ is the identity for all $a$ in $I$. Second,
$T_{a\leq c}\leq T_{b\leq c}T_{a\leq b}$ 
 for all $a\leq b\leq c$ in $I$.  For example, every functor  is a lax functor,
as in this case the equality $T_{a\leq c}= T_{b\leq c}T_{a\leq b}$ holds.

Let  $T\colon I\rightsquigarrow \text{Posets}$ be a lax functor.
 Define $\text{Gr}_I T : =\{(a,f)\  |\  a\in I,  f\in T_a\}$. For  $(a,f)$ and $(b,g)$ in $\text{Gr}_IT $, set $(a,f)\leq (b,g)$ if $a\leq b$ in $I$ and 
$T_{a\leq b}(f)\leq g$ in $T_b$. For example, $(a,f)\leq (a,g)$ if and only if
$f\leq g$ in $T_a$. The conditions required to be satisfied by lax functors guarantee  the transitivity of the relation $\leq$ on  
$\text{Gr}_IT$.
The poset $(\text{Gr}_I T,\leq)$
is called \textbf{Grothendieck construction}. In the case $T$ is the constant functor with value $J$, then $\text{Gr}_I T$  is isomorphic to the product $I\times J$. 

The function $\pi\colon \text{Gr}_IT\to I$, mapping $(a,f)$ to $a$, is a functor called the \textbf{standard projection}. For $a$ in $I$, the function $\text{in}_a\colon T_a\to \text{Gr}_IT$, mapping $f$ to $(a,f)$, is also a functor called 
the \textbf{standard inclusion}.  The functor $\text{in}_a$ is a subposet  inclusion (see~\ref{sdgsdgjdghkfhjk}).

\end{point}
\begin{prop}\label{adfgsgdhjdg} Let $I$ be a poset,   $T\colon I\rightsquigarrow \text{\rm Posets}$  a lax functor, and  $\emptyset\not=S\subset \text{\rm Gr}_IT$.
\begin{enumerate}
    \item Assume  $\pi(S)=\{a\}$. Then an element $(b,g)$ is a sup of
    $S$ in $\text{\rm Gr}_IT$ if and only if $b=a$ and $g$ is a  sup of $\{f\ |\   (a,f)\in S\}$ in $T_a$.
    \item Let $b$ be a sup  of $\pi(S)$ in $I$. Then an element $(b,m)$ is a sup of $S$
    in $\text{\rm Gr}_IT$ if and only if $m$ is a sup of 
    $\{T_{a\leq b}f\ |\ (a,f) \text{ in } S\}$ in $T_b$.
\end{enumerate}
\end{prop}

\begin{proof}
Statement (1) is  a particular case of (2).
Statement (2) is  a consequence of the equivalence:  $T_{a\leq b}f\leq g\leq m$ in $T_b$, for every $(a,f)$ in $S$, if and only if $(a,f)\leq (b,g)\leq (b,m)$ in $\text{\rm Gr}_IT$,
for every $(a,f)$ in $S$.
\end{proof}

\begin{point}\label{dadgdfhg}
Let $f\colon I\to J$ be a function (not necessarily a functor). Define 
$J[f]\subset 2^I$ to be  the subposet whose elements are 
subsets of the form $(f\leq a)\subset I$, for $a$ in $J$.
Let  $f\!\!\leq\colon J\to J[f]$ be  the function mapping $a$ to $f\leq a$. Note that if $a\leq b$ in $J$, then $(f\leq a)\subset (f\leq b)$,
which means $f\!\!\leq$ is a functor. An important property of being a functor is that its fibers $(f\!\!\leq)^{-1}(f\leq a)\subset J$, for $a$ in $J$, satisfy the following property: if $x\leq y$ in $J$ belong to $(f\!\!\leq)^{-1}(f\leq a)$, then so does any $z$ in $J$ such that $x\leq z\leq y$. Recall that such subposets of $J$ are called \textbf{intervals} or \textbf{convex}.
\end{point}

\begin{point}\label{werqergdfghjg}
Let $f\colon I\to J$ be a functor of posets and $\mathcal{C}$  a category. 
Recall that a functor $F\colon J\to \mathcal{C}$ is called a \textbf{left Kan extension} of $G\colon I\to \mathcal{C}$ along $f$  (see~\cite{MR1712872}), if there is a natural transformation 
 $\alpha\colon G\to Ff$ satisfying the following universal property: for every functor $H\colon J\to \mathcal{C}$, 
the function $\text{Nat}_{J}(F,H)\to \text{Nat}_I(G,Hf)$,
 mapping $\phi\colon F\to H$ to  the following composition, called the \textbf{adjoint} to $\phi$, is a bijection:
\[
\begin{tikzcd}
G\ar{r}{\alpha} & Ff\ar{r}{\phi^f} & 
Hf
\end{tikzcd}
\]
 This universal property  has two consequences. 
 The first is the  uniqueness of the left Kan extension: if $F$  and $F'$ are two left Kan extensions of $G$ along $f\colon I\to J$, 
and  $\alpha\colon G\to Ff$ and  $\alpha'\colon G\to F'f$ are  natural transformations satisfying the above universal property, then 
 there is a unique isomorphism $\psi\colon F\to F'$ for which
$\alpha' =\psi^f \alpha$. Because of this uniqueness, if it exists, a left Kan  extension
of $G\colon I\to \mathcal{C}$  along $f\colon I\to J$ is denoted by $f^kG\colon J\to \mathcal{C}$. 
Functoriality of left Kan extensions is another consequence of the universal property. 
Let $\phi\colon G\to G'$ be a natural transformation between two functors $G,G'\colon I\to\mathcal{C}$.
If these functors admit left Kan extensions along $f$, then there is a unique
natural transformation $f^k\phi\colon f^kG\to f^kG'$ for which the  following square commutes:
\[\begin{tikzcd}[row sep = 15pt]
G\ar{d}[swap]{\alpha}\ar{rr}{\phi} && G'\ar{d}{\alpha'}\\
(f^kG)f \ar{rr}{(f^k\phi)^f}&& (f^kG')f
\end{tikzcd}\]
If all functors in $ \text{Fun}(I,\mathcal{C})$ admit
left Kan extensions along $f$, then  
$f^k\colon \text{Fun}(I,\mathcal{C})\to \text{Fun}(J,\mathcal{C})$ is a functor which is left adjoint to the  precomposition with $f$  functor 
 $(-)^f\colon \text{Fun}(J,\mathcal{C})\to \text{Fun}(I,\mathcal{C})$. For example, this happens if
 $\mathcal{C}$ is closed under finite colimits and 
 $f\colon I\to J$ is  of finite-type. In this case, 
 the left Kan extension of $G$ along $f$ is given by
 $a\mapsto \text{colim}_{f\leq a} G$ (see~\cite{MR1712872}).
This description of the left Kan extension, in the case $f$ is of finite-type and $\mathcal{C}$ is   closed under finite colimits,
has several consequences. For example,  $f^kG$
is isomorphic to a functor given by a composition (see~\ref{dadgdfhg}): 
\[
\begin{tikzcd}
J\ar{r}{f\leq}  &J[f]\ar{r} &\mathcal{C}
\end{tikzcd}
\]
The restrictions of $f^kG$ to the  fibres  $(f\!\!\leq)^{-1}(f\leq a)= \{b\in J\ |\ (f\leq b) =(f\leq a)\}\subset J$, for $a$ in $J$, are  therefore isomorphic to  constant functors.
When  $I$ is finite, then, since $J[f]$ is also finite, the left Kan extension $f^kG$ factors through a finite poset $J[f]$. Such functors have been considered for example in~\cite[Definition 2.11]{EMiller1}.
\end{point}

\begin{point}\label{asfdfhgdhj}
Let $f\colon I\to J$ be a functor between posets. 
A functor $G\colon I \to\mathcal{C}$ is called $f$-\textbf{right invertible}
(or right invertible with respect to $f$)
if its left Kan extension $f^k G\colon J\to \mathcal{C}$ exists and the natural transformation $G\to (f^k G)f$, adjoint to the identity $f^kG\to f^kG$, is an isomorphism. 
A functor $F\colon J\to\mathcal{C}$ is called $f$-\textbf{left invertible} (or left invertible with respect to $f$) if it is the left Kan extension of its restriction $Ff\colon I\to\mathcal{C}$, i.e. if the natural transformation $f^k(Ff)\to F$, adjoint to the identity $Ff\to Ff$, is an isomorphism. 
Note that if $G\colon I \to\mathcal{C}$ is $f$-right invertible, then its left Kan extension $f^k G\colon J\to \mathcal{C}$ is $f$-left invertible.
Similarly, if $F\colon J\to\mathcal{C}$ is $f$-left invertible, then its restriction $Ff\colon I\to\mathcal{C}$ is $f$-right invertible. 

If $G\colon I \to\mathcal{C}$ is $f$-right invertible, then, for every $H\colon I \to\mathcal{C}$ for which $f^kH\colon J\to \mathcal{C}$ exits, the  function $f^k\colon \text{Nat}_I(H,G)\to \text{Nat}_J(f^kH,f^kG)$, mapping $\phi$ to $f^k\phi$, is a bijection.
Similarly, if $F\colon J\to\mathcal{C}$ is $f$-left invertible, then, for every $H\colon J \to\mathcal{C}$, the function $(-)^f\colon \text{Nat}_J(F,H)\to \text{Nat}_I(Ff,Hf)$, mapping $\phi$ to $\phi^f$, is a bijection.
It follows that a natural transformation $\phi\colon H\to G$, between $f$-right invertible functors $H, G\colon I \to\mathcal{C}$, is an isomorphism if and only if its left Kan extension $f^k\phi\colon f^kH\to f^kG$ is an isomorphism. 
Similarly, a natural transformation $\phi\colon F\to H$, between $f$-left invertible functors $F, H\colon J \to\mathcal{C}$, is an isomorphism if and only if its restriction $\phi^f\colon Ff\to Hf$ is an isomorphism.

A category $\mathcal{C}$ is called  $f$-\textbf{right-invertible} 
if every functor $G\colon I \to\mathcal{C}$ is $f$-right invertible.
If every functor $F\colon J\to \mathcal{C}$ is $f$-left invertible, then $\mathcal{C}$ is called $f$-\textbf{left invertible}.
If a category $\mathcal{C}$ is $f$-left invertible, then it is also $f$-right invertible. 
For example, if $\mathcal{C}$ is a category closed under finite colimits, then it is right invertible with respect to 
every subposet inclusion $f\colon I\subset J$, with $I$ poset of finite-type (see~\ref{sdgsdgjdghkfhjk}).
This is because every $a$ in $I$ is the terminal object in  $f\leq f(a)$, and consequently, for every functor $G\colon I \to\mathcal{C}$, the morphism $G(a)\to \text{colim}_{f\leq f(a)} G= (f^kG)(f(a))$ is an isomorphism. 
\end{point}

\begin{prop}\label{dfvbsfdgbsfgb} 
Let $g\colon I\to L$ and $h\colon L\to J$ be poset functors, and let 
 $F\colon J\to\mathcal{C}$ be the left Kan extension of  $G\colon I\to\mathcal{C}$
along the composition $f=hg$.
Then $F$ is  $h$-left invertible if and only if 
$Fh\colon L\to\mathcal{C}$ is $h$-right invertible.
\end{prop}
\begin{proof}
If $F$ is $h$-left invertible, then, directly from the definition, $Fh$ is  $h$-right invertible. 

Assume $Fh$ is $h$-right invertible.
This means that the left Kan extension $h^k(Fh)$ exists and the natural transformation 
$\varphi\colon Fh\to (h^k(Fh))h$, adjoint to the identity $\text{id}\colon h^k(Fh)\to h^k(Fh)$, is an isomorphism. 
To show that $F$ is $h$-left invertible, we need to prove that the natural transformation $\mu\colon h^k(Fh)\to F$, adjoint to the identity
$\text{id}\colon Fh\to Fh$, is also an isomorphim. 
 Since the composition
of $\varphi$ with the restriction $\mu^h\colon (h^k(Fh))h\to Fh$ is the identity,
this restriction $\mu^h$ is  an isomorphism. 
Consider next the following commutative diagram of sets where the horizontal arrows are the  bijections
expressing the fact that $F$ is the left Kan extension of $G$:
\[\begin{tikzcd}
\text{Nat}_I(G,\  (h^k(Fh))hg)\ar{d}[swap]{\text{Nat}_I(G,\ \mu^{hg})}\arrow[r, phantom, sloped, "\cong"]  &
\text{Nat}_J(F,\  h^k(Fh))\ar{d}{\text{Nat}_I(F,\ \mu)}
\\
\text{Nat}_I(G, Fhg) \arrow[r, phantom, sloped, "\cong"]& 
\text{Nat}_J(F,F)
\end{tikzcd}\]
Since $\mu^h$  is an isomorphism, the left vertical arrow in the above diagram is a bijection. 
It follows that so is the right vertical arrow, and consequently there is $\psi\colon F\to h^k(Fh)$
for which the composition $\mu\psi\colon F\to F$ is the identity. 
Consider the composition $\psi\mu\colon h^k(Fh)\to h^k(Fh)$.
It is a natural transformation between $h$-left invertible functors,  since $Fh$ is assumed to be $h$-right invertible.
Moreover its restriction
$(\psi\mu)^h$ is an isomorphism since $\mu^h$ is an isomorphism. 
This implies that $\psi\mu$ is also an isomorphism.
\end{proof}

\section{Dimension and parental dimension}\label{dsfsefseff}
Let $(I,\leq)$ be a poset and $x$ be an element in $I$.
In this section we propose   two ways of measuring complexity of expressing $x$ as a sup (see~\ref{dsdgdfhfgjh}).
\begin{point}\label{adsgdgfjhgh}
 Let
\[\Phi_I(x):=\left\{U\ |\ 
\begin{subarray}{c} U\subset I \text{ is finite, has a proper ancestor, $x$ is a sup of it, and }\\
\text{ $x$ is not  a sup of any  $S$ such that $\emptyset\not= S\subsetneq U$}\end{subarray}
\right\}.\]
Since $\emptyset$ belongs to $\Phi_I(x)$, this collection is non-empty. 
The following  extended (containing $\infty$) number is called the \textbf{dimension} of  $x$:
\[ \text{dim}_I(x):= \text{sup}\{ |U|\  |\   U\in \Phi_I(x)\}.\]
An element $U$ in $\Phi_I(x)$, for which $\text{dim}_I(x)=|U|$, is said to \textbf{realise} $\text{dim}_I(x)$.

Since any subset $U\subset I$, for which $x$ is a sup, is contained in $I\leq x$, there is an inequality $\text{dim}_I(x)\leq |I\leq x|$.
In particular, if $I$ is of finite-type (see~\ref{dsdgfsfghjdghj}), then  $\text{dim}_I(x)$ is finite for every $x$.

The empty set realises $\text{dim}_I(x)$ if and only if 
$\text{dim}_I(x)=0$, which happens if and only if $x$ is minimal in $I$ (see~\ref{dsdgdfhfgjh}).
The set $\{x\}$ is the only set realising $\text{dim}_I(x)$ if and only if $\text{dim}_I(x)=1$.
If  $\text{dim}_I(x)>1$, then every  set realising $\text{dim}_I(x)$ cannot contain $x$. 

Here are some examples:
\begin{itemize}
\item In the poset $[n]$, 
$\text{dim}_{[n]}(k)=\begin{cases}
0 &\text{ if } k=0,\\
1 & \text{ if } k>0.
\end{cases}$
\item In the poset $[-1,0]$,  
 $\text{dim}_{[-1,0]}(x)=\begin{cases}
0 &\text{ if } x=-1,\\
1 & \text{ if } x>-1.
\end{cases}$
\item In the poset $\mathbb{R}$, $\text{dim}_{\mathbb{R}}(x)=1$ for every $x$ in  $\mathbb{R}$.
\item In the poset $2^S$ (see~\ref{afadfhfgsh}), $\text{dim}_{2^S}(\sigma)=|\sigma|$ and 
 this dimension is realised by, for example, $\{\{x\}\ |\ x\in \sigma\}$.
 \item In the suspension  poset $\Sigma S$ (see~\ref{afadfhfgsh}),
 \[\text{dim}_{\Sigma S}(x)=\begin{cases}
2 &\text{ if  $x =\infty$ and $ |S|\geq 2$},   \\
1 & \text{ if  either  $x\in S$, or $x = \infty$ and $ |S|<2$},\\
0 & \text{ if $x= -\infty$}.
\end{cases}\]
If $|S|\geq 2$, then every subset of $S$  of size $2$ realises  
$\text{dim}_{\Sigma S}(\bigvee S)=2$.
 \end{itemize}
\end{point}

 The dimension of an element $x$ depends on the global properties of the subposet $I\leq x$. 
 This dimension can be approximated by a more tractable notion of dimension, called the parental dimension, which, for elements in a finite-type poset, depends only on their parents (see Proposition~\ref{sfgasg}). 
\begin{point}\label{aDFGSDFFHDHGJ}
Let
\[\Psi_I(x):=
\left\{ U\  |\     \begin{subarray}{c} U\subset (I<x) \text{ is finite, has an ancestor, and }\\
\text{$x$ is a sup of every two element subset of $U$}\end{subarray}\right\}.
\]
Since $\emptyset$ belongs to $\Psi_I(x)$, this collection is non-empty.
The following extended number is called the $\textbf{parental dimension}$ of $x$:
\[ \text{par-dim}_I(x):=  \text{sup}\{ |U|\  |\    U\in \Psi_I(x)\}.\]
An element $U$ in $\Psi_I(x)$, for which $\text{par-dim}_I(x)=|U|$, is said to \textbf{realise} $\text{par-dim}_I(x)$.
Here are some examples:
\begin{itemize}
\item If $I$ is $[n]$, or $[-1,0]$, or  $\mathbb{R}$, then $\text{par-dim}_{I}(x)=\text{dim}_{I}(x)$
for all $x$ in $I$.
\item In the poset $2^S$, for every $\sigma\subset S$,
$\text{par-dim}_{I}(x)=|\sigma|=\text{dim}_{I}(x)$ and the parental  dimension is realised, for example, by
$\{\sigma\setminus \{x\}\ |\ x\in \sigma\}$.
\item In the suspension poset $\Sigma S$,  $\text{par-dim}_{\Sigma S}(\infty) = \text{max}\{|S|,1\}$, and, if $S$ is non-empty, then $S$ realises
$\text{par-dim}_{\Sigma S}(\infty)$. 
If $S$ has at least three elements, then
there is a strict inequality
$\text{dim}_{\Sigma S}(\infty)=2< |S|=\text{par-dim}_{\Sigma S}(\infty)$.
\end{itemize}
 
 \end{point}
 
 \begin{prop} \label{hryfjhxfh}
Let $I$ be a poset. 
\begin{enumerate}
    \item $\text{\rm dim}_I(x)=0$ if and only  if $\text{\rm par-dim}_I(x)=0$.
    \item $\text{\rm dim}_I(x)=1$ if and only if  $\text{\rm par-dim}_I(x)=1$. 
    \item $\text{\rm dim}_I(x)\geq 2$ if and only  if $\text{\rm par-dim}_I(x)\geq 2$. 
\end{enumerate}
\end{prop}

\begin{proof}
The  empty set realises $\text{par-dim}_I(x)$ if and only if $\text{par-dim}_I(x)=0$, which happens if and only if $x$ is minimal in $I$.
This proves (1).
To show (3), assume $\mathrm{dim}_I(x)\ge 2$, and let $S$ in $\Phi_I(x)$ (see~\ref{adsgdgfjhgh}) be such that $|S|\geq 2$.
Choose $s$ in $S$. Since $x$ is not a sup of $S\setminus\{s\}$, there is
$z<x$ for which $(S\setminus \{s\})\subset (I\leq z)$. Moreover, $s<x$ and
$x$ is a sup of $\{s,z\}$, showing that $\{s,z\}$ belongs to $\Psi_I(x)$.
Thus, $\text{par-dim}_I(x)\ge 2$.
Since (1) and (3) hold, (2) must also hold.
\end{proof}

One reason why parental dimension is easier to calculate is the following:

\begin{lemma}\label{dadgsdsgfhn}
Let $I$ be a poset, $x$  its element, and $U$ be in $\Psi_I(x)$. 
Assume $\alpha\colon U\to I$ is a function such that $u\leq \alpha(u)<x$,
for every $u$ in $U$. Then $\alpha$ is injective, in particular $|\alpha(U)|=|U|$, and 
its image $\alpha(U)$  belongs to $\Psi_I(x)$.
\end{lemma}
\begin{proof}
Let $s$ and $u$ be elements in $U$. If $\alpha(s)=\alpha(u)$, then the relations $s\leq \alpha(s)=\alpha(u)\geq u$ and $\alpha(s)<x$ imply that $x$ is not a sup of $\{s,u\}$. This can happen only if $s=u$, as $U$ belongs to $\Psi_I(x)$.
Thus
$\alpha$ is injective
and  $|U|=|\alpha(U)|$.
If $\alpha(s)\not = \alpha(u)$, then a consequence of the relations $s\leq \alpha(s)<x>\alpha(u)\geq u$ and  the fact that $x$ is a sup of $\{s,u\}$ is that  $x$ is  also a sup of 
$\{\alpha(s),\alpha(u)\}$. Furthermore, any ancestor of $U$ is also an ancestor of $\alpha(U)$. Consequently $\alpha(U)$ is in $\Psi_I(x)$.
\end{proof}

Lemma~\ref{dadgsdsgfhn} can be used to prove the following proposition  which
 is the reason behind choosing the name parental dimension.
\begin{prop}\label{sfgasg}
 If $I$ is a poset of finite-type, then, for every $x$ in $I$,
\[
\text{\rm par-dim}_I(x)=\text{\rm max}\{|S|\  |\  S\subset {\mathcal P}(x)\text{ and $S$ has an ancestor}\}.
\]
\end{prop}
\begin{proof}
Since $x$ is a sup of every two elements subset of ${\mathcal P}(x)$, the right side of the claimed equality is smaller or equal than the left side.
To show the opposite inequality, consider $U$ in  $\Psi_I(x)$. 
Since $I$ is of finite-type, for all $u$ in $U$, there is $\alpha(u)$ in $\mathcal{P}(x)$ for which $u\leq \alpha(u)$. 
According to Lemma~\ref{dadgsdsgfhn},
$|\alpha(U)|=|U|$, which concludes the proof.
\end{proof}

There are two numbers  assigned to $x$ in  $I$, its dimension $\text{dim}_I(x)$  and
its parental dimension $\text{par-dim}_I(x)$. 
According to Proposition~\ref{hryfjhxfh}, if one of these dimensions is $0$ or $1$ then so is the other. 
In general, the parental dimension always bounds the dimension, as the following result shows.

\begin{prop}\label{adgfdsfhgf}
 For every  $x$ in a poset $I$, $\text{\rm dim}_I(x)\le \text{\rm par-dim}_I(x)$.
\end{prop}

\begin{proof}
The cases  $\text{\rm dim}_I(x)$ is $0$ and $1$  follow from Proposition~\ref{hryfjhxfh}.
Let $U$ be  in $\Phi_I(x)$, with $\lvert U\rvert \ge2$.  Since, for every $u$ in $U$, the element $x$
 is not a sup of $U\setminus\{u\}$, there is $s_u$ in $I$ such that
 $u'\leq s_u<x$ for every $u'$ in $U\setminus\{u\}$. If $u_0$ and $u_1$ are different elements in 
 $U$, then  $x$ being a sup of $U$ implies that  $s_{u_0}$ and $s_{u_1}$ are also different and $x$ is a sup of $\{s_{u_0},s_{u_1}\}$. Moreover, any proper ancestor of $U$ 
is also an ancestor of $S:=\{s_u\ |\ u\in U\}$.
 Thus,  $S$ has the same size as $U$ and  belongs to  $\Psi_I(x)$. 
For every element $U$ in $\Phi_I(x)$,
we  have constructed an element $S$ in  $\Psi_I(x)$ of the same size as $U$. 
\end{proof}

The suspension example in~\ref{aDFGSDFFHDHGJ} shows that the inequality in Proposition \ref{adgfdsfhgf} can be strict.

\begin{prop}\label{ghjkghjk}
Let $I$ and $J$ be  posets. For every $x$ in $I$ and $y$ in $J$:
\begin{enumerate}
    \item $\text{\rm dim}_{I \times J}(x,y) = \text{\rm dim}_I(x) + \text{\rm dim}_J(y)$,
    \item $\text{\rm par-dim}_{I \times J}(x,y) = \text{\rm par-dim}_I(x) + \text{\rm par-dim}_J(y)$.
\end{enumerate}
\end{prop}
\begin{proof}
(1)\quad 
First, we show $\text{\rm dim}_{I \times J}(x,y) \ge \text{\rm dim}_I(x) + \text{\rm dim}_J(y)$.
If $\text{\rm dim}_I(x)=0$, then $x$ is minimal in $I$.
In this case, every element of $\Phi_{I\times J}(x,y)$ is of the form $\{x\}\times U$, where $U$ belongs to $\Phi_{J}(y)$ and hence the inequality is clear.
 Assume $\text{\rm dim}_I(x)\geq 1$ and  $\text{\rm dim}_J(y)\geq 1$.
Let $S$ and $U$ belong  respectively to $\Phi_{I}(x)$  and $\Phi_{J}(y)$. Let $x'$ be a proper ancestor of $S$ and $y'$ be a proper ancestor of $U$. Then
$\lvert S\rvert +\lvert U\rvert =\lvert (S\times \{y'\})\cup (\{x'\}\times U)\rvert$. Since $(S\times \{y'\})\cup (\{x'\}\times U)\subset I\times J$ belongs to $\Phi_{I\times J}(x,y)$, we get  $\text{\rm dim}_{I \times J}(x,y)\geq \lvert S\rvert +\lvert U\rvert$, which gives  the desired inequality.

To show $\text{\rm dim}_{I \times J}(x,y) \leq \text{\rm dim}_I(x) + \text{\rm dim}_J(y)$, consider $W$ in $\Phi_{I\times J}(x,y)$. For every element $(a,b)$  in  $W$, the set $W$
cannot contain a subset of the form $\{(a',b),(a,b')\}$, where $a'\not=a$
and $b'\not = b$, otherwise $(x,y)$ would be a sup of $W\setminus\{(a,b)\}$.
 Consequently,  $\lvert W\rvert \le \lvert \text{pr}_I(W)\rvert + \lvert \text{pr}_J(W)\rvert$, where
 $\text{pr}_I$ and $\text{pr}_J$ denote the projections. Moreover, $\text{pr}_I(W)$ and $\text{pr}_J(W)$ have, respectively, $x$ and $y$ as sup, so $\lvert \text{pr}_I(W)\rvert + \lvert \text{pr}_J(W)\rvert\le \text{dim}_I(x)+\text{dim}_J(y)$.
\smallskip

\noindent
(2)\quad As before, we start by showing  $\text{\rm par-dim}_{I \times J}(x,y) \ge \text{\rm par-dim}_I(x) + \text{\rm par-dim}_J(y)$. We argue in the same way. 
Let $S$ and $U$ belong  respectively to $\Psi_{I}(x)$  and $\Psi_{J}(y)$.
 Then $(S\times \{y\})\cup (\{x\}\times U)\subset I\times J$ belongs to $\Psi_{I\times J}(x,y)$ and is of size $|S|+|U|$, which gives the desired inequality.
 
 To show the inequality $\text{\rm par-dim}_{I \times J}(x,y) \leq \text{\rm par-dim}_I(x) + \text{\rm par-dim}_J(y)$, consider $W$ in $\Psi_{I\times J}(x,y)$. Then every $(a,b)$ in $W$ is such that $(a,b)< (x,y)$ and we define  $\alpha(a,b):=(a,y)$ if $a<x$,  and $\alpha(a,b):=(a,b)$ if $a=x$. These elements are chosen so that 
 $(a,b)\le \alpha(a,b)<(x,y)$ for every $(a,b)$ in $W$.
  By Lemma~\ref{dadgsdsgfhn}, $\alpha (W)$ also belongs to $\Psi_{I\times J}(x,y)$. Since $\lvert \alpha(W)\rvert\le \lvert \text{pr}_{I} \alpha (W)\rvert+\lvert \text{pr} _{J}\alpha(W)\rvert$, we get the desired inequality.
\end{proof}

Let $S$ be a finite set. Since  $[-1,0]^{ S}$
and $[-1,0]^{ |S|}$ are isomorphic (see~\ref{afadfhfgsh}), Proposition~\ref{ghjkghjk} gives  $\text{dim}_{[-1,0]^{ S}}f=\text{par-dim}_{[-1,0]^{S}}f=
\lvert \{s\in S\mid f(s)>-1\}\rvert$.

\begin{point}
Neither the dimension nor the parental dimension  are  monotonic  in the following sense.  If  $I\subset J$ is a subposet inclusion (see~\ref{sdgsdgjdghkfhjk}) and  $x\in I$, then in general the following inequalities may fail: $\text{dim}_I(x)\leq \text{dim}_J(x) $ and $\text{par-dim}_I(x)\leq \text{par-dim}_J(x)$.
For example, consider:
\[I=\{(0,0),(1,0),(0,1),(2,2)\}\subset \mathbb{N}^2\]
\[J=\{(0,0),(1,0),(0,1),(1,1), (2,2)\}\subset \mathbb{N}^2\] Then
$\text{dim}_I(2,2)= \text{par-dim}_I(2,2)=2$, $\text{dim}_J(2,2)=\text{par-dim}_J(2,2)=1$. For this monotonicity to hold, additional assumptions need to be made, for example: 
\end{point}
\begin{prop}\label{aDFDFHFHJ}
Let $I\subset J$ be a subposet inclusion (see~\ref{sdgsdgjdghkfhjk}).
Assume  an element $x$ in $I$ has the following property: for every finite subset $S\subset I$, if $x$ is a sup   of  $S$ in $I$, then $x$  is   a sup  of $S$  in $J$. Under this assumption
$\text{\rm dim}_I(x)\leq \text{\rm dim}_J(x)$ and $\text{\rm par-dim}_I(x)\leq \text{\rm par-dim}_J(x)$.
\end{prop}
\begin{proof}
The  assumption on $x$ implies that  $\Phi_I(x)\subset \Phi_J(x)$ (see~\ref{adsgdgfjhgh}) and $\Psi_I(x)\subset \Psi_J(x)$ (see~\ref{aDFGSDFFHDHGJ}),  which gives  the claimed inequalities. 
\end{proof}

\section{Realising posets}

In this section we introduce a construction  that transforms  posets into posets, mimicking
the relation between $\mathbb{N}^r$ and $[0,\infty)^r$. 
Let $(I,\leq )$ be a poset.

\begin{point}\label{adgsdgfhjhg}
Consider $a$ in $I$ and its set of parents $\mathcal{P}(a)$. 
Table~\ref{zdfhsfgjhdgjdgh} describes how  $a$  partitions  $I$ into three disjoint subsets of $a$-\textbf{inconsistent}, $a$-\textbf{dependent}, and  $a$-\textbf{independent} elements. 
We say that  $x$ is  $a$-independent if it is not  $a$-dependent ($a\not\leq x$), but it is $p$-dependent ($p\leq x$) for some  $p$ in  $\mathcal{P}(a)$.  Observe that every parent $p$ of $a$ is $a$-independent. An element which is either $a$-dependent or $a$-independent is called $a$-\textbf{consistent}.
 This terminology is inspired by linear algebra. 
 We think of this partition as solutions of certain poset conditions, which in the case of   an upper semilattice (see~\ref{adfgsgfjfjk}) behave as solutions of a linear system: an inconsistent system has no solutions, an independent system has a unique solution, and a dependent system has many solutions.
\begin{table}[H]
\centering
\begin{tabular}{|c|c|c|}  
\hline
\multirow{ 2}{*}{$x$ is $a$-consistent} & $x$ is $a$-dependent & $a\leq x$  \\ \cline{2-3}
&  $x$ is $a$-independent  & $a\not\leq x$ and $(\mathcal{P}(a)\leq x)\not = \emptyset$\\ \hline
\multicolumn{2}{|c|}{$x$ is $a$-inconsistent}&  $a\not\leq x$ and $(\mathcal{P}(a)\leq x) = \emptyset$  \\ \hline
\end{tabular}
\caption{}
\label{zdfhsfgjhdgjdgh}
\end{table}

Figure~\ref{partition} illustrates the partitions 
of $\mathbb{N}^2$,  $\Sigma S$, and a poset of size $5$, into  $a$-inconsistent, $a$-dependent, and  $a$-independent blocks for some choices of $a$. 
\begin{figure}
    \centering
    \includegraphics[width=10cm]{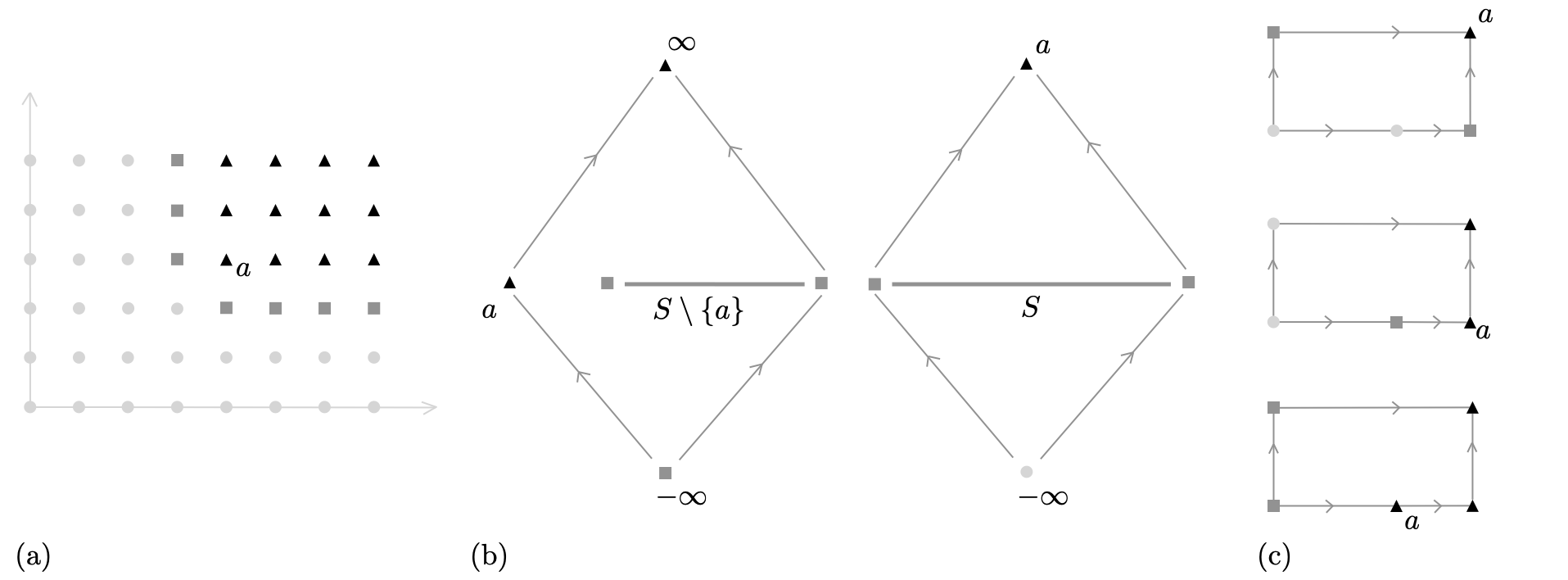}
    \caption{Black triangles are $a$-dependent, dark grey squares are $a$-independent, and light grey  dots are $a$-inconsistent}
    \label{partition}
\end{figure}

A consequence of the transitivity of a partial order is that being dependent is preserved by the poset relation:  if $a\leq b$, then a $b$-dependent element is also
$a$-dependent. Being independent does not have this property:  a $b$-independent element may fail to be $a$-independent for some $a\leq b$.
For example, consider the subposet $I=\{(0,0), (2,0),  (3,0),  (0,2), (3,2)\}\subset [0,\infty)^2$. Then, in $I$, the element $(0,2)$ is $(3,2)$-independent, however it  is not $(3,0)$-independent. The element $(0,2)$ is neither $(3,0)$-dependent and hence neither consistency is  in general preserved by the poset relation.   
\end{point}

\begin{point}\label{adgdgfjdghj}
Define  $I$ to be \textbf{consistent} if, for every $a\leq b$ in $I$, 
every $b$-independent element, having a common ancestor with $a$, is also $a$-consistent.
Thus, $I$ being consistent is equivalent to  the following condition for every $a\leq b$ in $I$: 
a $b$-consistent element $x$ is $a$-consistent if and only if there is an element $y$ in $I$ such that $x\geq y\leq a$.

Distributive upper semilattices of finite-type are key examples of 
 consistent posets (see~\ref{asfsgdhgfhjnhgn}). For example,
the poset $\mathbb{N}^r$ is consistent. If $|S|\geq 3$, then 
the suspension $\Sigma S$ (see~\ref{afadfhfgsh}) is an example of a consistent 
upper semilattice which is not distributive (see~\ref{dfhieir}).
\end{point}

\begin{point}\label{sadsfgdfghjhn}
Define $\mathcal{G}(I):=\coprod_{a\in I} [-1,0]^{\mathcal{P}(a)}$. 
We identify elements of $\mathcal{G}(I)$ with pairs $(a,f)$ consisting of an element $a$ in $I$ and a function $f\colon \mathcal{P}(a)\to [-1,0]$. 
The set $\text{supp}(f):=\{p\in \mathcal{P}(a)\  |\  f(p)<0\}$ is called the 
\textbf{support} of $f$.
If $\mathcal{P}(a)$ is empty, then the set $[-1,0]^{\mathcal{P}(a)}$
contains only one element, denoted by $0$, in which case $(a,0)$
is the only element in $\mathcal{G}(I)$ with the first coordinate equal to $a$.

For   $(a,f)$ and $(b,g)$ in $\mathcal{G}(I)$, set $(a,f)\leq (b,g)$ if
\begin{itemize}
\item[(a)] $a\leq b$ in $I$;
\item[(b)] all elements in  $\text{supp}(g)$ are $a$-consistent: if  $g(x)<0$, then either
$x$ is $a$-dependent ($a\leq x$) or  $x$ is $a$-independent ($a\not\leq x$ and $ (\mathcal{P}(a)\leq x) \not = \emptyset$);
\item[(c)] 
$\bigwedge_{[-1,0]}\{f(y)\ |\  y\in(\mathcal{P}(a)\leq x)\}\leq g(x)$ for all 
$a$-independent  $x$ in  $\mathcal{P}(b)$.
\end{itemize}

For example, $(a,f)\leq (a,g)$ in $\mathcal{G}(I)$ if and only if $f\leq g$ in $[-1,0]^{\mathcal{P}(a)}$. Thus the inclusion
$[-1,0]^{\mathcal{P}(a)}\subset \mathcal{G}(I)$, assigning to $f$ the pair $(a,f)$, is a subposet inclusion
(\ref{sdgsdgjdghkfhjk}).
If $0$ is the constant function with value $0$, then $(a,0)\leq (b,0)$ in $\mathcal{G}(I)$ if and only if $a\leq b$ in $I$. Thus the inclusion $I\subset \mathcal{G}(I)$, assigning to $a$ the pair $(a,0)$, is also  a subposet inclusion. The projection $\pi\colon  \mathcal{G}(I)\to I$, assigning to a pair $(a,f)$ the element $a$, is a functor of posets.

 We think about an element $(a,f)$ in $\mathcal{G}(I)$ as an intermediate point
between $a$ and its parents, where the value of $f(p)$ describes the time needed to go back in the direction of the parent $p$.
\end{point}

\begin{prop}\label{asdgsfghg}
 If $(a,f)\leq (b,g)$ in $\mathcal{G}(I)$, then  every  $w\leq a$ which is  an ancestor of $\text{\rm supp}(f)$ is also
an ancestor of $\text{\rm supp}(g)$.
\end{prop}

\begin{proof}
Let $w\leq a$  be an ancestor of  $\text{\rm supp}(f)$, and 
$x$ be in $\text{supp}(g)$. Then $x$ is $a$-consistent, and, thus, either $a\le x$, implying  $w\le x$, or there exists $y$ in $\mathcal{P}(a)\le x$. 
Since $\bigwedge\{f(y)\ |\  y\in(\mathcal{P}(a)\leq x)\}\leq g(x)<0$, we have $f(y)< 0$.
For such a $y$, $w\le y\le x$.
\end{proof}

To understand the relation $\leq$ on $\mathcal{G}(I)$ it is 
convenient to describe it   in an alternative way using  translations. 

\begin{Def}\label{asfadfhs} 
For $a$ in $I$, define $T_a:=[-1,0]^{\mathcal{P}(a)}$ (see~\ref{afadfhfgsh}). 
For $a\leq b$ in $I$, define  $T_{a\leq b}\colon  T_a\to T_b$ to map 
$f\colon \mathcal{P}(a)\to [-1,0]$ to $T _{a\le b}f\colon \mathcal{P}(b)\to [-1,0]$ where
\[
(T _{a\le b}f)(x):=\begin{cases} -1 & \text{if $x$ is $a$-dependent,} \\ 
\bigwedge_{[-1,0]}\{f(y)\ |\ y\in (\mathcal{P}(a)\leq x)\} & \text{if $x$ is $a$-independent,} \\ 
0 & \text{if $x$ is $a$-inconsistent.}
\end{cases}
\]
The function $T_{a\leq b}$ is called  \textbf{translation} along $a\leq b$.
\end{Def}

\begin{lemma}\label{asdgdsgshdgf}
Let    $a\le b \le c$  be in $I$.
\begin{enumerate}
\item  $T_{a\leq a}\colon T_a\to T_a$ is the identity,
\item  $T_{a\leq b}\colon T_a\to T_b$ 
is a functor of posets (see~\ref{sdgsdgjdghkfhjk}),
\item  $T_{a\le c}  \leq T _{b\le c}T_{a\le b}$. 
\end{enumerate}
\end{lemma}
\begin{proof}
Statements (1) and (2) are direct consequences of the definition. 

To prove statement (3) choose a function $f\colon\mathcal{P}(a)\to [-1,0]$ and an element $x$ in $\mathcal{P}(c)$.
If $x$  is $a$-dependent, then   $(T_{a\le c}f)(x) =-1\leq (T _{b\le c}T_{a\le b}f)(x)$.
Assume $x$ is not $a$-dependent. 
Then $x$ is not $b$-dependent.
If $x$ is $b$-inconsistent, then $(T_{a\le c}f)(x)\leq 0= (T _{b\le c}T_{a\le b}f)(x)$. Let $x$ be $b$-independent.
Then $(T _{b\le c}T_{a\le b}f)(x) = \bigwedge\{(T_{a\le b}f)(y)\ |\ y\in (\mathcal{P}(b)\leq x)\} $, and the desired inequality  $(T_{a\le c}f )(x) \leq (T _{b\le c}T_{a\le b}f)(x)$ would then follow if $(T_{a\le c}f )(x)\leq (T_{a\le b}f)(y)$ for all $y$ in $\mathcal{P}(b)\leq x$. 
Let $y$ be in   $\mathcal{P}(b)\leq x$. Then $y$ cannot be $a$-dependent, since  $x$ is  not $a$-dependent.
If $y$ is $a$-independent, then $x$ is also $a$-independent, and since
$(\mathcal{P}(a)\leq y)\subset (\mathcal{P}(a)\leq x)$, then
 $(T_{a\le c}f)(x) = \bigwedge\{f(z)\ |\ z\in (\mathcal{P}(a)\leq x)\}  \leq \bigwedge\{f(z)\ |\ z\in (\mathcal{P}(a)\leq y)\} =(T_{a\le b}f)(y)$.
 Finally, if $y$ is $a$-inconsistent, then $(T_{a\le c}f)(x)\leq 0 = (T_{a\le b}f)(y)$.  
\end{proof}

According to Lemma~\ref{asdgdsgshdgf}, translations form a lax functor (see~\ref{asdfsgadfsfhjejk}), denoted by 
$T\colon I\rightsquigarrow \text{Posets}$ and called  \textbf{$I$-translation}.
One should be aware   that in general $T$ may fail to be a functor: 
\begin{example}
Consider the  subposet $I=\{x=(0,0), a =(2,0),  b=(3,0),  y=(0,2), c=(3,2)\}\subset [0,\infty)^2$. 
Then $\mathcal{P}(a)=\{x\}$, $\mathcal{P}(b)=\{a\}$, and $\mathcal{P}(c)=\{y,b\}$. Note that $y$ is $b$-inconsistent.    
If we denote by $-0.5$ the function $\mathcal{P}(a)\to [-1,0]$ mapping $x$ to $-0.5$, then
\begin{itemize}
    \item $T_{a\le b}(-0.5)\colon \mathcal{P}(b)\to [-1,0]$ is the constant function $-1$, 
    \item $T_{b\le c}(T_{a\le b}(-0.5))
    \colon \mathcal{P}(c)\to [-1,0]$ maps  $y$ to $0$ and $b$ to $-1$,
    \item  $T_{a\le c}(-0.5)\colon \mathcal{P}(c)\to [-1,0]$   maps   $y$ to $-0.5$ 
    and  $b$ to $-1$.
\end{itemize}
In this case, $T_{a\le c}(-0.5)<T_{b\le c}T_{a\le b}(-0.5)$.
\end{example}

Since $T\colon I\rightsquigarrow \text{Posets}$  is a lax functor, we can form its 
Grothendieck construction $\text{Gr}_IT$.
As sets, $\mathcal{G}(I)$ and $\text{Gr}_IT$
are identical. The next lemma states  that the relation $\leq$ on $\mathcal{G}(I)$ 
coincides with the poset  relation  on the Grothendieck construction $\text{Gr}_IT$ described in~\ref{asdfsgadfsfhjejk}.

\begin{lemma}\label{asdgdfhfgjhg}
Let  $(a,f)$ and  $(b,g)$ be elements in $\mathcal{G}(I)$. 
\begin{enumerate}
    \item  $(a,f)\leq (b,g)$  if and only if $a\leq b$ in $I$  and $T _{a\le b}f\leq g$ in $[-1,0]^{\mathcal{P}(b)}$.
    \item If $a\le b$, then $(a,f)\leq (b,T_{a\leq b}f)$.
    \item Elements in $\text{\rm supp}(T_{a\leq b}f) $ are $a$-consistent.
\end{enumerate}
\end{lemma}
\begin{proof}
Statement (3) is a direct consequence of (2) and   (2) is a direct consequence of (1). 
To show (1), choose $x$ in $\mathcal{P}(b)$. If  $x$ is $a$-dependent, then $(T _{a\le b}f)(x)=-1\leq g(x)$. 
If  $x$ is $a$-independent, then $(T _{a\le b}f)(x)=\bigwedge\{f(y)\ |\ y\in (\mathcal{P}(a)\leq x)\}$ and hence $(T _{a\le b}f)(x)\leq g(x)$ is equivalent to condition (c) in~\ref{sadsfgdfghjhn}.  If $x$ is $a$-inconsistent, then 
$T _{a\le b}f(x)=0$ and hence $T _{a\le b}f(x)\leq g(x)$  if and only if $x$ is not in the support of $g$. These equivalences are  what is needed to show (1).
\end{proof}

Since the posets  $\mathcal{G}(I)$ and  $\text{Gr}_IT$ coincide,
Proposition~\ref{adfgsgdhjdg} gives the following explicit description
of some sup elements in $\mathcal{G}(I)$. 

\begin{cor}\label{dbhfjgyhdklsadj}
Let $S\subset \mathcal{G}(I)$ be a non-empty  subset.
\begin{enumerate}
    \item If $\pi(S)=\{a\}$, then $S$ has a coproduct: $\bigvee_{\mathcal{G}(I)} S=(a,g)$, with $g =\bigvee_{[-1,0]^{\mathcal{P}(a)}}\{f\ |\   (a,f)\in S\}$.
    \item If $b$ is a sup of $\pi(S)$ in $I$, then there is a unique sup of $S$ in 
    $\mathcal{G}(I)$ of the form $(b,m)$, and  $m=\bigvee_{[-1,0]^{\mathcal{P}(b)}}\{T_{a\leq b}f\ |\ (a,f) \text{ in } S\}$.
\end{enumerate}
\end{cor}

In this article, we are primarily interested in a subposet of $\mathcal{G}(I)$, called the realisation of $I$.

\begin{Def}\label{drjfhhj}
For an element $a$ in $I$, define $\mathcal{R}_a(I)$ to be the 
 subposet of $ \mathcal{G}(I)$ consisting of all the  pairs $(a,f)$ satisfying the following conditions:
\begin{itemize}
    \item $f(x)>-1$, for every $x$ in $\mathcal{P}(a)$;
    \item  $\text{supp}(f)$  is finite;
    \item  $\text{supp}(f)$ has an ancestor in $I$.
\end{itemize}

The subposet $\mathcal{R}(I):=\bigcup_{a\in I}\mathcal{R}_a(I)\subset \mathcal{G}(I)$
is called the \textbf{realisation} of $I$.
\end{Def}

For example, for every $a$ in $I$, the element $(a,0)$ belongs to  $\mathcal{R}_a(I)\subset \mathcal{R}(I)$. The function $I\to \mathcal{R}(I)$, mapping
$a$ to $(a,0)$, is a subposet inclusion, and we identify $I$ with its image in  $\mathcal{R}(I)$.

Let $(a,f)$ be in $\mathcal{R}(I)$  and  $(a,f)\leq (b,g)$ in $\mathcal{G}(I)$. Then, according to Proposition~\ref{asdgsfghg}, $\text{supp}(g)$ has an ancestor. 
Thus, to show that $(b,g)$ is in $\mathcal{R}(I)$ only the first two conditions 
in Definition~\ref{drjfhhj} need to be verified:  the values of $g$ are strictly bigger than $-1$ and the support of $g$ is finite. 

Let $a$ be in $I$.
For a finite subset   $S\subset \mathcal{P}(a)$ having an ancestor and a function $f\colon S\to (-1,0)$, denote by $\overline{f}\colon \mathcal{P}(a)\to (-1,0]$
the  extension:
\[\overline{f}(x):=\begin{cases}
f(x) &\text{ if } x\in S\\
0 &\text{ if } x\not\in S
\end{cases}\]
 The element   $(a,\overline{f})$ belongs to  $\mathcal{R}_a(I)\subset \mathcal{R}(I)$. The function $(-1,0)^{S}\to  \mathcal{R}_a(I)$,
mapping $f$ to   $(a,\overline{f})$, is a subposet inclusion. Since these subposets are disjoint,  as a set,  $\mathcal{R}_a(I)$ can be identified with the disjoint union 
 $\coprod_S (-1,0)^{S}$  where $S$ ranges over finite subsets of $\mathcal{P}(a)$ that have  ancestors. Thus $\mathcal{R}(I)$ can be identified with  the disjoint union 
 $\coprod_{a} \coprod_{S}(-1,0)^{S}$ where $a$ ranges over all elements of $I$ and, for each such $a$, $S$ ranges over all finite subsets of $\mathcal{P}(a)$ that have  ancestors. 

\begin{point}\label{dsgsdfhfj}
The reason we are interested in the realisation is that it generalises the relation $\mathbb{N}^r\subset [0,\infty)^r$.
Consider   $\alpha \colon \mathcal{R}(\mathbb{N}^r) \to [0,\infty)^r$ mapping $(a,f)$  
to $a+\sum_{(a-e_i)\in \mathcal{P}(a)}f(a-e_i)e_i$ 
 where $e_i$ is the $i$-th vector in the standard basis 
of $\mathbb{R}^r$.
 This function is  an isomorphism of posets. We use this  to identify
the realisation $\mathcal{R}(\mathbb{N}^r)$ with $[0,\infty)^r$. The composition
of the function $\mathbb{N}^r\to  \mathcal{R}(\mathbb{N}^r)$ mapping $a$ to $(a,0)$ with  $\alpha$ is the inclusion $\mathbb{N}^r\subset [0,\infty)^r$.
The same formula $(a,f)\mapsto a+\sum_{(a-e_i)\in \mathcal{P}(a)}f(a-e_i)e_i$ gives a poset isomorphism between  $\mathcal{R}([1]^r)$ and $[0,1]^r$.
\end{point}
\begin{point}
Let $P$ be a set.
Consider the inclusion poset $2^P$ of  subsets of $P$ (see~\ref{afadfhfgsh}). Since $2^P$ has a global minimum (see~\ref{dsdgfsgjhdgh}), given by the empty subset, 
all  subsets of $2^P$ have an ancestor.
Every parent of   $S$ in  $2^P$ is of the form  $S\setminus\{x\}$ for  $x$  in $S$, and hence the
 function mapping $x$ in $S$ to $S\setminus\{x\}$  is a bijection between $S$ and the set   of its parents $\mathcal{P}_{2^P}(S)$ in $2^P$.
We use this bijection to identify $\mathcal{P}_{2^P}(S)$ with $S$ (see~\ref{sdfgbsdhgfsh}).  
The realisation  $\mathcal{R}(2^P)$ can be identified with a  subposet of $[-1,0]^P$. 
By definition, $\mathcal{R}(2^P)$  consists of pairs
$(S,f\colon S\to (-1,0])$  (here we use the
identification between $S$ and $\mathcal{P}_{2^P}(S)$) where the support of $f$ is  finite. For such an element in $\mathcal{R}(2^P)$, define $\underline{f}\colon P\to [-1,0]$:
\[\underline{f}(x) :=\begin{cases}
f(x) &\text{ if } x\in S\\
-1 &\text{ if } x\not\in S
\end{cases}\]
The function $\mathcal{R}(2^P) \to [-1,0]^P$, mapping $(S,f\colon S\to (-1,0])$ to $\underline{f}$, is a subposet inclusion (see~\ref{sdgsdgjdghkfhjk}).
Its image consists of  functions $g\colon P\to [-1,0]$ for which  $g^{-1}((-1,0))=\{x\in P\ |\  -1<g(x)<0\}$ is finite.
The subposet in $ [-1,0]^P$ of such functions  is therefore isomorphic to  $\mathcal{R}(2^P)$. For example, if $P$ is finite, then
$\mathcal{R}(2^P)$ is isomorphic to $ [-1,0]^P$: the realisation of the discrete cube of finite dimension is 
isomorphic to the  geometric cube of the same dimension (see~\ref{afadfhfgsh}). 
\end{point}

\begin{point}\label{sadgsfghb}
The realisation $\mathcal{R}(I)$  has the following extension  property, which is convenient for constructing functors indexed by it. 
Let $U\colon I\to 2^Y$ be a  functor of posets (see~\ref{sdgsdgjdghkfhjk}) where $Y$ is a set and $2^Y$ is the inclusion poset of all subsets of $Y$ (see~\ref{afadfhfgsh}). Choose  a distance $d$ on $Y$ whose values are bounded by a real number $m$.  
 When all the values of $U$ are non-empty, we are going to use this distance  to extend $U$ along  $I\subset \mathcal{R}(I)$ to form a commutative diagram of poset functors:
\[\begin{tikzcd}[row sep=small]
I\ar[hook]{rr}\ar{rd}[swap]{U} & & \mathcal{R}(I)\ar{ld}{\overline{U}}\\
& 2^{Y}
\end{tikzcd}\]
For $(a,f)$ in  $\mathcal{R}(I)$, define:
\[\overline{U}(a,f):=U(a)\cap\bigcap_{p\in \mathcal{P}(a)}
B(U(p),(1+f(p))m)\]
where $B(V,r)=\{y\in Y\ |\ d(x,y)<r\text{ for some 
$x$ in $V$}\}$. Note that $\overline{U}(a,0)=U(a)$.
To prove the functoriality of  $\overline{U}$ all the conditions
describing  the poset relation on $\mathcal{R}(I)$ (see~\ref{sadsfgdfghjhn})  are needed.
\end{point}
\begin{prop}\label{sasDFGSDGJFHG}
Let $Y$ be a set and $U\colon I\to 2^Y$ be a functor of posets whose values are  non-empty subsets  of $Y$.
If $(a,f)\subset (b,g)$ in
$\mathcal{R}(I)$, then $\overline{U}(a,f)\subset \overline{U}(b,g)$ in $2^Y$.
Moreover   all values of $\overline{U}$ are non-empty.
\end{prop} 
\begin{proof}
Let $u$ be in $\overline{U}(a,f)$. By definition, this means:
$u$ is in $U(a)$ and, for all $p$ in $\mathcal{P}(a)$, there is  $u_p$ in $U(p)$ for which $d(u,u_p) < (1+f(p))m$.
We need to show that $u$ is in $\overline{U}(b,g)$, which is equivalent to: $u$ being in $U(b)$ and,  for all $x$ in $\mathcal{P}(b)$, there is  $v_x$ in $U(x)$ for which $d(u,v_x) < (1+g(x))m$.

The relation $(a,f)\leq  (b,g)$ yields $a\leq b$ (condition (a) in~\ref{sadsfgdfghjhn}), and hence $U(a)\subset U(b)$,
implying $u$ is in $U(b)$.
Choose $x$ in $\mathcal{P}(b)$. If $g(x)=0$, then since all distances  in $Y$ are bounded by $m$, we can choose $v_x$ to be any element in $U(x)$ (which exists 
since $U(x)$  is non-empty by assumption).
Assume $g(x)<0$.
By condition (b) in~\ref{sadsfgdfghjhn}, either $x$ is $a$-dependent ($a\leq x$)  or $x$ is $a$-independent ($a\not\leq x$ and $(\mathcal{P}(a)\leq x)\not=\emptyset$).
If $a\leq x$, then we can take $v_x=u$, since $U(a)\subset U(x)$.
If $x$ is $a$-independent, then the relation $(a,f)\leq (b,g)$ implies the existence of
$p$ in $\mathcal{P}(a)\leq x$ for which $f(p)\leq g(x)$ (here we use the fact that
$\text{supp}(f)$ is finite so that $\bigwedge_{[-1,0]}\{f(y)\ |\  y\in(\mathcal{P}(a)\leq x)\}$ is realised by $p$, see condition (c) in ~\ref{sadsfgdfghjhn}). In this case we can take $v_x= u_p$, as
$d(u_p,x)\leq (1+f(p))m\leq (1+g(x))m$.

To show non-emptiness of $\overline{U}(a,f)$, note that $\text{supp}(f)$ has an ancestor $w$. Thus $U(w)\subset U(a)$ and $U(w)\subset U(p)$ for every
$p$ in $\text{supp}(f)$ and consequently $U(w)\subset \overline{U}(a,f)$.
\end{proof}

Since $\mathcal{R}(I)$ is a subposet of $\mathcal{G}(I)$, one may also construct some sup elements in $\mathcal{R}(I)$ analogously to 
Corollary~\ref{dbhfjgyhdklsadj}.
\begin{cor}\label{bcvfhdxncvb}
Let $S\subset \mathcal{R}(I)$ be a non-empty  subset.
\begin{enumerate}
    \item If $\pi(S)=\{a\}$, then
    $S$ has a coproduct, $\bigvee_{\mathcal{R}(I)} S$, in $\mathcal{R}(I)$ and it is equal to $(a,g)$, where $g =\bigvee_{[-1,0]^{\mathcal{P}(a)}}\{f\ |\   (a,f)\in S\}$.
    \item Assume $I$ is a poset whose every element has a finite set of parents.
    If $b$ is a sup  of $\pi(S)$ in $I$, then there is a unique sup of $S$ in 
    $\mathcal{R}(I)$ of the form $(b,m)$, where $m=\bigvee_{[-1,0]^{\mathcal{P}(b)}}\{T_{a\leq b}f\ |\ (a,f) \text{ in } S\}$.
\end{enumerate}
\end{cor}
\begin{proof}
Since  $\mathcal{R}(I)$ is a subposet of $\mathcal{G}(I)$,  to prove the result it is enough to show that $(a,g)$ and $(b,m)$    
belong to $\mathcal{R}(I)$.
\smallskip

\noindent
(1):\quad For every $(a,f)$ in $S$, 
since $f\leq g$, then $g(x)>-1$ for all $x$ in $\mathcal{P}(a)$, and  $\text{supp}(g)\subset \text{supp}(f)$. Thus,  $\text{supp}(g)$ is finite and  every
ancestor of $\text{supp}(f)$ is also an ancestor of $\text{supp}(g)$. 
\smallskip

\noindent
(2):\quad 
The equality  $m(x)=-1$ holds if and only if $(T_{a\leq b}f)(x)=-1$
for all $(a,f)$ in $S$. Thus, $m(x)=-1$ if and only if 
 $a\leq x$ ($x$ is $a$-dependent), for every $a$ in $\pi(S)$. The equality  $m(x)=-1$   would then contradict the assumption that $b$ is a sup  of $\pi(S)$, since $x$ is in $\mathcal{P}(b)$.
The  values of the function  $m$  belong, therefore, to  $(-1,0]$.
The finiteness of $\text{supp}(m)$  is guaranteed by the finiteness assumption on the sets of parents in  $I$, and the existence of its ancestor is guaranteed by 
Lemma~\ref{asdgsfghg}.
\end{proof}

According to Corollary~\ref{bcvfhdxncvb}, if the set of parents of every element in $I$ is finite,
then elements  in $\mathcal{R}(I)$ of the form $(b,\bigvee\{T_{a\leq b}f\ |\ (a,f) \text{ in } S\})$, where $b$ is a sup of $\pi(S)$ in $I$, are sups of $S\subset \mathcal{R}(I)$. 
One should be aware, however, that  an element $(c,h)$ in $\mathcal{R}(I)$ may be a sup of $S$ even though $c$ is not a sup  of $\pi(S)$ in $I$ as Example~\ref{svdgsfghb} illustrates.

\begin{point}\label{sdgdfhkjkl}
An important aspect of realisations is that they admit 
explicit grid-like discretisations. These discretisations  are particularly useful to describe
properties of tame functors, such as 
their homological dimension (see~\ref{sdrtyhgf}).

For  subposets $D\subset I$ and $V\subset (-1,0)$, denote by 
$\mathcal{R}_D(I, V)\subset \mathcal{R}(I)$ the following subposet:
\[\mathcal{R}_D(I, V):=\{(a,f)\in \mathcal{R}(I)\ |\ a\in D\text{ and } f(\text{supp}(f))\subset V\}.\] 
An element  $(a,f)$ in $\mathcal{R}(I)$ belongs to $\mathcal{R}_D(I, V)$ if and only if $a$ is in $D$ and
all non-zero values of $f$ belong to $V$. In particular, 
$(a,0)$ belongs to $\mathcal{R}_D(I, V)$ if and only if $a$ is in $D$, which means that the following inclusions hold:
\[\begin{tikzcd}[row sep=small]
D\ar[hook]{r}\ar[hook]{d} & \mathcal{R}_D(I, V)\ar[hook]{d}\\
I\ar[hook]{r} & \mathcal{R}(I)
\end{tikzcd}
\]
Recall that $\mathcal{R}(I)$ can be identified with the disjoint union
 $\coprod_{a} \coprod_{S}(-1,0)^{S}$, where for every $a$ in $I$, $S\subset \mathcal{P}(a)$ is finite and have an ancestor (see the paragraph after~\ref{drjfhhj}). 
 Via this identification,
$ \mathcal{R}_D(I, V)$ corresponds to  $\coprod_{a} \coprod_{S}V^{S}$,
where $a$ ranges over all elements in $D$ and, for each such $a$,  $S$ ranges over all finite subsets of $\mathcal{P}(a)$ that have ancestors.

If $D=\{a\}$, then $\mathcal{R}_D(I, V)$ is also denoted as
$\mathcal{R}_a(I, V)$, if $D=I$, then  $\mathcal{R}_D(I, V)$ is also denoted as
$\mathcal{R}(I, V)$, and if $V=(-1,0)$, then 
 $\mathcal{R}_D(I, V)$ is also denoted as
$\mathcal{R}_D(I)$. 
For example, choose 
$V=\{-0.5\}\subset (-1,0)$.  Recall that  the  realisation
 $\mathcal{R}(\mathbb{N}^r)$ can be identified with 
 the poset
 $[0,\infty)^r$
(see~\ref{dsgsdfhfj}). 
Via this  identification, 
$\mathcal{R}(\mathbb{N}^r, V)\subset \mathcal{R}(\mathbb{N}^r) $ 
corresponds to $0.5 \mathbb{N}^r\subset
[0,\infty)^r.$

If $D$ and $V$ are  finite, then so is $\mathcal{R}_D(I, V)$. 
Furthermore, for any finite subposet  $S\subset \mathcal{R}(I)$,  there is  a finite $D\subset I$
 and a finite $V\subset (-1,0)$ for which $S\subset \mathcal{R}_D(I, V)$. 
 In case $I$ is of finite-type, we could choose $D$ to be of the form $I\leq a$, for some $a$ in $I$.
\end{point}

If $I $ is of finite-type, then the dimension and the parental dimension of every element in the realisation $\mathcal{R}(I)$ coincide and can be calculated in an analogous way as in  Proposition~\ref{sfgasg}. 

\begin{thm}\label{afashgrjh}
Let $V\subset (-1,0)$ be a subset, $I$ be a poset of finite-type, and $D\subset I$   a subposet such that  $(I\leq d)\subset D$ for every $d$ in $D$.
Assume $(a,f)$ is in $\mathcal{R}_{D}(I, V)$ and there is $\varepsilon$ in $V$ such that $\varepsilon<f(x)$, for all $x$ in $\mathcal{P}(a)$.
Then the  numbers
$\text{\rm par-dim}_{\mathcal{R}_a(I,V)}(a,f)$,
$\text{\rm par-dim}_{\mathcal{R}_{D}(I,V)}(a,f)$, 
$\text{\rm dim}_{\mathcal{R}_a(I,V)}(a,f)$ and 
$\text{\rm dim}_{\mathcal{R}_{D}(I,V)}(a,f)$
coincide and are equal to
\[
\text{\rm max}\{|S|\  |\   \text{\rm supp}(f)\subset S\subset \mathcal{P}(a) \text{ and $S$ has an ancestor}\}.
\]
\end{thm}
\begin{proof}
The vertical inequalities in the following diagram follow from Proposition~\ref{adgfdsfhgf}. The horizontal inequalities are a consequence of
Proposition~\ref{aDFDFHFHJ}, where the required assumption is given by Corollary~\ref{bcvfhdxncvb}.(1):
\[\begin{tikzcd}[column sep = 0.7em, row sep = 0.7em]
 \text{\rm par-dim}_{\mathcal{R}_{D}(I,V)}(a,f)\arrow[r,symbol=\geq]\arrow[d,symbol=\geq]   & \text{\rm par-dim}_{\mathcal{R}_a(I,V)}(a,f)\arrow[d,symbol=\geq] \\
 \text{\rm dim}_{\mathcal{R}_{D}(I,V)}(a,f)\arrow[r,symbol=\geq]  &\text{\rm dim}_{\mathcal{R}_a(I,V)}(a,f)
 \end{tikzcd}\]
 
 Let $\mathcal{A}:=\{S\  |\   \text{\rm supp}(f)\subset S\subset \mathcal{P}(a) \text{ and $S$ has an ancestor}\}$.
Next we show $\text{\rm dim}_{\mathcal{R}_a(I,V)}(a,f)\geq \text{\rm max}\{|S|\  |\  S\in \mathcal{A}\}$. Since $I$ is of finite-type, the collection $\mathcal{A}$ is finite and  all its elements are finite.
Let $S$ be in $\mathcal{A}$.
 Define  $\varepsilon_S\colon \mathcal{P}(a)\to (-1,0]$ and $f_p\colon \mathcal{P}(a)\to (-1,0]$, for every $p$ in $S$,  by the formulas:
\[
\varepsilon_S(x) = \begin{cases}
\varepsilon & \text{ if } x\in S\\
0 & \text{ if } x\not\in S
\end{cases}
\ \ \ \ \  \ \ \ \ 
f_p(x)=
\begin{cases}
f(x) & \text{ if } x=p\\
\varepsilon &\text{ if } x\in S\setminus\{p\}\\
0 & \text{ if } x\not\in S 
\end{cases}
\]

These functions are chosen so that, for every $p$ in $S$,  $-1<\varepsilon_S< f_p<f$.
Since $\text{supp}(f_p)\subset  S= \text{supp}(\varepsilon_S)$,  both sets
 $\text{supp}(f_p)$ and $\text{supp}(\varepsilon_S)$ have ancestors in $I$ and thus  
 $(a,f_p)$ and $(a,\varepsilon_S)$ belong to $\mathcal{R}_a(I,V)$.
 The relation $\varepsilon_S<f_p$, for every $p$, implies that $(a,\varepsilon_S)$ is a proper ancestor of  $U:=\{(a,f_p)\ |\ p\in S\}$  in $\mathcal{R}_a(I,V)$. If $p\not = q$ in $S$, then $f_p\not=f_q$,
 and hence $|U|=|S|$. Furthermore, $\bigvee_{\mathcal{R}_a(I,V)} U = (a,f)$ and 
 $\bigvee_{\mathcal{R}_a(I,V)} (U\setminus\{(a,f_p)\}) < (a,f)$ for every $p$ in $S$. The set $U$ belongs, therefore, to 
  the collection $\Phi_{\mathcal{R}_a(I,V)}(a,f)$  used to define the  dimension (see~\ref{adsgdgfjhgh}), and,  consequently, $\text{dim}_{\mathcal{R}_a(I,V)}(a,f)\geq |U|=|S|$.  As this happens for every $S$ in $\mathcal{A}$, we get the desired inequality.

To finish the proof we show 
 $\text{\rm max}\{|S|\  |\  S\in \mathcal{A}\}\geq \text{\rm par-dim}_{\mathcal{R}_{D}(I,V)}(a,f)$.
Let $U$  be an element in $\Psi_{\mathcal{R}_D(I,V)}(a,f)$ (see~\ref{aDFGSDFFHDHGJ}).
In particular, $U$ has an ancestor $(c,h)$ in $\mathcal{R}_D(I,V)$.
For $(b,g)$ in $U$, consider $\alpha(b,g):=(T_{b\leq a}g)\vee \varepsilon$. 
We claim  $(a,\alpha(b,g))$ belongs to $\mathcal{R}_D(I,V)$ and 
$\alpha(b,g)< f$. This is clear if $b=a$. If $b<a$, since $I$ is of finite-type, there is a parent $p$ of $a$ such that   $b\leq p<a$. In this case $T_{b\leq a}g(p)=-1$, and, consequently, $\alpha(b,g)(p)=\varepsilon<f(p)$. 
Furthermore, 
any ancestor of $\text{supp}(h)$  is also an ancestor of $\text{supp}(\alpha(b,g))$. In particular,  $S:=\bigcup_{(b,g)\in U} \text{supp}(\alpha(b,g))$  has an ancestor, and it contains $\text{supp($f$)}$, as $(T_{b\le a}g\vee \varepsilon) \le f$. Thus, $S$ is in $\mathcal{A}$.

For every $(b,g)$ in $U$, let $p_{(b,g)}$ be a parent of $a$ in $I$ for which $\alpha(b,g)(p_{(b,g)})<f(p_{(b,g)})$.
If $\alpha(b,g)\not =\alpha(b',g')$, then, since $\alpha(b,g)\vee\alpha(b',g') = f$, the elements
$p_{(b,g)}$ and $p_{(b',g')}$ have to be different.  Thus $|\{p_{(b,g)}\ |\ (b,g)\in U\}|=|U|$, and  consequently $|S|\geq |U|$ as  $\{p_{(b,g)}\ |\ (b,g)\in U\}\subset S$.  As this happens for every $U$ in  $\Psi_{\mathcal{R}_D(I,V)}(a,f)$, we get $\text{max}\{|S|\ |\ S\in \mathcal{A}\}\geq \text{par-dim}_{\mathcal{R}_D(I,V)}(a,f)$.
\end{proof}

Note that the considered dimensions in  Theorem~\ref{afashgrjh}  do not depend on
the choice of $D$ and $V$. Consequently,

\begin{cor}\label{sdfwefwefwsghj}
Let $I$ be a poset of finite-type. Then, for every element $(a,f)$ in $\mathcal{R}(I)$, the numbers
$\text{\rm par-dim}_{\mathcal{R}_a(I)}(a,f)$,
$\text{\rm par-dim}_{\mathcal{R}(I)}(a,f)$, 
$\text{\rm dim}_{\mathcal{R}_a(I)}(a,f)$ and 
$\text{\rm dim}_{\mathcal{R}(I)}(a,f)$
coincide and  are equal to
\[
\text{\rm max}\{|S|\  |\   \text{\rm supp}(f)\subset S\subset \mathcal{P}(a) \text{ and $S$ has an ancestor}\}.\]
\end{cor}

\begin{point}
Let  $I$ be a poset of finite-type and $(a,f)$ be   in  $\mathcal{R}(I)$. 
By Theorems~\ref{sfgasg} and~\ref{afashgrjh}, $ \text{par-dim}_I(a)=
\text{par-dim}_{\mathcal{R}(I)}(a,0)$ and 
$\text{ par-dim}_{\mathcal{R}(I)}(a,f)\leq \text{par-dim}_I(a)$. This last inequality can be sharp. For example, consider the subposet $I=\{(1,0,0), (1,1,0), (1,0,1), (0,1,1), (1,1,1)\}\subset \mathbb{R}^3$. Let   $a=(1,1,1)$ and $f\colon\mathcal{P}(a)\to (-1,0]$
be the function which is $0$  except for  $f(0,1,1)\neq 0$. In this case, 
$\text{ par-dim}_{\mathcal{R}(I)}(a,f)=1$ and $\text{par-dim}_I(a)=2$.
\end{point}

\begin{point}\label{ertyuiuygf}
We now aim to give a visual representation of realisations of some finite posets. 
Consider the following elements in $\mathbb{R}^2$:
$a=(0,0)$, $b=(3,0)$, $c=(0,2)$, $d=(3,2)$, $h=(2,0)$, and $k=(1,2)$, and the following subposets of  $\mathbb{R}^2$: $I_1 =\{a,b,c,d\}$, $I_2=\{a,b,c,d,h\}$, and 
$I_3=\{a,b,c,d,h,k\}$. Then their realisations are isomorphic to subposets of
$\mathbb{R}^2$ illustrated in Figure~\ref{dasgsdghsdfhjgdhj}. Note that in the second subfigure, the points in the grey square are not comparable to the points in the line segment from $h$ to $b$.  In the third subfigure, the points in the grey square are not comparable to the points in the line segments from $c$ to $k$ and from $h$ to $b$.

\begin{figure}[h]
\includegraphics[width=9cm]{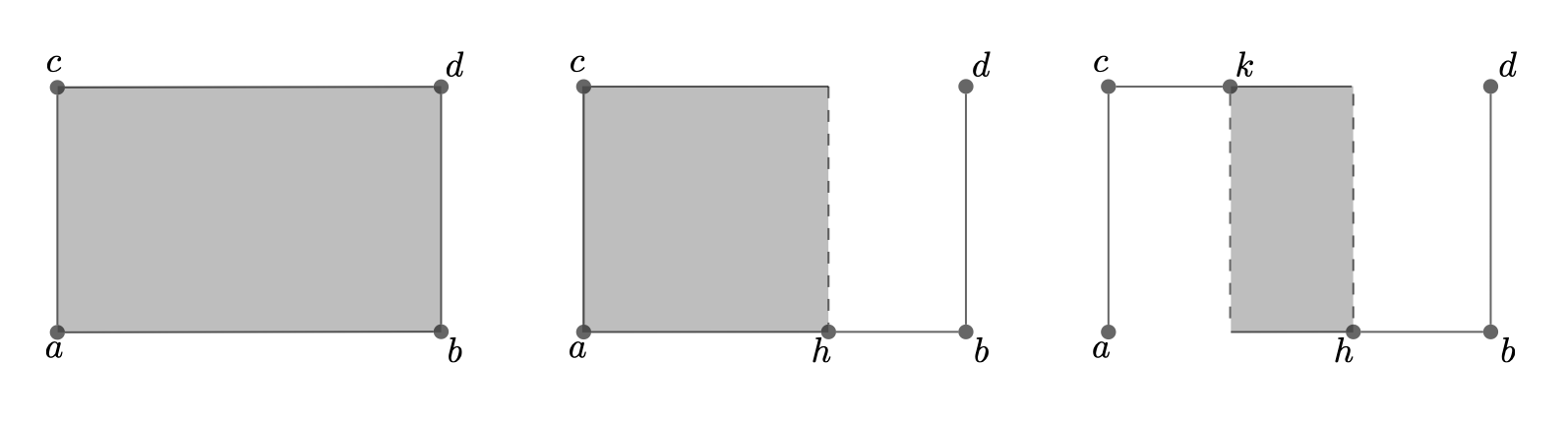}
\caption{The dashed lines are not part of the realisation poset.
}
\label{dasgsdghsdfhjgdhj}
\end{figure}
\end{point}

 \section{Upper semilattices}  
\begin{point}\label{dfhieir}
By definition an \textbf{upper semilattice} is a poset whose every  non-empty finite subset has  a coproduct (see~\ref{dsdgdfhfgjh}). To verify if a poset  is an upper semilattice it is enough to prove the existence of  the coproduct of every  two non-comparable elements.
 
An upper semilattice $I$ is called \textbf{distributive} if, for every $a$, $b$, $x$ in $I$ for which the products $a\wedge x$ and $b\wedge x$ exist, the product
$(a\vee b)\wedge x$ also exists and  $(a\vee b)\wedge x=(a\wedge x)\vee (b\wedge x)$.
 
A poset   $I$ is an upper semilattice if and only if $I_\ast$ (see~\ref{dasgfsdfhhfk}) is an upper semilattice. 
Moreover, an upper semilattice $I$ is distributive if and only if $I_\ast$ is distributive.
 
The product $I\times J$  of posets  (see~\ref{afadfhfgsh}) is an upper semilattice if and only if  $I$ and $J$ 
are  upper semilattices.  If $I$ and $J$ 
are upper semilattices, then $(x_1,y_1) \vee (x_2,y_2) =(x_1 \vee x_2,\ y_1\vee y_2)$ for  $(x_1,y_1)$, $ (x_2,y_2)$ in  $I\times J$. The product  $I\times J$ is a distributive upper semilattice if and only if $I$ and $J$ 
are distributive upper semilattices (compare with~\ref{jfsilss}). 
If $I$ is a (distributive)  upper semilattice, then so is $I^S$ for every  set $S$.  

The poset $\mathbb{R}$ is a distributive  upper semilattice which is not unital.  The posets $[n]$, $\mathbb{N}^r$, $[-1,0]^S$, and $[0,\infty)^r$ are  unital distributive upper semilattices.  

Let $S$ be a set. 
The suspension $\Sigma S$ (see~\ref{afadfhfgsh}) is a consistent and unital upper semilattice. 
It is  distributive if and only if  $|S|\leq 2$. 
For example, if $S=\{s_1,s_2,s_3\}$, then $(s_1\vee s_2)\wedge s_3=s_3$, while $(s_1\wedge s_3)\vee (s_2\wedge s_3)=-\infty$.
Thus, if  $|S|>2$, then  $\Sigma S$  is a consistent unital 
upper semilattice which is not distributive. 

For every element $a$ in  $I$, the poset $\mathcal{R}_a(I)$ (see~\ref{drjfhhj})
is a distributive upper semilattice. 
Thus, the realisation
 $\mathcal{R}(I)=\bigcup_{a\in I}\mathcal{R}_a(I)$ is a union of  distributive upper semilattices.  The realisation 
$\mathcal{R}(I)$, however, may fail to be an upper semilattice, even if $I$ is an 
upper semilattice (see Example~\ref{svdgsfghb}).
\end{point}

\begin{point}\label{koikik}
Let $I$ be  an upper semilattice and  $J$ a poset. A function $f\colon I\to J$ is called 
a \textbf{homomorphism} if, for every non-empty  finite subset $S\subset I$,  the element $f(\bigvee_I S)$ is a sup of $f(S)$  in $J$.  

If $f\colon I\to J$  is a homomorphism, then it is  a functor since if $a\leq b$ in $I$, then
$f(b)=f(a\vee b)$ is a sup of $\{f(a),f(b)\}$, which implies $f(a)\leq f(b)$ in $J$. 

If  $J$  is an upper semilattice, then a function $f\colon I\to J$ is a homomorphism if and only if $f(x\vee y) =f(x)\vee f(y)$ for all  elements $x$ and $y$ in $I$. Homomorphisms between upper semilattices are necessarily  functors.
If every two elements in $I$ are comparable, then a function $f\colon I\to J$ is a homomorphism if and only if it is a functor.
For example  $f\colon\mathbb{N}\to [0,\infty)^r$ is a homomorphism if and only it is a functor.

The inclusion $I\subset I_\ast$ is  a homomorphism (see~\ref{dasgfsdfhhfk}). Furthermore, 
a function $f\colon I\to J$ is a homomorphism if and only if 
$f_\ast\colon I_\ast\to J_\ast$ is a homomorphism.  The standard inclusion $\text{in}\colon I\to \mathcal{R}(I)$ (see~\ref{drjfhhj}) is also a homomorphism.
For every  $a$ in $I$, the poset $\mathcal{R}_a(I)$ is 
an upper semilattice (see~\ref{dfhieir}) and 
according to Corollary~\ref{dbhfjgyhdklsadj} the inclusion
$\mathcal{R}_a(I)\subset \mathcal{R}(I)$ is a homomorphism. 
For every $a$ in $I$, $T_a=[-1,0]^{\mathcal{P}(a)}$ is an upper semilattice and, for every $a\le b$ in $I$, the translation  $T_{a\le b}$  (see~\ref{asfadfhs}) is a homomorphism. 
\end{point}

 \begin{example}\label{asdfghg} Let $f\colon [1]^2\to [0,\infty)^2$ be a function defined as follows:
\[f(0,0)=(0,0),\ \ f(1,0)=(1,0),\ \ f(0,1) = (0,1),\ \ f(1,1)=(2,2).\]
This function is a unital functor, but it is not a homomorphism since:
\begin{center}$\bigvee_{[1]^2}\{(1,0),(0,1)\} =(1,1)\ \  \ \ \ \ \  \bigvee_{[0,\infty)^2}\{(1,0),(0,1)\}=(1,1)$\end{center}
\begin{center}$f(\bigvee_{[1]^2}\{(1,0),(0,1)\})=f(1,1)=(2,2)\not= (1,1) = \bigvee_{[0,\infty)^2}\{(1,0),(0,1)\}
$\end{center}
\end{example}

\begin{point}\label{jfsilss}
To construct upper semilattices, 
the Grothendieck construction (see~\ref{asdfsgadfsfhjejk})  can be used.  
Let $I$ be an upper semilattice and  $T\colon I\to \text{Posets}$ be a functor (and not just a lax functor) such that $T_a$ is an upper semilattice and $T_{a\le b}$ is a homomorphism, for every $a\le b$ in $I$. 
Then $\text{Gr}_I T$ is also an upper semilattice. 
To see this, consider  $(a,x)$ and $(b,y)$ in $\text{Gr}_I T$. 
We claim  $(a\vee b, T_{a\le a\vee b}x\vee T_{b\le a\vee b}y)$ is their coproduct in $\text{Gr}_I T$. Let $ (a,x)\leq (c,z)\ge (b,y)$, which means  $a\leq c\geq b$ in $I$ and $T_{a\le c}x\le z\geq T_{b\le c}y$  in $T_c$. Since $I$ is an upper semilattice, $a\vee b\le c$, and, since $T$ is a functor, $T_{a\vee b\le c}T_{a\le a\vee b}x=T_{a\le c}x\le z\geq  T_{b\le c}y=T_{a\vee b\le c}T_{b\le a\vee b}y$. Consequently, $(T_{a\vee b\le c}T_{a\le a\vee b}x)\vee (T_{a\vee b\le c}T_{b\le a\vee b}y)\le z$. As 
$T_{a\vee b\le c}$ is a homomorphism,  $T_{a\vee b\le c}(T_{a\le a\vee b}x \vee T_{b\le a\vee b}y)\le z$, which implies $(a\vee b, T_{a\le a\vee b}x \vee T_{b\le a\vee b}y)\le (c,z)$.
\end{point}

\begin{point}
 The dimension (see~\ref{adsgdgfjhgh}) and the parental dimension
 (see~\ref{aDFGSDFFHDHGJ}) of an element $x$ in an upper semilattice $I$ can be described using coproducts:
 
 \begin{align*} 
 \text{dim}_I(x) &=\text{sup}\left\{ |U| \  |\  \begin{subarray}{c} U\subset I \text{ is finite, has a proper ancestor,  $\bigvee U=x$,}\\
\text{and $\bigvee S< x$ for every  set $S$ such that $\emptyset\not= S\subsetneq U$}\end{subarray}\right\},\\
\text{par-dim}_I(x) &=\text{sup}\left\{ |U|\  |\   \begin{subarray}{c} U\subset (I<x) \text{ is finite, has an ancestor,}\\
\text{and $a\vee b= x$, for every $a\not= b$ in $U$}\end{subarray} \right\}.
\end{align*}
 \end{point}
 
  \begin{point}\label{afdafhghj}
 Let $I$ be an upper semilattice  of  finite-type (see~\ref{dsdgfsfghjdghj})  and $S\subset I$ be non-empty.
 If the product of $S$ exists, then it is an ancestor of $S$. On the other hand, if $S$ has an ancestor, then this ancestor  belongs to 
 $\bigcap _{x\in S}(I\le x)$, which is finite  by the finite-type assumption. Thus, if $S$ has an ancestor, then the coproduct 
  $\bigvee \left(\bigcap _{x\in S}(I\le x)\right)$
 exists ($I$ is an upper semilattice), which implies the existence of  the product $\bigwedge S$  and the equality  $\bigwedge S=\bigvee \left(\bigcap _{x\in S}(I\le x)\right)$  (see the end of~\ref{dsdgfsfghjdghj}). In conclusion, the product of a subset of a finite-type  upper semilattice  exists if and only if this subset has an ancestor. In particular,  if $\bigwedge S$ exists in $I$, then so does  $\bigwedge U$ for every  non-empty subset $U\subset S$.
 
From this  discussion and Proposition~\ref{sfgasg} it follows that the parental dimension 
of an element $x$ in $I$ can be described using products:
\begin{center}
$\text{par-dim}_I(x)=\text{max}\{|S|\ |\   S\subset {\mathcal P}(x)\text{ for which $\bigwedge S$ exists}\}$.
\end{center}
\end{point}

 \begin{point}\label{afsgsdfhdfgjhdhgj} 
Let $J$ be a poset. A subposet inclusion $f\colon I\subset J$ is called a \textbf{sublattice} if
  $I$ is an upper semilattice and $f$ is a homomorphism (see~\ref{koikik}).  For example $[0,\infty)\subset \mathbb{R} \supset \mathbb{N}$ are  sublattices.  If $J$ is 
 an upper semilattice, then
 for every $a$ in $J$, the inclusion
$(J\leq a)\subset J$ is  a sublattice. 
A sublattice $I\subset J$ satisfies the assumption of Proposition~\ref{aDFDFHFHJ} and hence $\text{\rm dim}_I(x)\leq \text{\rm dim}_J(x)$ and $\text{\rm par-dim}_I(x)\leq \text{\rm par-dim}_J(x)$, for all $x$ in $I$.

If  $J$ is an upper semilattice, then 
 the intersection of   sublattices in $J$  is also a sublattice (note that this may fail if $J$ is not an upper semilattice). 
Thus, if $J$ is an  upper semilattice, then 
 the intersection of all the sublattices  of $J$ containing a  subset $U\subset J$ is the smallest   sublattice of $J$ containing $U$.
 This intersection is denoted by  $\langle U\rangle$ and called the sublattice \textbf{generated} by $U$.  A  sublattice generated by a subset 
 $U\subset J$  is only defined  when $J$ is an upper semilattice and this is automatically assumed whenever $\langle U\rangle$ is discussed. The sublattice $\langle U\rangle$ can be described explicitly $\langle U\rangle= U\cup\{\bigvee_J S\ |\  S\subset U \text{ is finite and non-empty}\}$.
Thus if $U$ is  finite, then so is $\langle U\rangle$.
\end{point}

In Proposition~\ref{adgfdsfhgf} it was shown that the parental dimension of an element bounds  its dimension. This bound can be strict as the suspension example
(see~\ref{aDFGSDFFHDHGJ}) illustrates. 
However, for  distributive upper semilattices of finite-type dimension and parental dimension are equal.
\begin{prop}\label{gtyuyhgfdt}
Assume $I$ is a  finite-type distributive upper semilattice. Then $\text{\rm dim}_I(x)= \text{\rm par-dim}_I(x)$, for every $x$ in $I$.
\end{prop}
\begin{proof}
The cases when $\text{\rm dim}_I(x)$ is $0$ and $1$  follow from Proposition~\ref{hryfjhxfh}.
Assume $\text{par-dim}_I(x)\geq 2$ and  
 $S\subset \mathcal{P}(x)$ realises $\text{par-dim}_I(x)$. This means  $S$ is a maximal subset of $ \mathcal{P}(x)$ whose
product $\bigwedge S$ exists. The assumption $\text{par-dim}_I(x)\geq 2$ guarantees  $\lvert S\rvert\ge 2$.
 For $s$ in $S$, set $u_s:=\bigwedge(S\setminus \{s\})$ and   $U:=\{u_s\ | \ s\in S\}$. We claim: (a) $\bigwedge U=\bigwedge S$, (b) $(\bigvee U)=x$, (c) $|U|=|S|$, and (d) $(\bigvee U')<x$ if  $\emptyset\not=U'\subsetneq U$.
These properties imply  that $U$ belongs to  $\Phi_I(x)$ (see~\ref{dsdgfsfghjdghj}) and hence 
$\text{dim}_I(x)\geq |U|=|S|=\text{par-dim}_I(x)$, which together with Proposition~\ref{adgfdsfhgf} gives the desired equality.

Property (a) is clear.  Assume  (b) does not hold and  there is  $p$ in $\mathcal{P}(x)$ such that $\bigvee U\le p<x$,  which means
$u_s\leq p$ for every $s$ in $S$. If $p$ is not in $S$, then $\bigwedge S\le \bigwedge(S\setminus \{s\})\le p$ would imply existence of the product $\bigwedge (S\cup\{p\})$  contradicting  the maximality of $S$. The element $p$ is therefore in  $S$, and  $\bigwedge(S\setminus \{p\})=u_p\le p$,
leading to $u_p\vee p=p$. Since $s\vee p=x$ for every $s$ in $S\setminus \{p\}$,  distributivity gives  a contradiction
$x=\bigwedge_{s\in S\setminus \{p\}} (s\vee p)\leq u_p\vee p =p$. 

To prove property (c) we  need to show that if $s\not=s'$, then $u_s\not= u_{s'}$. 
Assume  this is not the case and $u_s= _{s'}$.
Distributivity and the already proven property (b)
leads to a contradiction:
\begin{center}
$s'=x\wedge s' = (\bigvee U)\wedge s' = \bigvee_{s} (u_s \wedge s')= (u_{s'}\wedge s') \vee (\bigvee _{s\neq s'}(u_s\wedge s')) = (\bigwedge S) \vee  (\bigvee _{s\neq s'}u_s)=\bigvee _{s\neq s'}u_s= \bigvee _{s}u_s=\bigvee U=x
$\end{center}
Thus, $U$ and $S$ have the same size, which coincide with property (c).

It remains to show  (d): $\bigvee U'<x$   if  $\emptyset\not=U'\subsetneq U$. 
Assume by contradiction that there is $s'$ in $S$ such that 
$\bigvee _{s\neq s'}u_s=x$. 
Then, again, distributivity and property (b) leads to a contradiction:
\[
s' = x\wedge s'=(\bigvee _{s\neq s'}u_s)\wedge s' = \bigvee _{s\neq s'}(u_s\wedge s') = \bigvee _{s\neq s'}u_s=x. \qedhere
\]
\end{proof}

\begin{point}
We finish this section with
a characterisation  of 
finite-type posets whose elements have dimension at most one.  A poset $I$ is called 
 a \textbf{forest} if  every pair of  non-comparable elements $x$ and $y$  in $I$ has no common ancestor. Every subposet of a forest is still a forest. A forest is called a \textbf{tree} if it is connected. Therefore, every connected subposet of a forest is a tree.
 \end{point}

\begin{prop}\label{sdagdsfjthj}
\begin{enumerate}
    \item An upper semilattice whose every element has  dimension (see~\ref{adsgdgfjhgh}) at most $1$ is a tree. 
    \item A tree of finite-type is an  upper semilattice whose every element has dimension at most $1$.
\end{enumerate}
\end{prop}

\begin{proof}
\noindent
(1):\quad First, observe that upper semilattices are connected. 
Let $x$ and $y$ be non-comparable elements in an upper semilattice $I$ whose every element has dimension at most $1$. If $x$ and $y$ had a common ancestor, then
$\{x,y\}$ would belong to $\Phi_I(x\vee y)$ (see~\ref{adsgdgfjhgh}) leading to a contradiction  $\text{dim}_I(x\vee y)\ge 2$.
\smallskip

\noindent
(2):\quad
Let $I$ be a tree of finite-type and $x$, $y$ be its elements. 
We need to show the existence of $x\vee y$.
First, assume that $x$ and $y$ have a common descendent. Then, since $I$ is of finite-type,  $\{x,y\}$ has a sup. We claim that this sup is unique, otherwise  two such sup elements  would be non-comparable and have $x$ as a common ancestor contradicting 
the assumption  $I$ is a tree.  

Assume $x$ and $y$ do not have a common descendent.
Since $I$ is connected, there is a sequence $C=\{c_0=x,c_1,\dots, c_n=y\}$  in $I$ where
$c_i$ is comparable to $c_{i+1}$ for every $0\leq i<n$.
This set $C$ does not have a descendent, otherwise $\{x,y\}$ would also have it. Thus 
$C$  must have at least two maximal elements.
Among these, there are $c_i$ and $c_j$ such that $(C\le c_i)\cap (C\le c_j)\neq \emptyset$, otherwise $C$ would not be connected. The elements 
 $c_i$ and $c_j$ are not comparable and have no common ancestor. This contradicts again the assumption  $I$ is a tree. 
 
 Finally, for every $x$ in $I$, $\text{dim}_I(x)\leq 1$, since any subset of size at least  $2$ with $x$ as sup cannot have a common ancestor by hypothesis.
\end{proof}

\section{Realisations of upper semilattices}
The realisation 
of an upper semilattice  may fail to be an upper semilattice (see Example~\ref{svdgsfghb}). 
This section is devoted to discussing some of the reasons for this and what assumptions can eliminate them. Our aim  is to prove:
\begin{thm}\label{dasghfgjk}
If $I$ is a  distributive  upper semilattice of finite-type, then its 
 realisation $\mathcal{R}(I)$ (see~\ref{sadsfgdfghjhn}) is also a distributive  upper semilattice.
\end{thm}

An upper semilattice whose elements have dimension at most $1$ is  distributive. 
In particular, finite-type trees are distributive upper semilattices (see~\ref{sdagdsfjthj}.(2)). Theorem~\ref{dasghfgjk} implies therefore that the realisation of a finite-type tree is a distributive upper semilattice, and since its elements have  dimensions at most $1$ 
it is a tree (see~\ref{sdagdsfjthj}.(1)).

\begin{point}\label{adfgsgfjfjk}
Let $I$ be an upper semilattice and $a$, $x$ be in $I$.
 If   $\mathcal{P}(a)\leq x$ contains two different elements 
$y_1$ and $y_2$,  then $y_1\vee y_2=a$, which implies $a\le x$, and hence   $(\mathcal{P}(a)\leq x)=\mathcal{P}(a)$.
Thus, there are three possibilities: the set $\mathcal{P}(a)\leq x$ is empty, consist of only one element, or is the entire set $\mathcal{P}(a)$. 
Furthermore, if $\mathcal{P}(a)\leq x$ contains only one element and $a\not\leq  x$, then this element has to be the product $a\wedge x$.  
These facts are used  to give in  Table~\ref{vabsfgdbnsfgnds} a characterisation of the blocks in the partition of $I$  described in~\ref{adgsdgfhjhg}.

\begin{table}[htbp]
\centering
\begin{tabular}{|c|c|c|}  
\hline
\multirow{ 2}{*}{$x$ is $a$-consistent} & $x$ is $a$-dependent &  \parbox[t]{6cm}{ \centering $a\leq x.$ \\
In this case $(\mathcal{P}(a)\leq x) =  \mathcal{P}(a)$}\\ \cline{2-3}
&  $x$ is $a$-independent  &  \parbox[t]{6cm}{\centering $a\not\leq x$ and  the product $a\wedge x$ exists and  belongs to $\mathcal{P}(a).$
\\ In this case  $\lvert\mathcal{P}(a)\leq x\rvert =1.$}\\ \hline
\multicolumn{2}{|c|}{$x$ is $a$-inconsistent}&  \parbox[t]{6cm}{\centering$a\not\leq x$ and either the product $a\wedge x$ does not exist or it exists but does not belong to  $\mathcal{P}(a).$\\
In this case $(\mathcal{P}(a)\leq x)=\emptyset.$} \\ \hline
\end{tabular}
\caption{}
\label{vabsfgdbnsfgnds}
\end{table}
This table can be  used to describe the translation operation (see~\ref{asfadfhs}) more explicitly.
Let  $a\leq b$ in $I$ and $f\colon\mathcal{P}(a)\to [-1,0]$ be a function. Then:
\[
(T _{a\le b}f)(x)=\begin{cases} -1 & \text{if $x$ is $a$-dependent,} \\ 
f(a\wedge x) & \text{if $x$ is $a$-independent,} \\ 
0 & \text{if $x$ is $a$-inconsistent.}
\end{cases}
\]
\end{point}

Similarly, a simplified  characterisation of an upper semilattice of finite-type to be
consistent (see~\ref{adgsdgfhjhg}) can be given.
\begin{lemma}\label{sadfgadfhg}
Let $I$ be an upper semilattice of finite-type. Then 
 $I$ is consistent if and only if the following condition is satisfied for every $a\leq b$ in $I$: 
 if $x$ is $b$-independent and the  product $a\wedge x$ exists, then $x$ is $a$-consistent, in which case either
$a\wedge x=a$ ($x$ is $a$-dependent) or $a\wedge x$ is a parent of $a$ ($x$ is $a$-independent). 
 
\end{lemma}
\begin{proof}
Consider elements $a\leq b$ in $I$. 
Assume $I$ is a consistent upper semilattice and $x$ is $b$-independent.   
 If the product $a\wedge x$ exists, then $a\wedge x$ is a common ancestor of $a$ and $x$,
 and consistency of $I$ implies that $x$ is $a$-consistent.
 
  Assume  the condition in the lemma holds. Let  $x$ be $b$-independent and let it have a common ancestor with $a$.  Since $I$ is of finite-type, the product $a\wedge x$ exists.
  The $a$-consistency of $x$ is then given by the assumed condition.
\end{proof}

\begin{prop}\label{asfsgdhgfhjnhgn}
A distributive  upper semilattice of finite-type is consistent. 
\end{prop}
\begin{proof}
To prove the proposition we  use Lemma~\ref{sadfgadfhg}.
Let $a\leq b$  and $x$ be an element which is  $b$-independent 
for which the product $a\wedge x$ exists. Let  $p$ be a parent of $b$ such that $p\leq x$.   
Then  $p\leq p \vee (a\wedge x) \leq b$, and  since $b\not\leq x$, the element $b$ cannot be equal to  $p \vee (a\wedge x)$. Thus,   
$a\wedge x\leq p\vee (a\wedge x)=p$. There are three options: either $a\wedge x=p$, or $a\wedge x$ is a parent of $a$,   or there is $y$ such that
$a\wedge x<y<a$.  The first two options give $a$-consistency of $x$  proving the proposition. 
We claim that the third option is not possible. For $y$ such that $a\wedge x<y<a$, we have $y\not\leq p$, otherwise $y\leq x$ implying    $y\leq a\wedge x$. 
Since   $y\not\leq p$, then $p\vee y =b$. 
We can then use   distributivity  to obtain  a contradiction $a = b\wedge a = (p\vee y)\wedge a = (p\wedge a) \vee (y\wedge a)\leq (a\wedge x)\vee y\leq y$.
\end{proof}

The key reason why we are interested in consistent upper semilattices is the following:
\begin{lemma}\label{ddsgadsgh}
 If $I$ is a consistent upper-semilattice, then the $I$-translation is a functor: $T_{a\leq c}=T_{b\leq c}T_{a\leq b}$ for every $a\leq b\leq c$ in $I$.
\end{lemma}
\begin{proof}
Consider a function $f\colon\mathcal{P}(a)\to[-1,0]$.
According to Lemma~\ref{asdgdsgshdgf} we just need to show $(T_{a\leq c}f)(x)\geq (T_{b\leq c}T_{a\leq b}f)(x)$ for all $x$ in $\mathcal{P}(c)$. Let $x$ be in $\mathcal{P}(c)$.
 If $x$ is $a$-inconsistent, then $(T_{a\leq c}f)(x)=0$ and the desired inequality is clear. 

Assume $x$ is $a$-consistent. Then  $a\wedge x$ exists and either $a=a\wedge x$ or $a\wedge x$ is a parent of $a$.
Since $x$ is a parent of $c$, it is $c$-consistent. Consistency of $I$ applied to $b\leq c$ 
implies the $b$-consistency  of $x$. Thus, either $x$ is $b$-dependent, or
 $b$-independent. If $x$ is $b$-dependent, then  $(T_{b\leq c}T_{a\leq b}f)(x)=-1$ and the desired inequality is clear.
Assume $x$ is $b$-independent. Then $(T_{b\leq c}T_{a\leq b}f)(x)=(T_{a\leq b}f)(b\wedge x)$.
The product $b\wedge x$ is $a$-dependent if and only if $x$ is $a$-dependent. If this happens, then
$(T_{a\leq c}f)=-1=(T_{a\leq b}f)(b\wedge x)$ and the desired equality holds in this case as well.
The product $b\wedge x$ is $a$-independent if and only if $x$ is $a$-independent. If this happens, then
$(T_{a\leq c}f)=f(a\wedge x)=f(a\wedge b\wedge x)=(T_{a\leq b}f)(b\wedge x)$  and the desired equality also  holds.
\end{proof}

For every $a\leq b$ in a poset $I$, the translation $T_{a\leq b}$  is a homomorphism of upper semilattices. If, in addition, $I$ is a consistent upper semilattice, then
 by Lemma~\ref{ddsgadsgh}, the $I$-translation  is a functor, and consequently (see~\ref{jfsilss}):
 \begin{prop}\label{asfdfjk}
 If $I$ is a consistent upper semilattice, then the \allowbreak Grothendieck construction $\text{\rm Gr}_IT=\mathcal{G}(I)$ is an upper semilattice.
 \end{prop}
 
 For the realisation $\mathcal{R}(I)\subset \mathcal{G}(I)$ (see~\ref{sadsfgdfghjhn}) to be an upper semilattice
 we need a further  finiteness restriction:
 
\begin{prop}\label{sdfgfds}
Let $I$ be a consistent upper semilattice whose every element has a finite set of parents. Then  $\mathcal{R}(I)$ is  an upper semilattice. 
\end{prop}
\begin{proof}
Let $(a,f)$, $(b,g)$ be  in $\mathcal{R}(I)$ and $m:=\bigvee\{T_{a\leq a\vee b}f,T_{b\leq a\vee b}g\}$. The poset $\mathcal{G}(I)$ is an upper semilattice (see~\ref{asfdfjk}). Thus 
 $(a\vee b,m)$ is the coproduct of $(a,f)$ and $(b,g)$ in $\mathcal{G}(I)$ (see~\ref{dbhfjgyhdklsadj}).
Since 
$(a\vee b,m)$ belongs to $\mathcal{R}(I)$ (see~\ref{bcvfhdxncvb}),
it is  the  coproduct of $(a,f)$ and $(b,g)$   in $\mathcal{R}(I)$.
\end{proof}

\begin{cor}\label{asfhgfgkkl}
Let $I$ be a consistent upper semilattice whose every element has a finite set of parents. Choose an element $d$ in $I$ and a  subposet $V\subset (-1,0)$. Then the inclusion $\alpha\colon \mathcal{R}_{I\leq d}(I,V)\subset \mathcal{R}(I)
$  is a sublattice (see~\ref{afsgsdfhdfgjhdhgj}), i.e.,
$\mathcal{R}_{I\leq d}(I,V)$ is an  upper semilattice, and $\alpha$ is a
homomorphism.
\end{cor}
\begin{proof}
Just note that if $(a,f)$ and $(b,g)$ belong to $\mathcal{R}_{I\leq d}(I,V)$, then so does their coproduct $(a\vee b, m)$ described in the proof of~\ref{sdfgfds}.
\end{proof}

We are now ready to prove Theorem~\ref{dasghfgjk}.
\begin{proof}[Proof of Theorem~\ref{dasghfgjk}]
Let $I$ be a distributive  upper semilattice of finite-type. 
Its realisation $\mathcal{R}(I)$ is an upper semilattice because
of Propositions~\ref{asfsgdhgfhjnhgn} and~\ref{sdfgfds}.
It remains  to show that $\mathcal{R}(I)$ is distributive.
This requires  understanding of  products in $\mathcal{R}(I)$.
Consider the elements  $(a,f)\geq (c,m)\leq (x,h)$ in $\mathcal{R}(I)$.
Then $a\geq c\leq x$ and hence the product $a\wedge x$ exists (see discussion in~\ref{afdafhghj}). 
Let $y$ be a parent of $a\wedge x$. If $p_1$ and $p_2$ are different parents of
$a$ for which $p_1 \wedge x = y =p_2\wedge x$, then, by the distributivity, we would get a contradiction
$a\wedge x= (p_1\vee p_2)\wedge x =(p_1\wedge x)\vee (p_2\wedge x)=y$. 
Thus, the set $\{p\in \mathcal{P}(a)\  |\ p\wedge x=y\}$
 contains at most one element. If $\{p\in \mathcal{P}(a)\  |\ p\wedge x=y\}$ is non-empty, then its unique element is denoted by $yx^{-1}$.
 Define $R(f)\colon\mathcal{P}(a\wedge x)\to(-1,0]$
 \[R(f)(y):=
 \begin{cases} 
 f(yx^{-1}) &\text{ if }\{p\in \mathcal{P}(a)\  |\ p\wedge x=y\}\not=\emptyset\\
 0&\text{ otherwise}
 \end{cases}\]
 
 Let $w$ be an ancestor of $\text{supp}(f)$.
 For every $y$ in  $\text{supp}(R(f))$, the set  $\{p\in \mathcal{P}(a)\  |\ p\wedge x=y\}$ is non-empty and
 $R(f)(y)=f(yx^{-1})<0$. Thus,   $w\leq yx^{-1}$ which implies $w\wedge x\leq y$. The element $w\wedge x$
 is therefore an ancestor of $\text{supp}(R(f))$ and   $(a\wedge x, R(f))$ is in $\mathcal{R}(I)$.
 
We claim that  the following relations hold $(c,m)\leq (a\wedge x,R(f))\leq (a,f)$.
 Let  $y$ be in  $\text{supp}(R(f))$. In particular,  $f(yx^{-1})<0$. The relation
 $(c,m)\leq (a,f)$  has two  consequences. First,  the product $yx^{-1}\wedge c=y\wedge c$ exists and 
 $y$ is  $c$-consistent (see Table~\ref{vabsfgdbnsfgnds}). All the elements of  $\text{supp}(R(f))$ are
 therefore $c$-consistent.  Second, if $y$ in addition is $c$-independent,
  then $yx^{-1}$ is also $c$-independent and $m(y\wedge c)=m(yx^{-1}\wedge c)\leq f(yx^{-1})=R(f)(y)$. These two consequences give the  relation $( c,m)\leq (a\wedge x,R(f))$.
 By direct calculation:
 \[T_{a\wedge x\leq a}R(f)(p)=
 \begin{cases}
 -1 &\text{ if  $p$ is $a\wedge x$-dependent}\\
 f(p)&\text{ if  $p$ is $a\wedge x$-independent}\\
 0 & \text{ if  $p$ is $a\wedge x$-inconsistent}
 \end{cases}
 \]
 If $p$ in $\mathcal{P}(a)$ is $a\wedge x$-inconsistent, then it is also $c$-inconsistent (see Table~\ref{vabsfgdbnsfgnds}). Thus, $T_{a\wedge x\leq a}R(f)\leq f$, as all the elements in $\text{supp}(f)$ are $c$-consistent. 
 This shows $(a\wedge x,R(f))\leq (a,f)$ (see~\ref{asdgdfhfgjhg}).
 
In an analogous way
$(c,m)\leq (a\wedge x,R(h))\leq (x,h)$ and  the obtained relations 
 imply:
  \[\begin{tikzcd}[column sep=10,row sep=10]
 (a,f) \arrow[r,symbol=\geq]  & (a\wedge x,\bigwedge\{R(f),R(h)\})\arrow[r,symbol=\leq] \arrow[d,symbol=\geq] & (x,h)\\
 & (c,m)
 \end{tikzcd}\]
 This proves  $(a\wedge x,\bigwedge\{R(f),R(h)\})$ is the product of $(a,f)$ and $(x,h)$.
 In particular, the product of two elements in $\mathcal{R}(I)$ exists if and only if these
 elements have a common ancestor.
 
 We are now ready to prove distributivity of  $\mathcal{R}(I)$.
 Let $(a,f)$, $(b,g)$, $(x,h)$ be  in $\mathcal{R}(I)$ for which the products 
 $(a,f)\wedge (x,h)$ and $(b,g)\wedge (x,h)$
 exist. Since  $(a,f)\wedge (x,h)$ is an  ancestor of both $(a,f)\vee (b,g)$ and $(x,h)$, 
 the product $((a,f)\vee (b,g))\wedge (x,h)$ exists. Moreover, this product is of the form
 $((a\wedge x)\vee (b\wedge x),m)$, where:
 \[m= \bigwedge\{R(\bigvee\{T_{a\leq a\vee b}f, T_{b\leq a\vee b}g\}), R(h) \}.\]
 By direct verification $m= \bigvee \{\bigwedge\{R(T_{a\leq a\vee b}f), R(h)\}, 
 \bigwedge\{R(T_{b\leq a\vee b}g),$ $ R(h)\}\}$, which gives the desired 
 equality  $((a,f)\vee (b,g))\wedge (x,h)= ((a,f)\wedge (x,h)) \vee ((b,g)\wedge (x,h))$.
\end{proof}

\begin{figure}
\centering
\includegraphics[width=9cm]{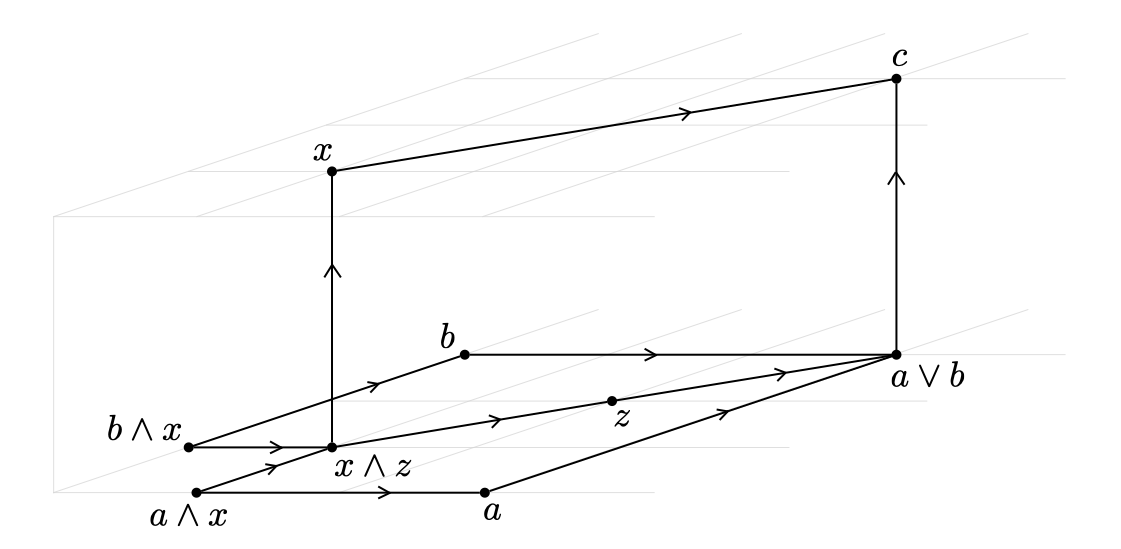}
\caption{The arrows go from smaller to bigger elements of the upper semilattice.}
\label{cons}
\end{figure}

\begin{examples}\label{svdgsfghb}
The subposet $I=\{a=(0,0), b=(3,0), c=(0,2), h=(2,0), d=(3,2)\}\subset [0,\infty)^2$ is an example of a non-consistent  upper semilattice,  which according to Proposition~\ref{asfsgdhgfhjnhgn}, cannot be distributive. Indeed, $b$ is $d$-independent, but it is not $c$-consistent, even though $c\le d$ and $b\wedge c=a$ exists. However, its realisation is still an upper semilattice (see~\ref{ertyuiuygf}). 

Consider the subposet $I\subset [0,\infty)^3$ given by the points $\{a\wedge x=(1,0,0), b\wedge x=(0,1,0), a=(3,0,0), b=(0,3,0), x\wedge z=(1,1,0), z=(2,2,0), a\vee b=(3,3,0), x=(1,1,2), c=(3,3,2)\}$, where the names of the points are chosen to help understanding the relations between them (see Figure~\ref{cons}). 
Analogously, the poset  $I$ is a non-consistent upper semilattice. Furthermore, its realisation $\mathcal{R}(I)$ is not an upper semilattice. This can be seen by taking the points $(a,f\colon\{a\wedge x\}\to (-1,0])$ and $(b,g\colon\{b\wedge x\}\to (-1,0])$ of $\mathcal{R}(I)$, where $f(a\wedge x)=-0.25$ and $g(b\wedge x)=-0.5$. We can define $(a\vee b,m\colon\{a,b,z\}\to (-1,0])$, where $m(a)=0$, $m(b)=0$ and $m(z)=f(a\wedge x)\vee g(b\wedge x)=-0.25$. Consider also $(c,h\colon\{a\vee b,x\}\to (-1,0])$, with $h(a\vee b)=0$ and $h(x)=-0.25$. The two elements $(a\vee b,m)$ and $(c,h)$ are both sups of $\{(a,f),(b,g)\}$. They are indeed not comparable, since $x$ in $\text{supp($h$)}$ is $(a\vee b)$-inconsistent.
\end{examples}

\bigskip
\noindent
\textbf{Part II, tameness.}  In this part we introduce tame  functors indexed by posets. An important property of a tame functor is that to describe it only a finite amount of  information 
is needed. 
For example, every functor indexed by a finite poset is tame.  
Tameness is therefore  an interesting property for functors indexed by infinite posets, for example by realisations.

\section{Transfers and left Kan extensions}\label{asfsdfhfhgjmd}

 \begin{point}\label{sDGSDFHFGJ}
 Let $I$ be an upper semilattice (see~\ref{dfhieir}) and $f\colon I\to J$  a function of finite-type (see~\ref{dsdgfsfghjdghj}).
 Consider the unital function $f_\ast\colon I_\ast\to J_\ast$ (see~\ref{dasgfsdfhhfk}).
 Since $I_\ast$ is also an upper semilattice and $f_\ast \leq  a$ is non-empty and finite, 
 we can form the coproduct
$\bigvee_{I_\ast} (f_\ast \leq  a)$ for all $a$ in $J_\ast$. Note that
\[\text{$\textstyle   \bigvee_{I_\ast} (f_\ast \leq  a)$}=
\begin{cases}
-\infty &\text{ if  $a=-\infty$, or $a$ is in $J$ and $(f\leq a)=\emptyset$,}\\
\bigvee_{I} (f \leq  a) & \text{ if $a$ is in $J$ and $(f\leq a)\not=\emptyset$.}
\end{cases}\]

If $a\leq b$ in $J_\ast$, then $\bigvee_{I_\ast} (f_\ast \leq  a)\leq \bigvee_{I_\ast} (f_\ast \leq  b)$ since
$(f_\ast \leq  a)\subset (f_\ast \leq  b)$.
The function mapping $a$ in $J_\ast$ to  $\bigvee_{I_\ast} (f_{\ast} \leq  a)$  in $I_\ast$ is therefore a unital functor, which we denote by $f^{!}\colon J_\ast\to I_\ast$ and call  the  \textbf{transfer} of $f$. The  transfer  $f^{!}\colon J_\ast\to I_\ast$ is only defined when
 $f\colon I\to J$  is a   function of finite-type between an upper semilattice $I$ and a poset $J$. Whenever the transfer
  $f^{!}$ is considered,  we automatically assume that these conditions are satisfied. 
 
Since $f_{\ast} \leq f_{\ast}(a)$  contains $a$, 
 the relation $a\leq \bigvee_{I_\ast} (f_\ast \leq  f_\ast (a))=f^{!}f_\ast (a)$
 holds for every $a$ in $I_\ast$. Consequently,  $\text{id}_{I_\ast}\leq f^{!} f_\ast$ holds in the poset ${I_\ast}^{I_\ast}$.
  
 In general,  $f^!(a)= \bigvee_{I_\ast} (f_\ast\leq a)$ may fail to belong to $f_\ast\leq a$, i.e., the relation $f_\ast f^{!}(a)=f_\ast ( \bigvee_{I_\ast} (f_\ast \leq a))\leq a$ may fail to hold. 
 For example, when $f\colon [1]^2\to [0,\infty)^2$ is given as in Example~\ref{asdfghg} and $a=(1,1)$ in $[0,\infty)^2$.
However, if  $f$  is a homomorphism  of finite-type, then:
 \end{point}
\begin{prop}\label{asfadfhgfgmn}
 Let $f\colon I\to J$  be a  homomorphism (see~\ref{koikik}) of finite-type between 
  an upper semilattice $I$ and a poset  $J$ (see~\ref{koikik}).
 \begin{enumerate}
\item $f_\ast f^{!}\leq \text{\rm id}_{J_\ast}$  in ${J_\ast}^{J_\ast}$.
\item For every $a$ in $J_\ast$, the element $f^!(a)$ belongs to $f_\ast\leq a$ and it is its global maximum.
\item $f^{!}f_\ast f^{!} = f^{!}$.
\item  $f_\ast f^{!}f_\ast  = f_\ast $.
\item  For every $a$ in $J_\ast$, the element $f_\ast f^!(a)$ belongs to $\{b\in J_\ast\ |\ (f_\ast\leq a)=(f_\ast\leq b)\}$ and it is its global minimum.
\item For $a$ and $b$ in $J$, $(f\leq a)=(f\leq b)$ if and only if $f^!(a)=f^!(b)$.
\item If $f$ is injective, then  $f^{!}f_\ast= \text{\rm id}_{I_\ast}$.
 \item If $f$ is surjective, then $f_\ast f^{!}  = \text{\rm id}_{J_\ast}$.
\end{enumerate}
 \end{prop}
 \begin{proof}
 \noindent
 (1):\quad Let $a$ be in $J_\ast$. Since $f$ is a homomorphism, then so is $f_\ast$, and therefore  
 $f_\ast f^{!}(a) = f_\ast (\bigvee_{I_\ast} (f_\ast \leq a))$ is a sup of
 $\{f_\ast (x)\  |\ x\in ( f_\ast \leq a)\}$ in $J_{\ast}$. The desired relation,
 $f_\ast f^{!}(a)\leq a$, is then
 a consequence of the inclusion  
 $\{f_\ast (x)\  |\ x\in ( f_\ast \leq a)\}\subset (J_\ast\leq a)$.
 
\noindent
(2):\quad  By (1),  $f_\ast f^!(a)\leq a$, and hence $f^!(a)$ is in   $f_\ast\leq a$.
The element  $f^!(a)$ is the global maximum of   $f_\ast\leq a$
since it belongs to it and it is its coproduct. 
 
 \noindent
 (3):\quad  
 The relation $ f^{!}\leq f^{!}f_\ast f^{!}$ holds in general and does not require $f$ to be a homomorphism.  The reverse relation
 $f^{!}f_\ast f^{!}\leq f^{!}$ follows from  (1) and the fact that $f^{!}$ is a functor (see~\ref{koikik}).
 
 \noindent
 (4):\quad 
 The relation $\text{id}_{I_\ast}\leq f^{!}f_\ast $ and $f_\ast$ being a functor imply $f_\ast\leq f_\ast  f^{!}f_\ast$.
 The reverse relation $f_\ast f^{!}f_\ast   \leq f_\ast$ is a particular case of   (1). 
 
 \noindent
 (5):\quad 
 The inclusion $(f_\ast \leq f_\ast f^!(a))\subset (f_\ast \leq a)$
 follows from (1). If $x$ in $I_\ast$ is such that $f_{\ast}(x)\leq a$, then, since $f^!$ and $f_\ast$ are homomorphisms, $f_\ast f^!f_{\ast}(x)\leq f_\ast f^!(a)$.
 Statement (4) implies therefore $f_{\ast}(x)\leq f_\ast f^!(a)$.
This shows  $(f_\ast \leq a)\subset (f_\ast \leq f_\ast f^!(a))$ and consequently $(f_\ast \leq a)= (f_\ast \leq f_\ast f^!(a))$, proving the first part of the statement.
According to (2), if $(f\leq a)=(f\leq b)$, then $f^!(a)$ belongs to $f\leq b$,
 and consequently $f_\ast f^!(a)\leq b$, showing the minimality  of  $f_\ast f^!(a)$. 
 
 \noindent
(6):\quad  If $(f\leq a)=(f\leq b)$, then $(f_\ast\leq a)=(f_\ast\leq b)$ and hence the  coproducts
$f^!(a)=\bigvee_I(f_\ast\leq a)$ and $f^!(b)=\bigvee_I(f_\ast\leq b)$
are equal. If $f^!(a)=f^!(b)$, then $f_\ast f^!(a)=f_\ast f^!(b)$
and consequently $(f_\ast\leq a)=(f_\ast\leq f_\ast f^!(a))$ and
$(f_\ast\leq b)=(f_\ast\leq f_\ast f^!(b))$ are also equal, which implies
$(f\leq a)=(f\leq b)$.
 
  \noindent
 (7, 8):\quad These two statements are direct consequences of (4).
\end{proof}

We can use Proposition~\ref{asfadfhgfgmn}  to interpret 
the transfer of a homomorphism  of 
  finite-type as a localization which is  a process 
  involving inverting  morphisms. 

\begin{cor}\label{asdagsdfhjjk}
Let $f\colon I\to J$  be a homomorphism  of finite-type between an upper semilattice $I$ and a poset  $J$, and 
$\mathcal{C}$  a category with an initial object. 
Then the following statements about a  functor
$F\colon J_\ast\to \mathcal{C}$ are equivalent:
\begin{itemize}
    \item for all $a$ in $J$, the morphism  $F(f_\ast f^!(a)\leq a)$ is an isomorphism ($F$ inverts morphisms of the form $f_\ast f^!(a)\leq a$);
    \item there is a  functor $G\colon I_\ast\to\mathcal{C}$ for which $F$ is isomorphic to $G f^!$.
\end{itemize}
 \end{cor}
 \begin{proof}
 Assume $F(f_\ast f^!(a)\leq a)$ is an isomorphism for all $a$ in $J$. Thus, the natural transformation  $(Ff_\ast)f^!\to F$, given by  $\{F(f_\ast f^!(a)\leq a)\}_{a\in J}$,  is  an isomorphism, and  we could take $Ff_\ast$ for $G$.

Assume $F$ is isomorphic to $Gf^!$. Since $f^!f_\ast f^!=f^!$ (see~\ref{asfadfhgfgmn}.(3)), for all $a$ in $J$, the morphism $Gf^!(f_\ast f^!(a)\leq a)$ is an isomorphism. The same is therefore true for $F(f_\ast f^!(a)\leq a)$.
 \end{proof}
 
 The transfer is convenient for  constructing certain left adjoints
 even in  cases when colimits cannot be performed.
 \begin{prop}\label{sdgdfhfh}
 Let $\mathcal{C}$ be a category with an initial object  and  $f\colon I\to J$  a  homomorphism of finite-type between 
  an upper semilattice $I$ and a poset  $J$. Then the functor  $(-)^{f^{!}}\colon
    \text{\rm Fun}_{\ast}(I_\ast,\mathcal{C})\to \text{\rm Fun}_{\ast}(J_\ast,\mathcal{C})$  is left adjoint to 
   $(-)^{f_\ast}\colon \text{\rm Fun}_{\ast}(J_\ast,\mathcal{C})\to \text{\rm Fun}_{\ast}(I_\ast,\mathcal{C})$.
 \end{prop}
 \begin{proof}
 Let  $G\colon I_\ast\to \mathcal{C}$ and $F\colon J_\ast\to \mathcal{C}$ be  unital functors. 
 We need to show  there is a bijection, functorial in $G$ and $F$, between the sets of natural transformations $\text{Nat}_{J_\ast}(Gf^{!},F)$ and $\text{Nat}_{I_\ast}(G,Ff_\ast)$. 
 For a natural transformation $\phi=\{\phi_a\colon Gf^{!}(a)\to F(a)\}_{a\in J_\ast}$
 in $\text{Fun}_{\ast}(J_\ast, \mathcal{C})$, define:
 \[\overline{\phi} :=\left\{\begin{tikzcd}[column sep = 5em]
 G(x)\ar{r}{G(x\leq f^{!}f_\ast (x))} & Gf^{!}f_\ast (x)\ar{r}{\phi_{f_\ast (x)}} &  Ff_\ast (x)
 \end{tikzcd}\right\}_{x\in I_\ast}\]
 Then $\overline{\phi} $ is a natural transformation between $G$ and $Ff_\ast$ in  $\text{Fun}_{\ast}(I_\ast, \mathcal{C})$.
 
 For a natural transformation $\psi=\{\psi_x\colon G(x)\to Ff_\ast (x)\}_{x\in I_\ast}$ 
 in $\text{Fun}_{\ast}(I_\ast, \mathcal{C})$, define:
 \[\widehat{\psi} :=\left\{\begin{tikzcd}[column sep = 5em]
 Gf^{!}(a) \ar{r}{\psi_{f^{!}(a)}}   &  Ff_\ast f^{!}(a)       \ar{r}{F(f_\ast f^{!}(a)\leq a)} & F(a)
 \end{tikzcd}\right\}_{a\in J_\ast}\]
Then $\widehat{\psi} $ is a natural transformation between $Gf^{!}$ and $F$ in  $\text{Fun}_{\ast}(J_\ast, \mathcal{C})$.

 Statements (2) and (3) of Proposition~\ref{asfadfhgfgmn} imply   $\widehat{\overline{\phi}} = \phi$ and $\overline{\widehat{\psi}} =\psi$.
 The assignments $\phi\mapsto \overline{\phi}$ and $\psi\mapsto \widehat{\psi} $ are therefore inverse bijections.
 \end{proof}

 Recall that there is a commutative square of functors where the vertical functors are isomorphism of categories (see~\ref{sdgsdgjdghkfhjk}):
 \[\begin{tikzcd}[column sep = 5em, row sep=15pt]
 \text{Fun}(J,\mathcal{C})\ar{r}{(-)^f} 
 \ar{d}[swap]{-_\ast}& \text{Fun}(I,\mathcal{C})\ar{d}{-_\ast}\\
  \text{\rm Fun}_{\ast}(J_\ast,\mathcal{C}) \ar{r}{(-)^{f_\ast}} & \text{\rm Fun}_{\ast}(I_\ast,\mathcal{C})
 \end{tikzcd}\]
 Using these vertical 
 isomorphisms, Proposition~\ref{sdgdfhfh} can be rephrased as: 
 \begin{cor}\label{afdsfgjjk}
 Let $\mathcal{C}$ be a category with an initial object and  $f\colon I\to J$  a  homomorphism of finite-type between 
  an upper semilattice $I$ and a poset  $J$.  Then 
   $(-)^f\colon \text{\rm Fun}(J,\mathcal{C})\to \text{\rm Fun}(I,\mathcal{C})$ has a left adjoint given by 
   $(-)^{f^{!}}$.
 \end{cor}

 \begin{point}\label{sdfadfhfghdfgj}
The left adjoint to  $(-)^f\colon \text{\rm Fun}(J,\mathcal{C})\to \text{\rm Fun}(I,\mathcal{C})$
  is also called the \textbf{left Kan extension} along $f$ (see~\ref{werqergdfghjg}).  
   Typically such left adjoints are constructed by via colimits in $\mathcal{C}$.
  However, according to Corollary~\ref{afdsfgjjk},  when $\mathcal{C}$ has an initial object $e$ and $f\colon I\to J$ is a   homomorphism of finite-type between 
  an upper semilattice $I$ and a poset  $J$,  the left Kan extension $f^kG\colon J\to \mathcal{C}$ of $G\colon I\to C$  along $f$ exists and can be constructed explicitly using the transfer
  (see also~\cite[Proposition 5.6]{Botnan2020ART}):
  \[f^kG(a) \text{ is isomorphic to } \begin{cases}
  e &\text{ if } (f\leq a)=  \emptyset,\\
  G(\bigvee_I(f\leq a))  &  \text{ if }   (f\leq a)\not=\emptyset.
  \end{cases}
  \]
\end{point}

The above explicit description of the left Kan extension has consequences that are generically not true for arbitrary left Kan extensions. 
For example, next result states that, under the hypothesis of this section on $f$, every category $\mathcal{C}$ with an initial object is $f$-right and left invertible even though $\mathcal{C}$ is not closed under finite colimits (confront with the discussion in~\ref{asfdfhgdhj}).

\begin{prop}\label{aSFEATGERYUU}
 Let   $f\colon I\to J$  be a  homomorphism of finite-type between 
  an upper semilattice $I$ and a poset  $J$, and   $\mathcal{C}$  a category with an initial object. 
  \begin{enumerate}
  \item If $f$ is injective, then $\mathcal{C}$ is $f$-right invertible.
  \item  If $f$ is  surjective, then $\mathcal{C}$ is $f$-left invertible.
  \item  If $\mathcal{C}$ is  closed under taking  (finite) limits, then  the left Kan extension functor  $f^k\colon \text{\rm Fun}(I,\mathcal{C})\to \text{\rm Fun}(J,\mathcal{C})$ 
  commutes with (finite) limits.
  \end{enumerate}
  \end{prop} 

  \begin{proof}
  Statement (1) follows from Proposition~\ref{asfadfhgfgmn}.(7), and (2)  from Proposition~\ref{asfadfhgfgmn}.(8).
  
  To prove (3), consider a functor
   $\phi\colon  D \to  \text{\rm Fun}(I,\mathcal{C})$
   indexed by a (finite) small  category  $D$. Since  $\mathcal{C}$ is closed  under taking  (finite) limits, then so is 
   $\text{\rm Fun}(I,\mathcal{C})$. 
   According to the above description of the left Kan extension and the fact that the limits in $\text{\rm Fun}(I,\mathcal{C})$ are constructed objectwise:
   \begin{align*} f^k(\text{lim}_D\phi)(a) & \text{ is isomorphic to }\begin{cases}
  e &\text{ if } (f\leq a)=  \emptyset\\
  (\text{lim}_D\phi)(\bigvee_I (f\leq a))  &  \text{ if }   (f\leq a)\not=\emptyset
  \end{cases}\\
 \text{lim}_D  (f^k\phi)(a) &\text{ is isomorphic to }
 \begin{cases}
  e &\text{ if } (f\leq a)=  \emptyset\\
  \text{lim}_D(\phi(\bigvee_I(f\leq a)) ) &  \text{ if }   (f\leq a)\not=\emptyset
  \end{cases}
  \end{align*}
  from which it follows that $ \text{lim}_D  (f^k\phi)$ and $  f^k(\text{lim}_D\phi)$ are isomorphic. 
  \end{proof}
 
Here is a characterisation of functors isomorphic to left Kan extensions along a homomorphism of finite-type based on Corollary~\ref{asdagsdfhjjk} and the discussion in \ref{sdfadfhfghdfgj}.
\begin{cor}\label{dsgsdghf}
 Let   $f\colon I\to J$  be a  homomorphism of finite-type between 
  an upper semilattice $I$ and a poset  $J$, and  $\mathcal{C}$  a category with an initial object $e$. 
Then the following statements about a functor
$F\colon J\to \mathcal{C}$ are equivalent:  
  \begin{itemize}
      \item for all $a$ in $J$, 
      if $(f\leq a)\not=\emptyset$, then
      $F(f(\bigvee_I(f\leq a))\leq a)$ is an isomorphism, and  if 
      $(f\leq a)=\emptyset$, then
      $F(a) $ is isomorphic  to $e$;
      \item there is a functor $G\colon I\to\mathcal{C}$ for which
      $F$ is isomorphic to $f^k G$;
      \item $F$ is isomorphic to $f^k(Ff)$;
      \item there is a unital functor $G\colon I_\ast\to\mathcal{C}$
      for which $F$ is isomorphic to the composition
      $J\xhookrightarrow[]{}J_\ast\xrightarrow[]{f^!}I_\ast \xrightarrow[]{G}\mathcal{C}$.
  \end{itemize}
\end{cor}

\begin{point}
Recall that every function $f\colon I\to J$ leads to a functor $f\!\!\leq\colon  J\to J[f]$, where $J[f]\subset 2^I$ is the subposet  of
subsets of the form $(f\leq a)\subset I$ for $a$ in $J$, and  $f\!\!\leq$  maps $a$ to $f\leq a$ (see~\ref{dadgdfhg}).
Since $f\!\!\leq $ is a functor, its fibers 
$(f\!\!\leq)^{-1}(f\leq a)\subset J$, for $a$ in $J$, are convex
subposets of $J$ (see~\ref{dadgdfhg}).

As noted in~\ref{werqergdfghjg}, if $f$ is a 
functor of finite-type and  $\mathcal{C}$ is a category closed under finite colimits, then, for every $G\colon I\to \mathcal{C}$,
the left Kan extension $f^k G\colon J\to \mathcal{C}$ is isomorphic to a functor that factors through 
$f\!\!\leq \colon J\to J[f]$. This implies that the restrictions of $f^k G$ to the fibers $(f\!\!\leq)^{-1}(f\leq a)\subset J$ of $f\!\!\leq$, for $a$ in $J$, are isomorphic to  constant functors.
In general these fibers, although being convex, can still be rather complex subposets
of $J$.

Assume  $f\colon I\to J$  is  a  homomorphism of finite-type between  an upper semilattice $I$ and a poset  $J$. In this case, the only assumption we need to make about $\mathcal{C}$ is that it has an initial object $e$.
Under these assumptions the fibers of $f\!\!\leq$ can be described using the transfer. Both functors,
the transfer $f^!\colon J_\ast\to I_\ast$
and $f\!\!\leq\colon J\to J[f]$, fit into the commutative diagram:
\[  \begin{tikzcd}[row sep=15pt]
 I\ar{r}{f} & J\ar{r}{f\leq}\ar[hook]{d} & J[f]\ar{d}{V}\\
&  J_\ast\ar{r}{f^!} & I_\ast
  \end{tikzcd}\]
  where $V(f\leq a):= f^!(a)=\bigvee_{I_\ast}(f_{\ast}\leq a)$. According to Proposition~\ref{asfadfhgfgmn}.(6), $V\colon J[f]\to I_\ast$ is a monomorphism. Since
  $f\!\!\leq$ is a surjection, $V$ is a poset isomorphism between $J[f]$ and the image $f^!(J)\subset I_\ast$.
  Moreover, the fibers of  $f\!\!\leq$ can be expressed using the fibers of the transfer:
\[(f\!\!\leq)^{-1}(I\leq a)=
\begin{cases}
 (f^!)^{-1}(-\infty)\cap J & \text{ if } (I\leq a)=\emptyset,\\
 (f^!)^{-1}(f^!(a)) & \text{ if } (I\leq a)\not=\emptyset.
\end{cases}\]  
Thus, if  $(I\leq a)\not=\emptyset$, then, according to Proposition~\ref{asfadfhgfgmn}.(5),   the fiber
 $(f\!\!\leq)^{-1}(I\leq a)=(f^!)^{-1}(f^!(a))$
 is a poset with  a  global minimum given by $f(\bigvee_I(f\leq a))$.
 In this case, for every $G\colon I\to \mathcal{C}$, the left Kan extension $f^kG\colon J\to \mathcal{C}$ is  constant on 
  subposets of  $J$ that are not only  convex but also have  global minima. 
\end{point}

\begin{point}\label{wertgf}
Assume $I$ is a consistent upper semilattice  of finite-type.  According to Proposition~\ref{sdfgfds}, its realisation $\mathcal{R}(I)$ is also an upper semilattice.  Choose an 
 element $d$ in $I$ and a finite subset $V\subset (-1,0)$. According to Corollary~\ref{asfhgfgkkl}, the inclusion  $\alpha\colon\mathcal{R}_{I\leq d}(I, V)\subset \mathcal{R}(I)$ is a sublattice 
 (see~\ref{afsgsdfhdfgjhdhgj}). It consists of these pairs $(a,f)$ for which $a\leq d$, and the non-zero values of $f$ are in  $V$.

We are going to describe the  transfer $\alpha^{!}\colon \mathcal{R}(I)_{\ast}\to \mathcal{R}_{I\leq d}(I, V)_{\ast}$.
Let $(a,f)$ be in $\mathcal{R}(I)$.  By definition (see~\ref{sDGSDFHFGJ}):
\[\alpha^{!}(a,f)=\begin{cases}
 -\infty & \text{ if } \mathcal{R}_{I\leq d}(I, V)\leq (a,f)\text{ is empty,}\\
  \bigvee \left( \mathcal{R}_{I\leq d}(I, V))\leq (a,f)\right) & \text{ otherwise. }
\end{cases}\]
Thus, to describe  $\alpha^{!}(a,f)$,  the set 
 $\mathcal{R}_{I\leq d}(I, V)\leq (a,f)$  needs to be discussed.   
  We claim that it is non-empty if and only if $d$ and  $\text{supp}(f)$ have a common ancestor. Here is a proof.
 If  there is 
 $(b,g)$ in $\mathcal{R}_{I\leq d}(I, V)$ for which $(b,g)\leq (a,f)$, then,
 according to Proposition~\ref{asdgsfghg}, any ancestor of $\text{supp}(g)$ is also
 an ancestor of  $\text{supp}(f)$. The relation $b\leq d$
 implies therefore that any common ancestor
 of $b$ and  $\text{supp}(g)$ is  also a common ancestor of $d$ and $\text{supp}(f)$.
 
 Assume $d$ and $\text{supp}(f)$ have a common ancestor. 
 If $V$ is empty, set $v_0:=0$, if $V$ is non-empty, set
 $v_0$ to be the minimal element in $V$.  Define 
 $S:=\{x\in \mathcal{P}(a)\ |\ f(x)< v_0\}$.  Since $S$ is a subset of $\text{supp}(f)$, it has a common ancestor with $d$, and consequently
 the following definition of an element $c$ in $I$ makes sense as the necessary products exist due to $I$ being an upper semilattice of finite-type:
 \[c :=\begin{cases}
(\bigwedge S)\wedge d & \text{if }S\neq \emptyset,\\
a\wedge d & \text{if }S=\emptyset.
\end{cases}
\]
Any  ancestor of  $\{d,\text{supp}(f),a\}$ is also an ancestor of $c$.
Let $y$ be a parent of $c$. Consider the subset $(y\leq \mathcal{P}(a))\subset \mathcal{P}(a)$.
There are two possibilities, either $c$ is an ancestor of the set $y\leq \mathcal{P}(a)$ or not.
In the second case, for every $x$ in $(y\leq \mathcal{P}(a))\setminus (c\leq I)$, two things happen:
first, $f(x)\geq v_0$ ($x$ does not belong to $S$) and, second, 
there is an equality $y = c\wedge x$. We use the first inequality to
define a function $h\colon\mathcal{P}(c)\to (-1,0] $ by  the following formula: 
\[
h(y) :=\begin{cases}
0 & \hspace{-37mm} \text{if $c$ is an ancestor of } y\le \mathcal{P}(a),\\
\text{max}\{v\in V\ |\ v\leq  f\left((y\leq\mathcal{P}(a))\setminus  (c\le I)\right) \} & \text{otherwise. }
\end{cases}
\]
According to the above formula, 
if $h(y)<0$, then there is $x$  in
$\mathcal{P}(a)$  for which $y=c\wedge x$ and  $f(x)<0$. Any  common ancestor of $c$ and $\text{supp}(f)$  is therefore also an ancestor of $y$. As this happens for every $y$ for which 
$h(y)<0$,  
the support $\text{supp}(h)$ has an ancestor, and hence  $(c,h)$ belongs to
the realisation $\mathcal{R}(I)$.
Since $h$ has values  in $V,$ the element  $(c,h)$ belongs  to  $\mathcal{R}_{I\leq d}(I, V)$. 

Consider the translation $T_{c\leq a}h$. Let $x$ be a parent of $a$ for which $f(x)<0$.
Consistency of $I$ and the fact  $c$ and $\text{supp}(f)$ have a common ancestor, imply that  the product $c\wedge x$ exists and is either $c$ or  is  a parent of $c$.
If $c\leq x$,
then $T_{c\leq a}h(x)=-1\leq f(x)$, and, if $c\wedge x$ is a parent of $c$, then
$T_{c\leq c}h(x)= h(c\wedge x)\leq f(x)$ by the formula defining $h$.  These relations imply
$(c,h)\leq (a,f)$, and hence
the set $\mathcal{R}_{I\leq d}(I, V)\leq (a,f)$  is non-empty. 

We are going to prove  $(c,h)= \bigvee \left( \mathcal{R}_{I\leq d}(I, V))\leq (a,f)\right)$,
which gives:
$\alpha^{!}(a,f)=(c,h)$.
Let $(b,g)$ be in $\mathcal{R}_{I\leq d}(I, V)\leq (a,f)$. We need to show $(b,g)\leq (c,h)$. The relation  $b\leq c$ holds since $b$ is an ancestor of $\{d, \text{supp}(f), a\}$.
 Consider 
$T_{b\leq c}g$.  Let $y$ be a parent of $c$ for which $h(y)<0$.
Then there is $x$ in $\mathcal{P}(a)$ such that $y=c\wedge x$ and
$f(x)<0$. By the assumption $T_{b\leq c}g(y)\leq f(x)$. Since the values of 
$T_{b\leq c}g$ are in $V$, then $T_{b\leq c}g(y)\leq h(y)$, which gives $(b,g)\leq (c,h)$.

Consider  $V=\{-0.5\}$ and $d=(2,2)$. Recall that the subposet inclusion $\alpha\colon \mathcal{R}_{\mathbb{N}^2\le d}(\mathbb{N}^2,V)\hookrightarrow \mathcal{R}(\mathbb{N}^2)$ can be identified with $0.5[2]^2\subset [0,\infty)^2$ (see~\ref{sdgdfhkjkl}).
Figure~\ref{transfer} illustrates the effect of  the transfer of
$\alpha$ via this identification.

\begin{figure}
    \centering
    \includegraphics[width=6cm]{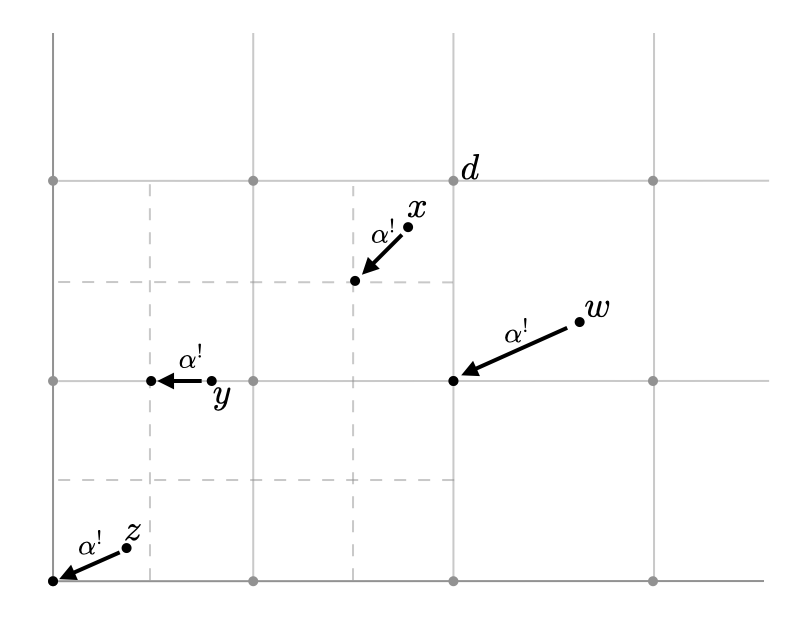}
    \caption{The intersections of the solid lines correspond to points in $\mathbb{N}^2$. The intersections of all lines correspond to points in $\mathcal{R}_{\mathbb{N}^2\le (2,2)}(\mathbb{N}^2,\{-0.5,0\})$.}
    \label{transfer}
\end{figure}
\end{point}

\section{Tame functors}\label{aoceavojcnigt}
    \begin{Def}\label{sdfubbeoru}
    Let $\mathcal{C}$ be a category and  $J$  a poset. 
  A functor $F\colon J\to \mathcal{C}$ is called \textbf{tame} or \textbf{discretisable} 
  if there is a  finite   poset  $I$ and  functors $f\colon I\to J$  and $G\colon I\to\mathcal{C}$ for which $F$ is isomorphic to the left Kan extension $f^kG$ of $G$ along $f$ (see~\ref{werqergdfghjg}), in which case $F$ is also said to be \textbf{discretised}  by $f$, and $G$ is called a \textbf{discretisation} of $F$. 
  
  The symbol $\text{Tame}(J, \mathcal{C})$ denotes the full subcategory of $\text{Fun}(J, \mathcal{C})$
  with tame functors as objects.
  \end{Def}
  
  If $J$ is a finite poset, then any functor $F\colon J\to \mathcal{C}$ is tame. 
  The notion of tameness is therefore restrictive and meaningful only when $J$ is
   infinite.  
   To determine if $F\colon J\to \mathcal{C}$ is tame  we need to look for two things:  a functor  $f\colon I\to J$ that discretises $F$ and a discretisation $G\colon I\to \mathcal{C}$ of $F$. 
   
   How can we search for a discretising $f$? In many situations it is enough to consider only finite subposet inclusions $I\subset J$. This is the case if for example $\mathcal{C}$ is closed under finite colimits (for instance  the category of vector spaces over a given field). The closure of $\mathcal{C}$ under finite colimits  guarantees the existence of the left Kan extension $f^k\colon \text{Fun}(I,\mathcal{C})\to \text{Fun}(J,\mathcal{C})$ for every functor  $f\colon I\to J$ with finite $I$. This has   the following consequences:

\begin{prop} \label{asdsfgadfshfgjh}
Let $\mathcal{C}$ be a category closed under finite colimits and $J$ a poset.
\begin{enumerate}
\item  Consider  the following  commutative diagram of poset functors with $I_0$ and $I_1$  finite: 
\[\begin{tikzcd} 
I_0\ar{r}{g}\ar[bend right]{rr}[description]{f_0}\ar{r}{g} & 
I_1\ar{r}{f_1} & J
\end{tikzcd}\]
  If  $F\colon J\to\mathcal{C}$ is discretised by $f_0$, then it is also discretised by $f_1$. 
\item  If $F\colon J\to\mathcal{C}$ is discretised by $f\colon I\to J$, then
it is also discretised by every subposet $I'\subset J$ for which   $f(I)\subset I'$.
\item   Let   $F_1,\ldots, F_n\colon J\to  \mathcal{C}$  be a finite sequence of tame functors. 
Then there is a finite subposet $I\subset J$ that  discretises $F_i$, for all $i$.
\item If $F\colon J\to \mathcal{C}$ is discretised by a subposet inclusion $f\colon I\subset J$ with $I$ finite, then $F$ is $f$-left invertible (the morphism $f^k(Ff)\to F$, adjoint to $\text{\rm id}_{Ff}$,  is an isomorphism, see~\ref{asfdfhgdhj}).
\item Let $L$ be a finite category and $\phi\colon  L \to  \text{\rm Fun}(J,\mathcal{C})$  a functor. Assume
$\phi(l)$ is tame for every object $l$ in $L$.  Then the functor  $\text{\rm colim}_L\phi$ is also tame.
\end{enumerate}
\end{prop}
\begin{proof}
\noindent (1)\quad  As long as they exist, Kan extensions commute with compositions:
$f_0^k$ and $f_1^kg^k$ are naturally isomorphic.  Thus, if $F$ is isomorphic to
$f_0^kG$ (is discretised by $f_0$), then it is also isomorphic to $f_1^k(g^kG)$ (is discretised by $f_1$). 
\smallskip

\noindent (2)\quad It  is a direct consequence of (1).

\smallskip

\noindent (3)\quad Let $I_i\subset J$ be a finite subposet discretising $F_i$. Then,
according to (2),
$\cup_{i=1}^{n}I_i\subset I$  discretises $F_i$, for every $i$.
\smallskip

\noindent (4)\quad
Consider $G\colon I\to \mathcal{C}$ such that there is an isomorphism $\alpha\colon f^k G\to F$. Since the category $C$ is closed under finite colimits, it is $f$-right invertible, and therefore  so is the functor  $G$. Consequently,  $f^kG$ is $f$-left invertible (see~\ref{asfdfhgdhj}).
\smallskip
\noindent (5)\quad 
Use (3) to choose a finite subposet inclusion $f\colon I\subset J$  that discretises $\phi(l)$ for every  $l$ in $L$.
By (4), the natural transformation $f^k (\phi f)\to \phi$, adjoint to the identity $\text{id}\colon \phi f\to \phi f$, is an isomorphism.
Kan extensions commute with colimits and hence $\text{\rm colim}_L\phi$ is isomorphic to  $f^k\text{\rm colim}_L (\phi f)$, proving its tameness. 
\end{proof}
  
 According to Proposition~\ref{asdsfgadfshfgjh}, if $\mathcal{C}$ is closed under finite colimits, the  search for discretising functors can be 
  restricted to finite subposet inclusions $I\subset J$, in which case
  a discretisation is given by restricting to $I$ (see~\ref{werqergdfghjg}). 
  That is a considerable simplification.
  For example, a tame functor indexed by a realisation $\mathcal{R}(I)$ with values in a category closed under finite colimits can always be discretised by a subposet of the form $\mathcal{R}_D(I,V)$, for some finite  $D\subset I$ and a finite  $V\subset (-1,0)$ (see~\ref{sdgdfhkjkl}).
 However not all the categories $\mathcal{C}$ for which we would like  to have a simpler way of verifying  tameness 
are closed under finite colimits.  For example, the category of simplicial complexes fails to have this property.
In such a case,  Corollary~\ref{dsgsdghf} can be used where the only assumption on $\mathcal{C}$ is that it has an initial object $e$. It gives a characterisation of functors which are discretised by a homomorphism $f\colon I\to J$ out of a finite upper semilattice $I$. According to this corollary, a functor 
$F\colon J\to\mathcal{C}$ is discretised by 
$f$ if and only if, for every $a$ in $J$, if $(f\leq a)\not=\emptyset$, then $F(f(\bigvee_I(f\leq a))\leq a)$ is an isomorphism, and if 
$(f\leq a)=\emptyset$, then $F(a) $ is isomorphic  to $e$.
 For an arbitrary $J$, not all tame functors are discretised by homomorphisms.  To make sure all tame functors are discretised by a homomorphism,  
$J$  itself needs to be  an upper  semilattice.
The rest  of this section is devoted to discussing tameness under this assumption on $J$ and when 
 $\mathcal{C}$ is a category with  an initial object $e$. 

\begin{prop}\label{sdgfgjhdgh}
Let $J$ be an upper semilattice and $\mathcal{C}$ a category with an initial object.
\begin{enumerate}
\item 
A functor  $F\colon J\to \mathcal{C}$ is tame if and only if it is
discretised by a finite sublattice $I\subset J$ (see~\ref{afsgsdfhdfgjhdhgj}).
\item Let $I_0\subset I_1\subset J$ be finite sublattices of $J$. If $F\colon J\to \mathcal{C}$  is discretised by $I_0\subset J$, then it is also discretised by $I_1\subset J$.
\item    Let   $F_1,\ldots, F_n\colon J\to  \mathcal{C}$  be a finite sequence of tame functors. Then there is a finite sublattice $I\subset J$ 
   that  discretises $F_i$, for all $i$.
\end{enumerate}
\end{prop}
\begin{proof}
\noindent
(1):\quad 
Let $F\colon J\to\mathcal{C}$ be tame. Choose a finite poset $I'$ and  functors $f\colon I'\to J$ and $G\colon I'\to \mathcal{C}$ for which there is an isomorphism $f^kG\to F$. 
 Let $\phi\colon G\to Ff$ be the natural transformation adjoint to this isomorphism.
 Since $f(I')\subset J$ is finite, the sublattice generated by this image $I:=\langle f(I')\rangle\subset J$ is also finite (see~\ref{afsgsdfhdfgjhdhgj}).
Let us denote the inclusion $I\subset J$ by $g$ and by $h\colon I'\to I$ the functor that maps $x$ to $f(x)$. 
Since $g$ is a homomorphism of finite-type between upper semilattices, the left Kan extension $g^k\colon \text{Fun}(I,\mathcal{C})\to \text{Fun}(J,\mathcal{C})$
 exists (see Section~\ref{asfsdfhfhgjmd}). 
In particular, $g^k(Fg)$ exists and, hence, by Proposition~\ref{asfadfhgfgmn}.(4), the natural transformation $Fg\to (g^k(Fg))g$ is an isomorphism, meaning that $Fg$ is $g$-right invertible. 
By Propostion~\ref{dfvbsfdgbsfgb}, $F$ is $g$-left invertible, implying that $F$ is discretised by the finite sublattice inclusion $g\colon I\subset J$.
 \smallskip
 
 \noindent
 (2):\quad  Exactly the same argument as in Proposition~\ref{asdsfgadfshfgjh}.(2) can be used to prove this statement since
 the left Kan extensions along $I_0\subset I_1$, $I_0\subset J$, and $I_1\subset J$ exist by Corollary~\ref{afdsfgjjk}. 
 \smallskip
 
 \noindent
 (3):\quad  
 For every $i$ choose a finite sublattice  $ I_i\subset J$   discretising  $F_i$. 
 Define $I:=\langle I_1\cup\cdots \cup I_n\rangle\subset J$ (see~\ref{afsgsdfhdfgjhdhgj}).
Since    $I$ is finite, according to (2),   $F_i$ is   discretised by the inclusion $I\subset J$ for every $i$.
 \end{proof}
 
For example,  let $J$ be an upper semilattice of finite-type.
Then for every finite subset $I\subset J$, there is $b$ (for example $\bigvee_JI$) for which $I\subset (J\leq b)\subset J$.  Consequently, according to Proposition~\ref{sdgfgjhdgh},
 a functor $F\colon J\to \mathcal{C}$ is tame if and only if 
 it is discretised by the sublattice $(J\leq b)\subset J$, for some $b$ in $J$. Thus, according to Corollary~\ref{dsgsdghf}, 
   $F$ is tame if and only if, there is $b$ in $J$ such that, for every $a$ in $J$, if $(J\leq b)\cap(J\leq a)\not=\emptyset$, then
   $F\left((\bigvee (J\leq b)\cap(J\leq a))\leq a\right)$ is an isomorphism, and if $(J\leq b)\cap(J\leq a)=\emptyset$, then $F(a)$ is isomorphic to $e$.
\vspace{5mm}
  
  If $ \mathcal{C}$ is closed under taking finite colimits and limits, then so is the category 
  of all functors $\text{Fun}(J,\mathcal{C})$ for every poset $J$.  
  In this case, the category $\text{Tame}(J,\mathcal{C})$ is also closed under finite colimits (see Proposition~\ref{asdsfgadfshfgjh}.(4)).  In general, we do not know if 
  $\text{Tame}(J,\mathcal{C})$ is closed under  finite limits. 
   A sufficient assumption to guarantee this is, again, that $J$ is an upper semilattice.

  \begin{prop}\label{afsadagdfhbg}
Let $J$ be an upper semilattice,  $ \mathcal{C}$  a category   closed under finite limits, $L$ a finite category, and $\phi\colon  L \to  \text{\rm Fun}(J,\mathcal{C})$  a functor. Assume
$\phi(l)$ is tame for every object $l$ in $L$.  Then the functor  $\text{\rm lim}_L\phi$ is also tame.
\end{prop}
\begin{proof}
  Let $f\colon I\subset J$ be a finite  sublattice that discretises $\phi(l)$ for every object $l$ in $L$
  (see Proposition~\ref{sdgfgjhdgh}.(3)).  Then the natural transformation $f^k (\phi f)\to \phi $, adjoint to the identity $\text{id}\colon \phi f\to \phi f$, is an isomorphism (see Corollary~\ref{dsgsdghf}). Consequently, $\text{\rm lim}_L\phi$ is isomorphic to $\text{\rm lim}_Lf^k (\phi f)$ and hence, 
  by Proposition~\ref{aSFEATGERYUU}, it is also isomorphic to $f^k (\text{\rm lim}_L\phi f)$. 
   \end{proof}

 According to Propositions~\ref{asdsfgadfshfgjh}.(4) and~\ref{afsadagdfhbg},
 if  $J$ is an upper semilattice and $ \mathcal{C}$ is  closed under finite limits and colimits,  then so is
the category of tame functors  $\text{Tame}(J, \mathcal{C})$. 
Furthermore, the inclusion   $\text{Tame}(J, \mathcal{C})\subset \text{Fun}(J, \mathcal{C})$ preserves finite limits and colimits.

 \section{Homotopy theory of tame functors}
 \label{afasdfdghfdgjh}
Let $\mathcal{M}$ be a model category as defined in~\cite{MR1361887}. This means  three classes of morphisms in $\mathcal{M}$ are chosen:
{\bf weak equivalences} ($\xrightarrow{\sim}$), {\bf fibrations} ($\twoheadrightarrow$), and {\bf cofibrations}
($\hookrightarrow$), which
are required to satisfy  the following axioms:
\begin{enumerate}
	\item[MC1] Finite limits and colimits exist in $\mathcal{M}$.
	\item[MC2] If $f$ and $g$ are morphisms in $\mathcal{M}$ for which $g f$ is defined and if two of the three morphisms $f$, $g$, $gf$ are weak equivalences, then so is the third.
	\item[MC3]   The three classes of morphisms are
	preserved by retracts.
	\item[MC4] Consider a  commutative square 
	in $\mathcal M$ consisting of the solid morphisms:
	\[
	\begin{tikzcd}[row sep=13pt]
	X \arrow[hook']{d}[swap]{\alpha} \arrow[r] 
	& E\arrow[two heads]{d}{\beta}
	\\
	Y\arrow[r] \arrow[ru,dotted]
	& B
	\end{tikzcd}
	\]
	Then a morphism, depicted by the dotted arrow,   making this diagram commutative, exists under either of the following two assumptions:  (i) $\alpha$ is a cofibration and a weak equivalence and $\beta$  is a fibration, or
	(ii) $\alpha$ is a cofibration and $\beta$ is a fibration and a weak equivalence.
	\item[MC5] Every morphism in $\mathcal{M}$ can be factored in two ways: 
	(i) $ \beta\alpha$, where $\alpha$ is a cofibration and $\beta$ is  a fibration and a weak equivalence, and 
	(ii) $\beta\alpha$, where $\alpha$ is  a cofibration and a weak equivalence and $\beta$ is a fibration.
\end{enumerate}

In particular, the axiom MC1 guarantees the existence of the initial object $\text{colim}_{\emptyset}F$, denoted by $e$, and of the terminal object $\text{lim}_{\emptyset}F$, denoted by $\ast$. 
An object $X$ in $\mathcal M$ is called {\bf cofibrant} if the morphism $e \to X$ is a cofibration. 
If the morphism $X\to \ast$ is a fibration, then $X$ is called {\bf fibrant}.

Assume  $I$ is a finite  poset. Then the following choices of weak equivalences, fibrations, and cofibrations in 
$\text{Fun}(I,\mathcal{M})$ gives a model structure (see for example~\cite{MR1361887,reedy}).
A natural transformation $\varphi:F \to G$ in $\text{Fun}(I,\mathcal{M})$  is
\begin{itemize}
 \item a weak equivalence (resp.\@  a fibration) if, for all $a$ in $I$, the morphism $\varphi_a: F(a) \to G(a)$ is a weak equivalence (resp.\@ a fibration) in
 $\mathcal{M}$;
 \item a cofibration if, for all $a$ in $I$,
 the morphism
 \[\colim\left(\colim_{I<a}G \leftarrow \colim_{I<a} F \rightarrow F(a)\right)\to G(a),\] 
 induced by the following commutative diagram, is a cofibration in $\mathcal{M}$:
\[\begin{tikzcd}
\colim_{I<a}F \ar{r} \ar{d}[swap]{\colim_{I<a} \varphi} & F(a) \ar{d}{\varphi_a}\\
\colim_{I<a}G \ar{r} & G(a)
\end{tikzcd}\]
\end{itemize}

The described model structures  on functors indexed by finite posets  are compatible in the following sense. 
Let  $f\colon I_0\to I_1$   be a functor between finite posets. If $\phi\colon F\to G$ is a cofibration (resp.\@ a cofibration and a weak equivalence)
in $\text{\rm Fun}(I_0,\mathcal{M})$, then so is $f^k\phi\colon f^kF\to f^kG$ 
in  $\text{\rm Fun}(I_1,\mathcal{M})$.
This  is a consequence of the universal property of  left Kan extensions and the axiom (MC4).
However,  left Kan extensions in general  fail to preserve weak equivalences and fibrations (compare with Proposition~\ref{asdgadfsghjdhjghk}).

To verify if  a natural transformation
in $\text{Fun}(I,\mathcal{M})$ is a cofibration, we need to perform colimits over
subposets $(I<a)\subset I$. 
In general, these subposets can be large and constructing  colimits over them
may require performing a lot of identifications. When $I$ is an upper semilattice however,
we can be  more efficient.  For $a$ in a finite  upper semilattice  $I$, define a sublattice:
\[I_a:=\{\text{$\bigwedge$} S\ |\ S\subset \mathcal{P}(a) \text{ has an ancestor}\}\subset (I<a).\]
For all   $b$  in $I<a$, the set $b\leq \mathcal{P}(a)$ is non-empty and $b\leq \bigwedge (b\leq \mathcal{P}(a))$. 
The element $ \bigwedge  (b\leq \mathcal{P}(a))$ is therefore  the  initial  object in  $b\leq I_a$ and 
consequently $b\leq I_a$  is contractible (see~\cite{MR0365573}). As this happens for all $b$ in $I<a$,  the inclusion $I_a\subset  (I<a)$ is cofinal (see~\cite{MR1712872, MR0365573}), which proves:

\begin{prop}\label{dsdfgdafhrgyhj}
Let $I$ be a finite upper semilattice and $a$  its element.   For every  functor $F\colon (I<a)\to\mathcal{M}$, the morphism $\text{\rm colim}_{I_a}F\to \text{\rm colim}_{I<a}F$ is an isomorphism and 
$\text{\rm hocolim}_{I_a}F\to \text{\rm hocolim}_{I<a}F$ is a weak equivalence.
\end{prop}

According to Proposition~\ref{dsdfgdafhrgyhj}, when $I$ is a finite upper semilattice,  to verify if a natural transformation in $\text{Fun}(I,\mathcal{M})$ is a cofibration, we need only to perform colimits over the subposets  $I_a\subset I$ for $a$ in $I$.  How can such colimits be calculated? One way is to consider, for $a$ in $I$, the subposet of the the discrete cube
(see~\ref{afadfhfgsh}):
\[C_a:=\{S\subset \mathcal{P}(a)\  |\  \mathcal{P}(a)\setminus S\text{ has an  ancestor}\}\subset 2^{\mathcal{P}(a)}\]
and the functor $\bigwedge^c\colon C_a\to I_a$ mapping $S$ to $\bigwedge (\mathcal{P}(a)\setminus S)$. Since, for every
$b$ in $I_a$,  the subset  $\mathcal{P}(a)\setminus(b\leq \mathcal{P}(a))\subset  \mathcal{P}(a)$ is the initial object in the poset $b\leq \bigwedge^c$ (see~\ref{sdgsdgjdghkfhjk}), this poset is contractible. The  functor 
$\bigwedge^c$ is  therefore cofinal. Thus, for every  $F\colon  I_a\to \mathcal{M}$, the morphism 
$\text{colim}_{C_a}(F\bigwedge^c) \to \text{colim}_{I_a}F$ is an isomorphism and
$\text{hocolim}_{C_a}(F\bigwedge^c) \to \text{hocolim}_{I_a}F$ is a weak equivalence.

The assumption on the indexing poset  being an upper semilattice is not only  helpful in
 verifying if a natural tranformation is  a cofibration. It is also  crucial  in proving that 
the following choices of weak equivalences, fibrations, and  cofibrations in $\text{Tame}(J,\mathcal{M})$  satisfy the axioms of a model structure.

\begin{Def}\label{sdfsdgsfghjd}
Let  $J$ be an upper semilattice and $\mathcal{M}$ a model category. A natural  transformation $\phi\colon F\to G$  in $\text{Tame}(J,\mathcal{M})$ is called:
\begin{itemize}
\item a weak equivalence (resp. fibration) if, for all $a$ in $J$,  the morphism $\phi_a\colon F(a)\to G(a)$ is  
a weak equivalence (resp. fibration) in $\mathcal{M}$;
\item  a cofibration if there is a finite subposet inclusion
$f\colon I\subset J$ discretising  $F$ and $G$, and  for which   $\phi^f\colon Ff\to Gf$ is a cofibration in
$\text{Fun}(I,\mathcal{M})$.
\end{itemize}
\end{Def}

\begin{thm}\label{afadfhsfgh}
Let  $J$ be an upper semilattice and $\mathcal{M}$ a model category. The choices described in Definition~\ref{sdfsdgsfghjd} satisfy the axioms of a model structure on $\text{\rm Tame}(J,\mathcal{M})$.
\end{thm}

Before we prove Theorem~\ref{afadfhsfgh}, we first show:

\begin{prop}\label{asdgadfsghjdhjghk}
Let $\mathcal{M}$ be a model category and $f\colon I\to J$  a homomorphism of finite-type
 between  upper semilattices.
If a natural transformation $\phi\colon F\to G$  in  $\text{\rm Fun}(I,\mathcal{M})$ is a 
weak equivalence  (resp.\@ fibration), then so is its left Kan extension $f^k\phi\colon
f^kF\to f^kG$ in $\text{\rm Fun}(J,\mathcal{M})$.
\end{prop}
\begin{proof}
Let $\phi\colon F\to G$ in  $\text{\rm Fun}(I,\mathcal{M})$  be a weak equivalence  (resp.\@ fibration).
Since $f$ is a finite-type homomorphism between upper semilattices,
 the morphism  $(f^k\phi)_a\colon
f^kF(a)\to f^kG(a)$  is isomorphic to either $\text{id}\colon e\to e$, if ($f\leq a)=  \emptyset$,
or to $\phi_{\bigvee_{I}(f\leq a)}\colon  F(\bigvee_{I}(f\leq a)) \to  G(\bigvee_{I}(f\leq a)) $, if
$ (f\leq a)\not=\emptyset$ (see~\ref{sdfadfhfghdfgj}). As these  morphisms are  weak equivalences  (resp.\@ fibrations), for all $a$, then so is the left Kan extension $f^k\phi$.
\end{proof}
\begin{proof}[Proof of Theorem~\ref{afadfhsfgh}]
Requirement MC1 follows from  Proposition~\ref{afsadagdfhbg}.  
MC2 and MC3 are clear since $\mathcal{M}$ satisfies them.
\smallskip

\noindent
MC4:\quad
 Let $\alpha\colon F\hookrightarrow G$ be a cofibration in $\text{\rm Tame}(J,\mathcal{M})$.
Choose a finite sublattice $I \subset J$ and a functor $f: I \to J$ that discretises both $F$ and $G$ and for which
$\alpha^f\colon Ff\to  Gf$
is a cofibration in $\text{\rm Fun}(I,\mathcal{M})$. 
Assume  $\alpha$ is part of a commutative square in $\text{\rm Tame}(J,\mathcal{M})$ 
depicted in~(\ref{adfhsgghjs})  on the left, where $\beta$ is a fibration. 
By applying $(-)^f$ to this square we get a commutative square in $\text{\rm Fun}(I,\mathcal{M})$, depicted  by the solid arrows square on the right in~(\ref{adfhsgghjs}),  where the natural transformations are a cofibration and a fibration in $\text{\rm Fun}(I,\mathcal{M})$. 
If, in addition,  $\alpha$ or $\beta$ is  a weak equivalence, then the lift, depicted by the dotted arrow on  the right of~(\ref{adfhsgghjs}), exists
since $\text{\rm Fun}(I,\mathcal{M})$ is a model category.
\begin{equation}
	\begin{tikzcd}[column sep = 40pt, row sep = 18pt]
	F \arrow[hook']{d}[swap]{\alpha} \arrow[r] 
	& E\arrow[two heads]{d}{\beta}
	\\
	G\arrow[r] 
	& B
	\end{tikzcd}
	\ \ \  \ \ \ \ \ \  \ \ \ \ \ \ 
	\begin{tikzcd}[column sep = 40pt, row sep = 18pt]
	Ff \arrow[hook']{d}[swap]{\alpha^f} \arrow[r] 
	& Ef\arrow[two heads]{d}{\beta^f}
	\\
	Gf\arrow[r]  \arrow[dotted]{ur}[description]{s}
	& Bf
	\end{tikzcd}
\label{adfhsgghjs}
\end{equation}
By applying the left Kan extension $f^k$ to the square on the right of~(\ref{adfhsgghjs})
and comparing the result to the original square on the left, we can form the following commutative diagram in $\text{\rm Tame}(J,\mathcal{M})$:
\begin{equation*}
	\begin{tikzcd}
	& F\ar{rrr} \arrow[hook']{dd}[swap, pos=0.2]{\alpha}& & & E\arrow[two heads]{dd}{\beta}\\
	f^k(Ff)\ar[crossing over]{rrr} \ar{ur} \arrow[hook']{dd}[swap]{(f^k\alpha)^f }& & &  f^k(Ef)\ar{ur}\\
	& G\ar{rrr} & & & B \\
	f^k(Gf)\ar{rrr} \ar{ur} \arrow[dotted,bend right=15,crossing over]{uurrr}[swap, description]{f^ks} & & & f^k(Bf)\ar{ur} \arrow[from=uu, two heads,crossing over, description, near start, "(f^k\beta)^f "]	
	\end{tikzcd}
\label{cdfgdfagdfgdfvons}
\end{equation*}
Since the natural transformation $f^k(Gf)\to G$ is an isomorphism, a desired lift exists
in the right square of ~(\ref{adfhsgghjs}) under the additional assumption that either $\alpha$ or $\beta$ is a weak equivalence.
\smallskip

\noindent
MC5:\quad
Let $\phi\colon F\to G$ be a natural transformation in $\text{\rm Tame}(J,\mathcal{M})$.
Choose a finite sublattice $f\colon I\subset J$ that discretises both $F$ and $G$. 
Factor the natural transformation 
$\phi^f\colon Ff\to Gf$  in $\text{\rm Fun}(I,\mathcal{M})$
as  $\phi^f=\beta\alpha$  where $\alpha\colon Ff\to H$ is a cofibration, $\beta\colon H\to Gf$ is a fibration, and either $\alpha$ or $\beta$ is a weak equivalence.  
By applying  the left Kan extension $f^k$  to these factorisations we obtain a commutative diagram
in $\text{\rm Tame}(J,\mathcal{M})$ where the vertical natural transformations are isomorphisms:
\[\begin{tikzcd}[row sep = small]
f^k(Ff)\ar{d}\ar{r}{f^k\alpha} & f^kH\ar{r}{f^k\beta} & f^k(Gf)\ar{d}\\
F\ar{rr}{\phi} & & G
\end{tikzcd}\] 
Since $(f^k\alpha)^f$ is isomorphic to $\alpha$, it is a cofibration 
in $\text{\rm Fun}(I,\mathcal{M})$ and consequently $f^k\alpha$ is a cofibration in 
$\text{\rm Tame}(J,\mathcal{M})$.  According to Proposition~\ref{asdgadfsghjdhjghk},
$f^k\beta$ is a fibration  in $\text{\rm Tame}(J,\mathcal{M})$.  The same proposition assures also that 
if $\alpha$ (resp. $\beta$) is a weak equivalence, then so is $f^k\alpha$
(resp. $f^k\beta$). This gives the desired factorisation of $\phi$.
\end{proof}

\section{Betti diagrams of vector spaces valued functors.}\label{sdgadfhsgdb}

In this section, we show that   upper semilattices and realisations share similar
 homological algebra properties.
For example, they both allow to reduce the computation of Betti diagrams in the minimal resolution of a tame functor to convenient finite subsets, and this can be done using Koszul complexes.
As a consequence, we show that, in both cases, the length of a minimal resolution of a functor is  bounded above by the parental-dimension of elements of a subposet of the indexing poset, providing a Hilbert syzygies-type theorem.

\begin{point}[\em Freeness]\label{dhjfgki}
Let $(J, \leq)$ be a poset and $\text{vect}_K$ the category of finite dimensional $K$-vector spaces. 
Consider the poset $(J,=)$ with the trivial poset  relation on the set $J$,
where two elements are related if and only if they are equal. 
The identity function $\iota\colon (J,=)\to J$, mapping $a$  to $\iota(a)=a$, 
is a functor.
A functor $V\colon (J,=)\to \text{vect}_K$ is just a sequence $\{V_a\}_{a\in J}$ of finite dimensional $K$-vector spaces. 
The set $\{a\in J\ |\ V_a\not = 0\}$ is called the support of $V$ and is denoted by $\text{supp}(V)$.

A functor $F\colon J\to \text{vect}_K$ is called  \textbf{free finitely generated} if it is isomorphic to the left Kan extension along $\iota\colon  (J,=)\to J$ of some  $V=\{V_a\}_{a\in J}$ whose support is finite. 
Since  all free functors considered in this article are  finitely generated, we  use the term \textbf{free}, without mentioning finite generation, to describe such functors.  
The name free is justified by the  universal property of the left Kan extension, which gives a linear isomorphism between $\text{Nat}_J(\iota^k V, H)$ and $\prod_{a\in J}\text{Hom}(V_a,H(a))$,
for every functor $H\colon J\to \text{vect}_K$. For example, for $b$ in $J$ and a finite dimensional vector space $U$, consider the simple functor 
$U[b]\colon  J\to \text{vect}_K$ defined as
\[
U[b](a):=\begin{cases}
0 &\text{ if } a\not = b,\\
U &\text{ if } a = b.
\end{cases}
\]
Then the vector spaces  $\text{Nat}_J(\iota^k V, U[b])$ and
$\text{Hom}(V_b,U)$ are isomorphic.
By varing  $b$ in $J$ and taking $U$ to be $K$, we can conclude that  free functors $\iota^k V$ and $\iota^k W$ are isomorphic  if and only if 
$V_a$ and $W_a$ are isomorphic for all $a$ in $J$. Thus, a free  functor $F$ determines a unique sequence $\beta F=\{(\beta F)_a\}_{a\in J}$ of vector spaces whose support is finite, called the \textbf{Betti diagram} of $F$, for which $F$ is isomorphic to the left Kan extension $\iota^k (\beta F)$.  
If $\text{supp}(\beta F)=\{a\}$, then
$F$ is called \textbf{homogeneous} and is also denoted by the symbol 
$F(a)[a,-)$. The restriction of  $F(a)[a,-)$
to the subposet $(a\leq J)\subset J$ is isomorphic to the constant functor with value $F(a)$. 
Its restriction to  $\{x\in J\ |\ a\not \leq x\}$
is isomorphic to the constant functor with value $0$.
If $F$ is free, then it is isomorphic to the direct sum
$\oplus_{a\in J}(\beta F)_a[a,-)$  and, consequently,  $F$ is isomorphic to the left Kan extension of $\{(\beta F)_a\}_{a\in \text{supp}(\beta F)}$ along $(\text{supp}(\beta F),=)\hookrightarrow J$. 
Thus, every free functor $F$ is tame and is  discretised by the subposet inclusion $\text{supp}(\beta F)\subset J$.
Moreover, every collection $V=\{V_a\}_{a\in J}$ with finite support is the Betti diagram of a free functor.

If $G\colon I\to \text{vect}_K$ is free, then, for every functor
$f\colon I\to J$, the left Kan extension $f^kG\colon J\to \text{vect}_K$ is also free and $(\beta f^k G)_a=0$ if $f^{-1}(a)$ is empty, and   $(\beta f^k G)_a$ is isomorphic to $\oplus_{b\in f^{-1}(a)}(\beta G)_b$, if $f^{-1}(a)$ is non-empty.
\end{point}
\begin{point}[\em Resolutions]
Let $J$ be a poset and $F\colon J\to \text{vect}_K$   a functor.
An exact sequence $P_{n-1}\to\cdots\to P_0\to F\to 0$ in $\text{Fun}(J,\text{vect}_K)$, with $P_i$ free for every $i$, is called an \textbf{$n$-resolution} of $F$. 
A $1$-resolution $P_0\to F\to 0$ is also called a \textbf{cover} of $F$.
A $2$-resolution $P_1\to P_0\to F\to 0$ is also called a \textbf{presentation} of $F$. 
An infinite exact sequence  $\cdots\to P_1\to P_0\to F\to 0$,   with $P_i$  free for every $i$, is called an  \textbf{$\infty$-resolution} of $F$.  
An $n$-resolution of $F$ is   also denoted, as a map of chain complexes, by $P\to F$, where $F$ is a chain complex concentrated in degree $0$, and either 
 $P=(P_{n-1}\to\cdots\to P_0)$ or 
 $P=(\cdots\to P_1\to P_0)$, depending if $n$ is finite.

A functor $F\colon J\to \text{vect}_K$ is called \textbf{$n$-resolvable} if it has an $n$-resolution.  Functors that are $1$-resolvable are also called \textbf{finitely generated}.  
Functors that are $2$-resolvable are also called \textbf{finitely presented}.
If the indexing poset $J$ is finite, then all functors are $\infty$-resolvable. 
This is because all such functors are finitely generated, and hence by  taking   covers  of successive  kernels, an $\infty$-resolution can be constructed.

If $J$ is infinite, then not all $n$-resolvable functors have to be $(n+1)$-resolvable.
For example, let $J=[0,2)\coprod \{a,b\}\coprod (2,3]$, with $a$ and $b$ incomparable and
$x<a<y$ and $x<b<y$, for $x$ in $[0,2)$ and $y$ in $(2,3]$. 
Consider the subposet $I= \{1, a,b\}\subset J$ and a functor $G\colon I\to \text{\rm Vect}_K$ where both 
$G(1<a)$ and $G(1<b)$ are given by  $K\to 0$.
Let $F\colon J\to  \text{\rm Vect}_K$ be the left Kan extension of $G$ along $I\subset J$.
Then $F$ is finitely presented ($2$-resolvable), however it is not 
$3$-resolvable.
\end{point}

Here is a characterisation of finitely generated and finitely presented functors.
 \begin{prop}\label{sgdgjhf}
 Let $J$ be a poset and  $F\colon J\to \text{\rm vect}_K$ a functor.
 \begin{enumerate}
     \item  $F$  is finitely generated if and only if 
     there is a finite poset $I$ and a functor  $f\colon I\to J$ for which the natural transformation
 $\mu\colon f^k(Ff)\to F$, adjoint to the identity $\text{\rm id}\colon Ff\to Ff$, is surjective.
 \item $F$  is finitely presented if and only if it is tame, i.e., if and only if there is a finite poset $I$ and a functor  $f\colon I\to J$ such that $F$ is $f$-left invertible (for which the natural transformation
 $\mu\colon f^k(Ff)\to F$, adjoint to the identity $\text{\rm id}\colon Ff\to Ff$, is an isomorphism).
 \end{enumerate}

\end{prop}
\begin{proof}
\noindent
(1):\quad Let $\pi\colon P_0\to F$ be a  cover. Since $P_0$ is tame, there is
a functor  $f\colon I\to J$, from a finite $I$, discretising $P_0$ (for example $I=\text{supp}(\beta P_0)\subset J$). 
The commutativity of the following square and the fact that left Kan extensions preserve surjections imply the surjectivity of $\mu$:
\[\begin{tikzcd}
f^k(P_0 f)\ar{r}\ar{d}[swap]{f^k (\pi^f)} & P_0\ar{d}{\pi}\\
f^k(Ff)\ar{r}{\mu} & F
\end{tikzcd}\]

Let $I$ be a finite poset,  $f\colon I\to J $  a functor  for which $\mu\colon f^k(Ff)\to F$ is  surjective, and 
 $\pi\colon P_0\to Ff$  a cover.
The surjectivity of $\mu$  implies the surjectivity of the 
 following composition, which is then a cover of $F$:
\[\begin{tikzcd}
f^k(P_0)\ar{r}{f^k\pi}& f^k(Ff)\ar{r}{\mu} & F
\end{tikzcd}\]
\smallskip

\noindent(2):\quad 
If $P_1\to P_0\to F\to 0$ is a $2$-resolution, then $\text{colim}(0\leftarrow P_1\to P_0)$ is  isomorphic to $F$. 
Since tameness is preserved by finite colimits (see Proposition~\ref{asdsfgadfshfgjh}.(4)), $F$ is  tame. 

Assume $F$ is discretised by a subposet inclusion $f\colon I\subset J$ with finite $I$. Consider a $2$-resolution
 $P_1\to P_0\to Ff\to 0$, which exists since $I$ is finite.
As before,  $f^k(Ff)$, and hence $F$, is  isomorphic to $\text{colim}(0\leftarrow f^kP_1\to f^kP_0)$. Since   $f^kP_1$ and $f^kP_0$
are free, $F$ is finitely presented.
\end{proof}

Proposition~\ref{sgdgjhf} characterises  $n$-resolvable functors for $n\leq 2$. 
We do not have a similar characterisation for $n>2$, only a  partial result:

\begin{prop}\label{fhjghgkyik}
Let $J$ be a poset. Assume $F\colon J\to \text{\rm vect}_K$ is discretised by 
a functor $f\colon I\to J$ from a finite poset $I$ and for which
$f^k\colon \text{\rm Fun}(I,\text{\rm vect}_K)\to \text{\rm Fun}(J,\text{\rm vect}_K)$  is exact.  Then $F$ is $\infty$-resolvable.
\end{prop}
\begin{proof}
Let $G\colon I\to \text{\rm vect}_K$ be such that $f^kG$ is isomorphic to $F$. Choose an $\infty$-resolution $\pi\colon  P\to G$. Exactness of $f^k$ means that the left Kan extension $f^k\pi \colon f^kP\to f^kG$ is an $\infty$-resolution 
$f^kG$. 
\end{proof}

Here is  a condition guaranteeing exactness of the left Kan extension:
\begin{lemma}\label{sjjhkiluilk}
Assume $f\colon I\to J$ is a functor of posets from a finite poset $I$ satisfying the following property:
for every $a$ in $J$, every pair of elements in $f\leq a$ that have an ancestor also has 
a descendent in $f\le a$ (such posets are called weakly directed in~\cite{beyondpersistence}).
Then  $f^k\colon \text{\rm Fun}(I,\text{\rm vect}_K)\to \text{\rm Fun}(J,\text{\rm vect}_K)$  is exact.
\end{lemma}
\begin{proof}
We  show   $\text{colim}_{f\leq a} \colon
\text{\rm Fun}(f\leq a,\text{\rm vect}_K)\to \text{\rm vect}_K$ is  exact for all $a$ in $J$. 
Let $M\subset (f\leq a)$ consists of all the maximal elements.
The assumption on $f$ implies  $M\subset (f\leq a)$ is cofinal and hence  $\text{colim}_{f\leq a} F$ is isomorphic to
$\bigoplus_{x\in M} F(x)$. The lemma follows from the exactness of  direct sums.
\end{proof}

\begin{cor}\label{asdfgadfgsdfhg}
Let $J$ be a poset.
\begin{enumerate}
    \item Every functor  $F\colon J\to \text{\rm vect}_K$
     discretised by a homomorphism $f\colon I\to J$ from a finite upper semilattice  $I$ is  $\infty$-resolvable.
    \item If $J$ is an upper semilattice, then  
    a functor  $F\colon J\to \text{\rm vect}_K$ is tame if and only if it is  $\infty$-resolvable.
\end{enumerate}
\end{cor}
\begin{proof}
In both (1) and (2) the poset $f\leq a$, for all $a$ in $J$, has a terminal object given by its  coproduct in $I$. Thus, by Lemma~\ref{sjjhkiluilk}, the conclusion of Proposition~\ref{fhjghgkyik} holds. 
\end{proof}

\begin{point}[\em Minimality and Betti diagrams]
An $n$-resolution $\pi\colon P\to F$ is called \textbf{minimal} if   every chain map $\phi\colon P\to P,$ for which the following triangle commutes, is an isomorphism:
\[\begin{tikzcd}[row sep = small]
P\ar{rr}{\phi}\ar{dr} & & P\ar{dl}\\
& F
\end{tikzcd}\]
If $P\to F$ and $Q\to F$ are minimal $n$-resolutions of $F$, then $P$ and $Q$ are isomorphic.
Thus, if  $P\to F$ is a minimal $n$-resolution, then the isomorphism type of  $P_i$, for $i< n$, is  uniquely determined by the isomorphism type of $F$, and its  Betti diagram
$\beta P_i$ is called the $i$-th \textbf{Betti diagram} of $F$, 
and is 
denoted by 
$\beta^{i} F$.  
For example if $F$ is free, then $\text{id}\colon F\to F$ is a minimal $\infty$-resolution of $F$ and, hence, for every $a$ in $J$, the vector space
$(\beta^{0} F)_a$ is isomorphic to $(\beta F)_a$, and $(\beta^{i} F)_a=0$ for $i>0$.

A minimal $1$-resolution $P_0\to F$ is also called a \textbf{minimal cover} of $F$.
An $n$-resolution $P\to F$ is  minimal if and only if, for every $i<n$, the following 
natural transformations are minimal   covers:
 $P_0\to F$, $P_1\to \text{ker}(P_{0}\to F)$, \ldots,
$P_{i}\to \text{ker}(P_{i-1}\to P_{i-2})$. Thus, if
$P_0\to F$ is a minimal cover of an $n$-resolvable functor, then, for $1\leq i<n$, $\beta^iF$ is isomorphic to $\beta^{i-1}\text{ker}(P_0\to F)$.
Also, if $P\to F$ is a minimal $n$-resolution, then
$\beta^0F$ is isomorphic to $\beta P_0$,
$\beta^1F$ is isomorphic to $\beta^0\text{ker}(P_{0}\to F)$, and for $i>1$, 
$\beta^iF$ is isomorphic to $\beta^0\text{ker}(P_{i-1}\to P_{i-2})$.
\end{point}

Our strategy for constructing minimal resolutions of $n$-resolvable functors 
indexed by an arbitrary poset is to reduce this problem to
 the case when the  indexing poset is finite.
 
\begin{prop}\label{afsgsdfshfd}
Let $J$ be a poset and $F\colon J\to  \text{\rm vect}_K $ a functor. 
Assume $P\to F$  is a $n$-resolution  which is  discretised 
by a subposet inclusion  $f\colon I\subset J$, with  $I$ finite. 
 If $Q\to Ff$ is a minimal $n$-resolution of $Ff$, then its adjoint $f^kQ\to F$ is a minimal resolution of $F$ and, for $i<n$ and $a$ in $J$,
\[(\beta^i F)_a \text{ is isomorphic to }\begin{cases}
    (\beta^i (Ff))_a= (\beta Q_i)_a & \text{ if } a\in I,\\
    0 & \text{ if } a\not \in I.
    \end{cases}\]
\end{prop}
\begin{proof}
By restricting the $n$-resolution $P\to F$ of $F$ along $f$, we obtain a resolution $Pf\to Ff$ of $Ff$. Let   $Q\to Ff$ be a minimal $n$-resolution. Then there are chain maps 
$\phi\colon Q\to Pf$ and $\psi\colon Pf\to Q$ for which the composition
$\psi\phi\colon Q\to Q$ is an isomorphism and the  diagram on the left commutes:
\[\begin{tikzcd}[row sep = small]
Q\ar{r}{\phi} \ar{dr} & Pf\ar{d}\ar{r}{\psi} & Q\ar{dl}\\
& Ff
\end{tikzcd}\ \ \ \ \ \ 
\begin{tikzcd}[row sep = small]
f^kQ\ar{r}{f^k\phi} \ar{dr} & f^k(Pf)\ar{d}\ar{r}{f^k\psi} & f^kQ\ar{dl}\\
& F
\end{tikzcd}\]
By taking the adjoints of the vertical maps in the left diagram, we obtain a commutative diagram
on the right. Since  $P=f^k(Pf)\to F$ is an $n$-resolution, then so is its retract $f^kQ\to F$. Its minimality follows from the  minimality of $Q\to Ff$ and the fact that  $f^k\colon \text{Nat}_I(Q_i,Q_i)\to \text{Nat}_J(f^kQ_i,f^kQ_i)$ is a
 bijection for every $i$ (see~\ref{asfdfhgdhj}).
\end{proof}

According to Proposition~\ref{afsgsdfshfd}, 
to  construct a minimal $n$-resolution of a functor  $F\colon J\to \text{\rm vect}_K$,
the first step  is to find  a finite subposet inclusion $f\colon I\subset J$ for which there is an  $n$-resolution $P\to F$ with $P$ being discretised by $f$ (the natural transformation $f^k(Pf)\to P$ is an isomorphism). 
For $n=1$, according to the proof of Proposition~\ref{sgdgjhf}.(1),  such a subposet inclusion is given by any $f\colon I\subset J$ for which $f^k(Ff)\to F$ is a surjection. For $n=2$, according to the proof of Proposition~\ref{sgdgjhf}.(2), such a subposet inclusion is given by any $f\colon I\subset J$ that discretises $F$. 
We do not have a similar statement for $n> 2$. 

The second step is to construct a minimal $n$-resolution $Q\to Ff$ of the restriction $Ff$. The adjoint  of this minimal resolution $f^kQ\to F$ is then the  desired minimal $n$-resolution of $F$. This process reduces finding a minimal $n$-resolution  of $F$ to finding a minimal $n$-resolution of $Ff$ which is a functor indexed by a finite poset. Constructing minimal resolutions of functors indexed by finite posets is standard and involves radicals (see for example~\cite{MR1476671, MR1731415}).

\begin{point}[\em Radical]
Let $I$ be a finite  poset.
The \textbf{radical} of $G\colon I\to \text{\rm vect}_K$  is a subfunctor $\text{rad}(G)\subset G$ given by  
$\text{rad}(G)(a)=\text{im} \left(\bigoplus_{s<a}G(s)\to G(a)\right)$ for  $a$ in $I$.
 The quotient functor $G/\text{rad}(G)$ is semisimple as it is   isomorphic to a  direct sum $\oplus_{a\in I} U_a[a]$ of simple functors (see~\ref{dhjfgki}), where $U_a:=(G/\text{rad}(G))(a)$. 
 For example, for a  free functor $G=\oplus _{a\in I} (\beta G)_a[a,-)$, the quotient
 $G/\text{rad}(G)$ is isomorphic to  $\oplus _{a\in I} (\beta G)_a[a]$.

 A key property of quotienting by the radical, when the indexing poset is finite, is the  surjectivity detection: a natural transformation $H\to G$ is surjective if and only if its composition with the quotient $G\to G/\text{rad}(G)$
 is surjective. 
 The surjectivity detection  can be used to construct minimal  covers. 
Consider the quotient $G/\text{rad}(G)$ of    $G\colon I\to \text{\rm vect}_K$. Set  $P_0:=\bigoplus_{a\in I}U_a[a,-)$, where
$U_a= (G/\text{rad}(G))(a)$. Note that 
there is an isomorphism $P_0/\text{rad}(P_0)\to G/\text{rad}(G)$. Let $\pi\colon P_0\to G$ be any natural transformation fitting into the following commutative square, where
the bottom horizontal arrow represents the chosen isomorphism:
\[\begin{tikzcd}[row sep = small]
P_0\ar{r}{\pi}\ar{d} & G\ar{d}\\
P_0/\text{rad}(P_0)\ar{r} & G/\text{rad}(G)
\end{tikzcd}\]
Such  $\pi$ exists since $P_0$ is free. The composition $P_0\to G/\text{rad}(G)$ is surjective, and hence so is $\pi$.  The same argument 
gives the surjectivity of every $\phi\colon P_0\to P_0 $ for which $\pi \phi = \pi$.
Since the values of $P_0$ are finite dimensional, every such  $\phi$ is therefore an isomorphism, and hence $\pi\colon P_0\to G $ is a minimal  cover.  This shows:
\end{point}
\begin{prop}\label{sadgfsfasddh}
Let $I$ be a finite poset and $G\colon I\to \text{\rm vect}_K$ a functor. A natural transformation
$P_0\to G$ is a minimal cover if and only if $P_0$ is free and the induced natural transformation
$P_0/\text{\rm rad}(P_0)\to G/\text{\rm rad}(G)$ is an isomorphism. Furthermore, $(\beta^0G)_a$
is isomorphic to $(G/\text{\rm rad}(G))(a)$ for all $a$.
\end{prop}

Since all the functors indexed by a finite poset $I$ are finitely generated, by taking minimal  covers of  successive kernels, every functor indexed by $I$ admits a minimal $\infty$-resolution. 
This inductive  step-wise construction of a minimal $\infty$-resolution can be used for a step-wise
inductive procedure of calculating the Betti diagrams. However, can the Betti diagrams be retrieved
in one step without the need of an inductive procedure? 
The answer to this question is: yes, it can standardly be done by using Koszul complexes.

\begin{point}
Let $I$ be a finite  poset and $a$ its element.  Choose a linear ordering $\prec$ on the set  of  parents $\mathcal{P}(a)$ of $a$. For every functor $G\colon I\to \text{vect}_K$, we  define a non-negative chain complex denoted by $\mathcal{K}_a G$ and called the \textbf{Koszul complex} of $G$ at $a$.
Let $k$ be a natural number. Define:
\[(\mathcal{K}_a G)_k:=\begin{cases}
G(a) & \text{ if } k=0,\\
\displaystyle  \bigoplus _{\substack{S\subset \mathcal{P}(a),\ |S|=k\\ S\text{ has an ancestor}}} \text{colim}_{\substack{\cap_{s\in S}(I\le s)}} G
& \text{ if } k>0.
\end{cases}
\]
For example, $(\mathcal{K}_a G)_1=\bigoplus _{s\in \mathcal{P}(a)}G(s) $ and
$(\mathcal{K}_a G)_k=0$ if $k>\text{par-dim}_I(a)$ (see Proposition~\ref{sfgasg}).

Define $\partial\colon (\mathcal{K}_a G)_{k+1}\to (\mathcal{K}_a G)_k$ as follows:
\begin{itemize}
\item If $k=0$, then $\partial\colon (\mathcal{K}_a G)_1\to (\mathcal{K}_a G)_0=G(a)$ is the linear function which on 
the summand  $G(s)$ in $(\mathcal{K}_a G)_1$, indexed by $s$ in $\mathcal{P}(a)$, is given by $G(s<a)$.
\item Let $k>0$. For $k\geq j\geq 0$, let 
$\partial_j\colon  (\mathcal{K}_a G)_{k+1}\to (\mathcal{K}_a G)_k$ 
  be the linear function mapping
the summand $\text{colim}_{\substack{\cap_{s\in S}(I\le s)}} G$  in $(\mathcal{K}_a G)_{k+1}$, indexed by
$S=\{s_0\prec\cdots\prec s_{k}\}\subset \mathcal{P}(a)$, to the summand
$\text{colim}_{\substack{\cap_{s\in S\setminus\{s_j\}}(I\le s)}} G$  in $(\mathcal{K}_a G)_{k}$,
indexed by $S\setminus\{s_j\}\subset \mathcal{P}(a)$, via the function
of the colimits $\text{colim}_{\substack{\cap_{s\in S}(I\le s)}} G\to \text{colim}_{\substack{\cap_{s\in S\setminus\{s_j\}}(I\le s)}} G$ induced by the poset inclusion
$\cap_{s\in S}(I\le s)\subset \cap_{s\in S\setminus\{s_j\}}(I\le s)$.
Define $\partial\colon (\mathcal{K}_a G)_{k+1}\to (\mathcal{K}_a G)_k$  to be the alternating sum
$\partial = \sum_{j=0}^k(-1)^j\partial_j$.
\end{itemize}
The linear functions $\partial$ form a chain complex as it is standard to verify that composition of two consecutive such functions is the $0$ function.

For a natural transformation 
$\phi\colon F\to G$, define:
\[
\begin{tikzcd}
(\mathcal{K}_a F)_k\ar{r}{(\mathcal{K}_a \phi)_k} & (\mathcal{K}_a G)_k
\end{tikzcd}:=
\begin{cases}
\phi_a & \text{ if } k=0,\\
\displaystyle \bigoplus _{\substack{S\subset \mathcal{P}(a),\ |S|=k\\ S\text{ has an ancestor}}} \text{colim}_{\substack{\cap_{s\in S}(I\le s)}} \phi
& \text{ if } k>0.
\end{cases}
\]
These  linear  functions, for all $k$, form the chain map $\mathcal{K}_a \phi\colon \mathcal{K}_a F\to \mathcal{K}_a G$.
The association $\phi\mapsto \mathcal{K}_a \phi$ is a functor. 

We observe that the image of the  differential $\partial\colon (\mathcal{K}_aG)_1\to (\mathcal{K}_aG)_1=G(a)$
coincides with  $\text{rad}(G)(a)$.
Consequently, the vector spaces $H_0(\mathcal{K}_aG)$ and $(G/\text{rad}(G))(a)$ are isomorphic, and, according to Proposition~\ref{sadgfsfasddh}, they are also isomorphic to
$(\beta^0 G)_a$.
\end{point}
\begin{point}\label{dgdfgdfhsfgn}
Consider an element $a$ in a finite poset $I$ satifying: 
every  $S\subset \mathcal{P}(a)$ having an ancestor also has the product  $\bigwedge S$ in $I$. For example, if $I$ is an upper semilattice, then all its elements satisfy this property. Under this assumption, for every subset
$S\subset \mathcal{P}(a)$ that has an ancestor, the product $\bigwedge S$ is the  terminal object in the category 
$\cap_{s\in S}(I\le s)$ and consequently, for $k>0$,
\[
(\mathcal{K}_a G)_k= \displaystyle  \bigoplus _{\substack{S\subset \mathcal{P}(a),\ |S|=k\\ S\text{ has an ancestor}}}  G(\bigwedge S).
\]
\end{point}

\begin{point}
Let $I$ be a finite poset and $a$ its element.
Since colimits commute with direct sums, so does the functor $\mathcal{K}_a$, i.e., 
the natural transformation $\mathcal{K}_a F\oplus \mathcal{K}_a G\to \mathcal{K}_a(F\oplus G)$ is an isomorphism.

As the colimit operation is right exact, then so is our  Koszul complex construction: if 
$0\to F\to G\to H\to 0$ is an exact sequence in $\text{Fun}(I,\text{vect}_K)$, then 
$\mathcal{K}_aF\to \mathcal{K}_aG\to \mathcal{K}_aH\to 0$ is an exact sequence of chain complexes.
In general  colimits  do not preserve monomorphisms, and hence one does not expect  the  Koszul complex construction to preserve monomorphisms in general either.  Let $n$ be an extended positive natural number (containing $\infty$). An element $a$ in $I$ is called \textbf{Koszul $n$-exact} if, for every exact sequence of functors $0\to F\to G\to H\to 0$ and  $k< n$, the  following sequence of vector spaces
is exact:
\[
0\to (\mathcal{K}_aF)_k\to (\mathcal{K}_aG)_k\to (\mathcal{K}_aH)_k\to 0
\]
For example all elements in $I$ turn out to be   Koszul $2$-exact. 
\end{point}
\begin{prop}\label{asfgdfhjdfh}
Let $I$ be a finite poset. 
\begin{enumerate}
    \item Then every element  $a$ in $I$ is  Koszul $2$-exact.
    \item Assume an element $a$ in $I$ has the following property: every subset $S\subset \mathcal{P}(a)$ which has an ancestor has the product  $\bigwedge S$ in $I$. Then $a$ is Koszul $\infty$-exact.
    \item If $I$ is an upper semilattice, then all its elements are Koszul $\infty$-exact.
\end{enumerate}
\end{prop}
\begin{proof}
\noindent
Statement (1)  is a consequence of the fact that  direct sums preserves exactness.
Under the assumption of statement (2),  for a functor $F$  and $k>0$, the vector space
$ (\mathcal{K}_aF)_k$ can also be described as (see~\ref{dgdfgdfhsfgn}):
\[
\displaystyle  \bigoplus _{\substack{S\subset \mathcal{P}(a),\ |S|=k\\ S\text{ has an ancestor}}}  F(\bigwedge S).
\]
In this case, the statement also follows from the  exactness of direct sums. Finally statement (3) is a particular case of (2).
\end{proof}

Proposition~\ref{asfgdfhjdfh}  translates into the homology exact sequence:
\begin{cor}\label{aSDFGDFHJG}
Let $I$ be a finite poset and $0\to F\to G\to H\to 0$ be an exact sequence 
in $\text{\rm Fun}(I,\text{\rm vect}_K)$.
\begin{enumerate}
    \item  For every element $a$ in $I$, there is an exact sequence of vector spaces:
    \[\hspace{-10mm}\begin{tikzcd}
     & H_2(\mathcal{K}_a G)\ar{r} & H_{2}(\mathcal{K}_a H)\ar[out=0, in=180]{dll}\\
 H_1(\mathcal{K}_a F)\ar{r} & H_1(\mathcal{K}_a G)\ar{r}& H_{1}(\mathcal{K}_a H)\ar[out=0, in=180]{dll}\\
  H_0(\mathcal{K}_a F)\ar{r} & H_0(\mathcal{K}_a G)\ar{r}& H_{0}(\mathcal{K}_a H)\ar{r} & 0
    \end{tikzcd}\]
    \item Assume an element $a$ in $I$ has the following property: every subset $S\subset \mathcal{P}(a)$ which has an ancestor has the product  $\bigwedge S$ in $I$. Then 
    there is an exact sequence of vector spaces:
    \[\hspace{-10mm}\begin{tikzcd}
    & & \cdots \ar[out=0, in=180]{dll}\\
         H_2(\mathcal{K}_a F)\ar{r} & H_2(\mathcal{K}_a G)\ar{r}& H_{2}(\mathcal{K}_a H)\ar[out=0, in=180]{dll}\\
     H_1(\mathcal{K}_a F)\ar{r} & H_1(\mathcal{K}_a G)\ar{r}& H_{1}(\mathcal{K}_a H)\ar[out=0, in=180]{dll}\\
    H_0(\mathcal{K}_a F)\ar{r} & H_0(\mathcal{K}_a G)\ar{r}& H_{0}(\mathcal{K}_a H)\ar{r} & 0
    \end{tikzcd}\]
\end{enumerate}
\end{cor}

\begin{point}\label{cuhnoinjwcaeodj}
Consider a homogeneous free functor
$F=V[b,-)\colon I\to  \text{vect}_K$, where $I$ is a finite poset.
We claim that $\mathcal{K}_a F$ has the following homology:
\[H_i(\mathcal{K}_a F)\text{ is isomorphic to }
\begin{cases}
V & \text{ if } i= 0 \text{ and } b= a, \\
0 &\text{ otherwise. }
\end{cases}\]
To prove the claim, recall that $F(x)=0$ if $b\not\leq x$, and  $F$ restricted to
$(b\leq I)\subset I$ is isomorphic to the constant functor with value $V$. 
Thus, if $b\not\leq a$, then $\mathcal{K}_a F=0$, and the claim holds.
If $b=a$, then $(\mathcal{K}_a F)_0=F(a)=V$ and $(\mathcal{K}_a F)_i=0$ for $i>0$, and again the claim holds. Assume $b<a$. Then $(\mathcal{K}_a F)_0=F(a)$, which is isomorphic to $V$. Moreover, for a  subset that has an ancestor 
$S\subset \mathcal{P}(a)$ with $|S|>0$, the colimit $\text{colim}_{\substack{\cap_{s\in S}(I\le s)}} F $ is isomorphic to $\text{colim}_{\substack{\cap_{s\in S}(b\leq I\le s)}} F $,
 which is either isomorphic to $V$, in the case $S\subset (b\leq \mathcal{P}(a))$, or is $0$ otherwise. Consequently,  the complex
 $\mathcal{K}_a F$  is  isomorphic to $L\otimes V$, where $L$
 is  the augmented  chain complex of  the standard
$|b\leq \mathcal{P}(a)|$-dimensional simplex whose homology is trivial in all degrees:
\[L:=\left(\cdots\to\bigoplus _{\substack{S\subset (b\leq \mathcal{P}(a))\\  |S|=2}}K \to\bigoplus _{\substack{S\subset (b\leq \mathcal{P}(a))\\  |S|=1}}K \to  K\right).
\]
Since Koszul complexes commute with direct sums, if $F$ is free, isomorphic to
$\oplus_{b\in I}(\beta F)_b[b,-)$, then:
\[H_i(\mathcal{K}_a F)\text{ is isomorphic to }
\begin{cases}
(\beta F)_a& \text{ if } i=0,\\
0 &\text{ if } i>0.
\end{cases}
\]

\end{point}

We are now ready to state the key fact connecting the homology of the Koszul complexes of a functor with its Betti diagrams.

\begin{thm}\label{asdfgdfhfgjh}
Let $I$ be a finite poset and $G\colon I\to  \text{\rm vect}_K$ a functor.
\begin{enumerate}
    \item  For every $a$ in $I$ and $i=0,1,2$, the vector spaces $(\beta^i G)_a$ and   $H_i(\mathcal{K}_a G)$ are isomorphic.
    \item Assume an element $a$ in $I$ has the following property: every subset $S\subset \mathcal{P}(a)$ which has an ancestor has the product  $\bigwedge S$ in $I$. Then, for every $i$,  the vector spaces $(\beta^i G)_a$ and   $H_i(\mathcal{K}_a G)$ are isomorphic.
\end{enumerate}
\end{thm}
\begin{proof}
The proof relies on Corollary~\ref{aSDFGDFHJG}. Since the arguments for statements (1) and (2) are analogous, we show only (1).
The case $i=0$  follows from  Proposition~\ref{sadgfsfasddh} and the fact that
$H_0(\mathcal{K}_a G)$  is isomorphic to $(G/\text{rad}(G))(a)$.

Consider an  exact sequence
$0\to S_1\to P_0\xrightarrow{\pi} G\to 0$ where $\pi$ is a minimal  cover. 
It leads to an  exact sequence of homologies (see Corollary~\ref{aSDFGDFHJG}):
\[\begin{tikzcd}
     & H_2(\mathcal{K}_a P_0)\ar{r} & H_{2}(\mathcal{K}_a G)\ar[out=0, in=180]{dll}[description]{\alpha_2}\\
 H_1(\mathcal{K}_a S_1)\ar{r} & H_1(\mathcal{K}_a P_0)\ar{r}& H_{1}(\mathcal{K}_a G)\ar[out=0, in=180]{dll}[description]{\alpha_1}\\
  H_0(\mathcal{K}_a S_1)\ar{r} & H_0(\mathcal{K}_a P_0)\ar{r}{H_0(\mathcal{K}_a\pi)}& H_{0}(\mathcal{K}_a G)\ar{r} & 0
\end{tikzcd}\]
Minimality of $\pi$ is equivalent to  $H_0(\mathcal{K}_a\pi)$ being an isomorphism.
Since $P_0$ is free,   $H_1(\mathcal{K}_a P_0)=H_2(\mathcal{K}_a P_0)=0$ (see~\ref{cuhnoinjwcaeodj}). 
These two observations imply    $H_{1}(\mathcal{K}_a G)$ is isomorphic to $H_0(\mathcal{K}_a S_1)$, and $H_{2}(\mathcal{K}_a G)$ is isomorphic to  $H_1(\mathcal{K}_a S_1)$. 
By the already proven case $i=0$, $H_0(\mathcal{K}_a S_1)$ is isomorphic to
$(\beta^0S_1)_a$ which is isomorphic to $(\beta^1G)_a$. This gives the case $i=1$. Applying this case to $S_1$, we get that $H_1(\mathcal{K}_aS_1)$ is
isomorphic to $(\beta ^1 S_1)_a$, which is isomorphic to $(\beta ^2 G)_a$, and the case $i=2$ also holds.
\end{proof}

Here are some consequences of the presented statements, 
which are proven by the same strategy: first dicretise and then use the Koszul complex construction.

\begin{cor}\label{fjgdkvfkdcmv}
Let $J$ be a poset and $F\colon J\to  \text{\rm vect}_K$ a functor.
\begin{enumerate}
    \item Assume $F$ is $n$-resolvable and $f\colon I\subset J$  is a 
    subposet inclusion  with finite $I$ discretising an $n$-resolution of $F$.
    Then $\text{\rm supp}(\beta^iF)\subset I$ for all $0\leq i< n$.  Moreover, for $0\leq i< \text{\rm min}(3,n)$ and $a$ in $I$, $(\beta^i F)_{a}$ is isomorphic to 
    $H_i(\mathcal{K}_a(Ff))$.
\item Assume $F$ is discretised by a  subposet inclusion $f\colon I\subset J$ 
which is a homomorphism out of a  finite upper semilattice $I$. Then $F$
is $\infty$-resolvable and
    $\text{\rm supp}(\beta^iF)\subset I$ for  all  $i\geq 0$. Moreover,
    for $i\geq 0$ and $a$ in $I$, $(\beta^i F)_{a}$ is isomorphic to $H_i(\mathcal{K}_a(Ff))$.
\item Assume $J$ is an upper semilattice and $F$  a tame functor. 
Then $F$ is $\infty$-resolvable and, for $a$ in $J$ with $\text{\rm par-dim}_J(a)<i$, $(\beta^i F)_{a}=0$.
\end{enumerate}
\end{cor}

\begin{proof}
Statement (1) is a consequence of Proposition~\ref{afsgsdfshfd} and  Theorem~\ref{asdfgdfhfgjh}. Statement (2) is a consequence of Corollary~\ref{asdfgadfgsdfhg}.(1), and again Proposition~\ref{afsgsdfshfd} and  Theorem~\ref{asdfgdfhfgjh}. Statement (3) follows 
from statement (2) and Proposition~\ref{aDFDFHFHJ} since
tame functors indexed by an upper semilattice can be discretised by a finite sublattice  (see Proposition~\ref{sdgfgjhdgh}.(1)).
\end{proof}

We finish this article with our key theorem describing how to determine Betti diagrams of functors indexed by realisations of finite-type  posets that admit discretisable resolutions. 
For such functors, Koszul complexes can be used, similarly to the case of functors indexed by upper semilattices (see Corollary~\ref{fjgdkvfkdcmv}).   
The main reason for this is the following fact:
functors indexed by realisations have  natural grid-like  discretisations that can be refined in a way that every set of parents of an element having an ancestor
also has a product in the refinement.

\begin{thm}\label{sdrtyhgf}
Let $I$ be a finite-type poset and $F\colon \mathcal{R}(I)\to \text{\rm vect}_K$ an $n$-resolvable functor. Assume $d$ is an element in $I$ and $V$ is a finite subset of $(-1,0)$ for which 
$\text{\rm supp}(\beta^jF)\subset \mathcal{R}_{I\le d}(I,V)\stackrel{\alpha}{\hookrightarrow} \mathcal{R}(I)$ for all $j\leq n$.
\begin{enumerate}
\item 
Let $0\leq i<\text{\rm min}(3,n)$ and $(a,f)$ be in $\text{\rm supp}(\beta^iF)$.
Then  $(\beta^i F)_{(a,f)}$ is isomorphic to $H_i(\mathcal{K}_{(a,f)}(F\alpha))$.
\item Let $i<n$ and $(a,f)$ be in $\text{\rm supp}(\beta^iF)$ for which  there is 
$\varepsilon$ in $V$ such that $f(x)>\varepsilon$ for all $x$ in $\mathcal{P}(a)$.
Then  $(\beta^i F)_{(a,f)}$ is isomorphic to $H_i(\mathcal{K}_{(a,f)}(F\alpha))$.
\item  If $\text{\rm par-dim}_{\mathcal{R}(I)}(a,f)< i< n $, then 
$(\beta^i F)_{(a,f)}=0$.
\end{enumerate}
\end{thm}
\begin{proof}
\noindent
(1)\quad The assumption implies that a minimal $\text{min}(3,n)$-resolution of $F$ is discretised by
$\mathcal{R}_{I\le d}(I,V)\subset \mathcal{R}(I)$.
This statement is then a particular case of 
Corollary~\ref{fjgdkvfkdcmv}.(1).

\noindent
(2)\quad  The assumption implies that a minimal $i+1$ resolution of $F$ is discretised by
$\mathcal{R}_{I\le d}(I,V)\subset \mathcal{R}(I)$. Moreover the  product of every set of  parents of  $(a,f)$ in $\mathcal{R}_{I\le d}(I,V)$ exists.
This statement is then a particular case of 
Theorem~\ref{asdfgdfhfgjh}.(2).

\noindent
(3)\quad 
Since $\text{\rm par-dim}_{\mathcal{R}_{I\le d}(I,V)}(a,f)\leq \text{\rm par-dim}_{\mathcal{R}(I)}(a,f)$
for every $(a,f)$ in  $\mathcal{R}_{I\le d}(I,V)$,
this statement follows from (2).
\end{proof}

\bibliographystyle{plainurl}
\bibliography{Tameness/bibliography}

\begin{thebibliography}{10}

\bibitem{MR1476671}
Maurice Auslander, Idun Reiten, and Sverre~O. Smal\o.
\newblock {\em Representation theory of {A}rtin algebras}, volume~36 of {\em
  Cambridge Studies in Advanced Mathematics}.
\newblock Cambridge University Press, Cambridge, 1997.
\newblock Corrected reprint of the 1995 original.

\bibitem{beyondpersistence}
Mattia~G. Bergomi, Massimo Ferri, Pietro Vertechi, and Lorenzo Zuffi.
\newblock Beyond topological persistence: Starting from networks.
\newblock {\em Mathematics}, 9(23), 2021.
\newblock \href {https://doi.org/10.3390/math9233079}
  {\path{doi:10.3390/math9233079}}.

\bibitem{Botnan2020ART}
Magnus~Bakke Botnan, Justin Curry, and Elizabeth Munch.
\newblock A relative theory of interleavings.
\newblock {\em ArXiv}, abs/2004.14286, 2020.

\bibitem{MR0365573}
A.~K. Bousfield and D.~M. Kan.
\newblock {\em Homotopy limits, completions and localizations}.
\newblock Lecture Notes in Mathematics, Vol. 304. Springer-Verlag, Berlin-New
  York, 1972.

\bibitem{MR4323617}
Peter Bubenik and Nikola Mili\'{c}evi\'{c}.
\newblock Homological {A}lgebra for {P}ersistence {M}odules.
\newblock {\em Found. Comput. Math.}, 21(5):1233--1278, 2021.
\newblock \href {https://doi.org/10.1007/s10208-020-09482-9}
  {\path{doi:10.1007/s10208-020-09482-9}}.

\bibitem{MR2476414}
Gunnar Carlsson.
\newblock Topology and data.
\newblock {\em Bull. Amer. Math. Soc. (N.S.)}, 46(2):255--308, 2009.
\newblock \href {https://doi.org/10.1090/S0273-0979-09-01249-X}
  {\path{doi:10.1090/S0273-0979-09-01249-X}}.

\bibitem{MR1731415}
Henri Cartan and Samuel Eilenberg.
\newblock {\em Homological algebra}.
\newblock Princeton Landmarks in Mathematics. Princeton University Press,
  Princeton, NJ, 1999.
\newblock With an appendix by David A. Buchsbaum, Reprint of the 1956 original.

\bibitem{MR4057607}
Wojciech Chach\'{o}lski and Henri Riihim\"{a}ki.
\newblock Metrics and stabilization in one parameter persistence.
\newblock {\em SIAM J. Appl. Algebra Geom.}, 4(1):69--98, 2020.
\newblock \href {https://doi.org/10.1137/19M1243932}
  {\path{doi:10.1137/19M1243932}}.

\bibitem{algebraic_stability}
F.~Chazal, D.~Cohen-Steiner, M.~Glisse, L.~J. Guibas, and Steve~Y. Oudot.
\newblock Proximity of persistence modules and their diagrams.
\newblock {\em Proceedings of the 25th annual symposium on Computational
  geometry}, SCG ’09:237–246, 2009.

\bibitem{MR3873177}
Ren\'{e} Corbet and Michael Kerber.
\newblock The representation theorem of persistence revisited and generalized.
\newblock {\em J. Appl. Comput. Topol.}, 2(1-2):1--31, 2018.
\newblock \href {https://doi.org/10.1007/s41468-018-0015-3}
  {\path{doi:10.1007/s41468-018-0015-3}}.

\bibitem{MR1361887}
W.~G. Dwyer and J.~Spali\'{n}ski.
\newblock Homotopy theories and model categories.
\newblock In {\em Handbook of algebraic topology}, pages 73--126.
  North-Holland, Amsterdam, 1995.
\newblock \href {https://doi.org/10.1016/B978-044481779-2/50003-1}
  {\path{doi:10.1016/B978-044481779-2/50003-1}}.

\bibitem{MR4334502}
Woojin Kim and Facundo M\'{e}moli.
\newblock Generalized persistence diagrams for persistence modules over posets.
\newblock {\em J. Appl. Comput. Topol.}, 5(4):533--581, 2021.
\newblock \href {https://doi.org/10.1007/s41468-021-00075-1}
  {\path{doi:10.1007/s41468-021-00075-1}}.

\bibitem{Lesnick2015}
M.~Lesnick.
\newblock The theory of the interleaving distance on multidimensional
  persistence modules.
\newblock {\em Foundations of Computational Mathematics}, 15:613–650, 2015.

\bibitem{MR1712872}
Saunders Mac~Lane.
\newblock {\em Categories for the working mathematician}, volume~5 of {\em
  Graduate Texts in Mathematics}.
\newblock Springer-Verlag, New York, second edition, 1998.

\bibitem{EMiller1}
Ezra Miller.
\newblock Essential graded algebra over polynomial rings with real exponents.
\newblock {\em arXiv:2008.03819}, 2020.

\bibitem{OliverWojtek}
Gäfvert Oliver and Chach\'olski Wojciech.
\newblock Stable invariants for multidimensional persistence.
\newblock {\em arXiv:1703.03632}, 2017.

\bibitem{MR3975559}
Amit Patel.
\newblock Generalized persistence diagrams.
\newblock {\em J. Appl. Comput. Topol.}, 1(3-4):397--419, 2018.
\newblock \href {https://doi.org/10.1007/s41468-018-0012-6}
  {\path{doi:10.1007/s41468-018-0012-6}}.

\bibitem{reedy}
Christopher Reedy.
\newblock Homotopy theory of model categories.
\newblock {\em Unpublished manuscript}, 1974.

\bibitem{MR3735858}
Martina Scolamiero, Wojciech Chach\'{o}lski, Anders Lundman, Ryan Ramanujam,
  and Sebastian \"{O}berg.
\newblock Multidimensional persistence and noise.
\newblock {\em Found. Comput. Math.}, 17(6):1367--1406, 2017.
\newblock \href {https://doi.org/10.1007/s10208-016-9323-y}
  {\path{doi:10.1007/s10208-016-9323-y}}.

\bibitem{MR347936}
Ross Street.
\newblock Two constructions on lax functors.
\newblock {\em Cahiers Topologie G\'{e}om. Diff\'{e}rentielle}, 13:217--264,
  1972.

\end{thebibliography}

\end{document}